\documentclass[12pt]{article}
\usepackage{amsmath,amssymb,epsfig}
\newcommand{\be}{\begin{equation}}
\newcommand{\ee}{\end{equation}}
\def\H{{\cal H}}
\def\B{{\cal B}}

\newcounter{unnumber}

\newtheorem{thrm}{Theorem}
\newtheorem{lm}{Lemma}

\newtheorem{prop}{Proposition}

\newtheorem{de}{Definition}
\newtheorem{prf}[unnumber]{Proof}

\numberwithin{equation}{section}

\begin{document}

\thispagestyle{empty} \phantom{}\vspace{-1.5cm}
%\rightline{HU-EP-01/31} \rightline{hep-th/0108220}
\rightline{\today} \vspace{-0.5cm}

\begin{center}
  {\bf \LARGE Representations of the
ultrahyperbolic BMS group $\H \B.
%_c
$
\newline
III.  Determination of the representations induced from finite
little groups
%} \label{chap2} \label{lyglygyhgyuhglyughbyhug}
%\section{How $\H \B_c$ arises
%NEW GROUPS IN GRAVITATIONAL INSTANTON THEORY\\[1mm]
%          ULTRAHYPERBOLIC BMS GROUP
}
\end{center}

%{\bf \LARGE G-Fluxes and Non-Perturbative\\[2mm]
%              Stabilisation of Heterotic M-Theory}

%\vspace{1.3truecm} \centerline{Gottfried Curio
%$^{a,}$\footnote{curio@physik.hu-berlin.de} and Axel Krause
%$^{b,}$\footnote{krause@physics.uoc.gr}} \vspace{.6truecm}

%{\em \centerline{$^a$ Humboldt-Universit\"at zu Berlin,}
%\centerline{Institut f\"ur Physik, D-10115 Berlin, Germany}}
%\vspace{.3truecm} {\em \centerline{$^b$ Physics Department,
%University of Crete,} \centerline{71003 Heraklion, Greece}}

\vspace{1.3cm} \centerline{
%Patrick J. McCarthy $^{a,}$ \footnote{
%${\rm patrick}_{-}$${\rm j}_{-}$mcc@hotmail.com}  and
Evangelos Melas $^{a,}$ \footnote{emelas@econ.uoa.gr $\&$
evangelosmelas@yahoo.co.uk }} \vspace{.6cm}
%\vspace{1.3cm}
%\vspace{.6cm}
% {\em \centerline{$^a$ Horn Lane , Woodford
%Green,} \centerline{Essex, IG8 9AF, U.K. }} \vspace{.3cm}
%{\em \centerline{$^a$ Department of Mathematical Sciences, QMW }
%\centerline{ London, U.K.}}
{\em \centerline{$^a$ Department of Economics, Unit of Mathematics and Informatics,}
\centerline{ Athens, Greece}}

%Patrick J. McCarthy

%123 Horn Lane, Woodford Green, Essex IG8 9AF, U.K.

\vspace{1.0truecm}
%%%%%%%%%%%%%%%%%%%%%%%%%%%%%%%%%%%%%%%%%%%%%%%%%%%%%%%%%
%\vspace{.4truecm}
\begin{abstract}
The ordinary Bondi$-$Metzner$-$Sachs (BMS) group $B$ is the common
asymptotic symmetry group of all asymptotically flat Lorentzian
space-times. \ As such, $B$ is the best candidate for the
universal symmetry group of General Relativity. \ However, in
studying quantum gravity, space$-$times with signatures other than
the usual Lorentzian one, and complex space$-$times, are frequently
considered. \ Generalisations of $B$ appropriate to these other
signatures have been defined earlier. In particular,
%\ Here,
%the generalisation
$\H \B$, a variant of  BMS group appropriate to
the ultrahyperbolic signature (+,+,$-$,$-$),
%is
%described in detail, and the irreducible unitary representations
%(IRs) of $B(2,2)$ are analysed
has been defined in a previous paper where it was shown that all the
strongly continuous unitary irreducible representations (IRs) of $\H \B $
can be obtained with the Wigner$-$Mackey's inducing method and that all the little groups of
$\H \B $ are compact.
%. \ It is proved that all induced
%IRs of $B(2,2)$ arise from IRs of \textit{compact }``little
%groups''.
\ Here we describe in detail all the finite little groups of $\H \B$  and we find
the IRs of $\H \B$  induced by them.
%\ These little groups, which are closed subgroups of
%$K=SO(2)\times SO(2)$, are classified here in detail, with
%particular attention paid to those of infinite order.
\end{abstract}

\bigskip \bigskip
\newpage
\pagenumbering{arabic}

\section{Introduction}

\indent
%bigcirc \!\!\;\!\!\!\!s
%\section{A close-up for $\H \B_c$}

%\section{All Little Groups are Compact
%for $\H\B_{c}$
%}
The best candidate for the universal symmetry group of General
Relativity (G.R), in any signature, is the so called
Bondi$-$Metzner$-$Sachs (BMS) group $B$. These groups have recently been
described \cite{mac1} for all possible signatures and all possible
complex versions of GR as well.
%The induced irreducible
%representations (IRs) have also been classified and constructed
%for complex GR, and, in more detail, for ultrahyperbolic GR
%\cite{macmel}.

\noindent In earlier papers \cite{mac1,macmel,mel,mel1,mel2}  it has
been argued that the IRs of the BMS group and of its
generalizations in complex space$-$times as well as in space$-$times
with Euclidean or Ultrahyperbolic signature are what really lie
behind the full description of (unconstrained) moduli spaces of
gravitational instantons.
%, via the McKay
%correspondence between Lie algebras and the A,D,E series.
Kronheimer \cite{Kronheimer, Kronheimer1} has given a description
of these instanton moduli spaces for {\it Euclidean} instantons.
However, his description only partially describes the moduli
spaces, since it still involves {\it constraints}. Kronheimer does
not solve the constraint equations, but it has been argued
\cite{mac1,mel2} that IRs of BMS group (in the relevant signature)
give an {\it unconstrained} description of these same moduli
spaces.

The original BMS group $B$ was discovered by Bondi,
Metzner and Van der Burg \cite{Bondi}
for asymptotically flat space$-$times which were axisymmetric, and
by Sachs \cite{Sachs1} for general asymptotically flat
space$-$times, in the usual Lorentzian signature. The group $\mathcal H \mathcal B$
is a different generalised BMS group, namely  one appropriate
to the `ultrahyperbolic' signature, and asymptotic flatness in
null directions introduced in \cite{Mel}.

Recall that the ultrahyperbolic version of
Minkowski space is the vector space $R^{4}$ of row vectors with 4
real components, with scalar product defined as follows. \ Let
$x,y\in R^{4}$ have components $x^{\mu }$ and $y^{\mu }$
respectively, where $\mu =0,1,2,3$. \ Define the scalar product
$x.y$ between $x$ and $y$ by
\begin{equation}
x.y=x^{0}y^{0}+x^{2}y^{2}-x^{1}y^{1}-x^{3}y^{3}.
\end{equation}
Then the ultrahyperbolic version of Minkowski space, sometimes
written $R^{2,2}$, is just $R^{4}$ with this scalar product.

In \cite{Mel} it was shown that
\begin{thrm}
\label{omorpho} The group $\H\B
%_{c}
$   can be realised as \be
\H\B
%_{c}
= L^{2}(\mathcal P,\lambda ,R)
\bigcirc\!\!\;\!\!\!\;\!\!\!\!s \  _{T} G^{2} \ee with semi$-$direct
product specified by \be
\label{an}
(T(g,h)\alpha)(x,y)=k_{g}(x)s_{g}(x)k_{h}(w)s_{h}(w)\alpha
(xg,yh), \ee
\end{thrm}
where $\alpha \in L^{2}(\mathcal P,\lambda ,R)$,
%$L^{2}(\mathcal P,\lambda ,R)$ is
the separable Hilbert space of real$-$valued functions defined on $\mathcal P$,
 and $(x,y)\in \mathcal P$.
For ease of notation, we   write $\mathcal P$
for the torus ${\rm T} \simeq P_{1}(R)\times P_{1}(R)$, $P_{1}(R)$ is the one$-$dimensional real projective space,
and $\mathcal G$ for $G\times G$, $G=SL(2,R)$. In analogy to $B$, it is natural to choose a measure $\lambda$
on $\mathcal P$ which is invariant under
the maximal compact subgroup ${ \rm S}{\rm O}(2)\times {\rm S}{\rm
O}(2)$ of $\mathcal G $.
Moreover, if $g \in G$ is
\be
\left[ \begin{array}{lr}
   a    & b \\
   c    & d
\end{array}  \right],
\ee then the components $ x_{1},x_{2} $ of ${\bf x} \in R^{2}$ transform
linearly, so that the ratio $ x= {x_{1}}/{x_{2}} $ \  transforms
fraction linearly. Writing $ xg $ \ for the transformed ratio, \be
\label{aliki} xg=
\frac {({\bf x}g)_{1}}{({\bf x}g)_{2}}=\frac
{x_{1}a + x_{2}c} {x_{1}b+x_{2}d}=
\frac {xa+c}{xb+d}. \ee
%Hence we
%have \be \psi({\bf x}g,{\bf y}h)=r({\bf x}g)s(({\bf
%x}g)_{2})r({\bf y}h) s(({\bf y}h)_{2})\alpha(xg,yh). \ee \indent
%The action $T$ of $G^{2}$ on the $\psi$s in Theorem
%\ref{hvvmkuylyhvbgyug} ( Eq. ( \ref{giorti})) induces an action,
%also denoted $T$, of $G^{2}$ on the $\alpha$s, defined by \be
%(T(g,h)\psi)({\bf x},{\bf y})=r({\bf x})s(x_{2})r({\bf y})
%s(y_{2})(T(g,h) \alpha)(x,y). \ee The last two equations give \be
%\label{agonas}
%(T(g,h)\alpha)(x,y)=k_{g}(x)s_{g}(x)k_{h}(y)s_{h}(y)\alpha(xg,yh),
%\ee
%5 where
 The factors $k_{g}(x)$ and $s_{g}(x)$ on the right hand side of (\ref{an}) are defined by
\be
\label{donkey}
k_{g}(x)=
%\frac {r({\bf x}g)}{r({\bf x})}=\frac {|({\bf x}g)_{2}|\rho(xg)}
%{|x_{2}|\rho(x)}=
\left\{ \frac {(xb+d)^{2}+(xa+c)^{2}}{1+x^{2}} \right\}
^{\frac {1}{2}},
\ee

\be
\label{mule}
s_{g}(x)=
%\frac {s(({\bf x}g)_{2})}{s(x_{2})}=\frac {x_{1}b+x_{2}d}
%{|x_{1}b+x_{2}d|}\cdot \frac {|x_{2}|}{x_{2}}=
\frac {xb+d}{|xb+d|},
\ee
 with similar formulae for $yh$, $k_{h}(y)$ \ and $s_{h}(y)$.
  %\ Note that the
%last three formulae are expressed entirely in terms of local
%coordinates $(x,y)\in~P_{1}(R)~\times~P_{1}(R).$ \ Strictly
%speaking, three more local charts are needed to cover all of
%$P_{1}(R)~\times P_{1}(R)$ \ (cf. \cite{mac1}), related to the
%above via $(x,y^{-1}), \ (x^{-1},y)$ \  and $ (x^{-1},y^{-1})$,
%but the single one used above will be sufficient here.
%Summarizing, we have
%\begin{thrm}
%\label{omorpho} The group $\H\B_{c}$   can be realised as \be
%\H\B_{c}= C^{\infty}(P_{1}(R)^{2},R)
%\bigcirc\!\!\;\!\!\!\;\!\!\!\!s \  _{T} G^{2} \ee with semi-direct
%product specified by \be
%(T(g,h)\alpha)(x,y)=k_{g}(x)s_{g}(x)k_{h}(w)s_{h}(w)\alpha
%(xg,yh). \ee
%\end{thrm}

\indent
It is well known that the topological dual of a Hilbert space can be identified
with the Hilbert space itself, so that we
%now
have $ {L^{2}}^{'}(\mathcal P,\lambda ,R) \simeq
L^{2}(\mathcal P,\lambda ,R).$
%$ A^{'}(\mathcal N) \simeq
%A(\mathcal N).$
%L^{2}(\mathcal P,\lambda ,R)
~ In fact, given a continuous linear functional
 $\phi \in {L^{2}}^{'}(\mathcal P,\lambda ,R)$, we can write,
 for $\alpha \in L^{2}(\mathcal P,\lambda ,R)$
\be
(\phi,\alpha)=<\phi,\alpha>
\ee
where the function $\phi \in L^{2}(\mathcal P,\lambda ,R)$ on the right is uniquely determined
by (and denoted by the same symbol as) the linear functional $\phi \in {L^{2}}^{'}(\mathcal P,\lambda ,R)
$ on the left.
The representation theory of
$\H\B
%_{c}
$
is governed \cite{Mel} by the dual
action $T'$ of $\mathcal G$ on the topological dual
$ {L^{2}}^{'}(\mathcal P,\lambda ,R)$ of $L^{2}(\mathcal P,\lambda ,R).$
The dual action $T'$ is defined by:
%$A'(\mathcal N)\simeq
%A(\mathcal N)$, defined typically as :
%in \S\quad :
\be <T'(g,h)\phi,\alpha>=<\phi,T(g^{-1},h^{-1})\alpha>\cdot \ee A
short calculation
%, involving a simple change of variables,
gives
%\be <T'(g,h)\phi,\alpha>=\int_{\mathcal
%P}k_{g}^{-3}(x)s_{g}(x)k_{h}^{-3}(y)
%s_{h}(y)\phi(xg,yh)\alpha(x,y){\rm d}\lambda(x,y). \ee Since this
%holds for all $\phi \in A(\N)$,
\be \label{tzatziki}
(T'(g,h)\phi)(x,y)=k_{g}^{-3}(x)s_{g}(x)k_{h}^{-3}(y)s_{h}(y)\phi(xg,yh).
\ee
%It is convenient to have the corresponding transformation law
%for the absolute value  $|\phi(x,y)|$: \be |(T'(g,h)\phi)({\bf
%m},{\bf n})|=k^{-3}_{g}({\bf m})k_{h}^{-3}({\bf n}) |\phi({\bf
%m}g,{\bf n}h)|. \ee \indent
Now, this action $T'$ of $\mathcal G$
on ${L^{2}}^{'}(\mathcal P,\lambda ,R)$,
given explicitly above, is like the action $T$ of
$\mathcal G$ on $L^{2}(\mathcal P,\lambda ,R)$, continuous. The `little group' $L_{\phi}$
of any  $\phi \in {L^{2}}^{'}(\mathcal P,\lambda ,R)$ is the stabilizer \be
L_{\phi}=\{(g,h)\in \mathcal G
%G\times G
\ | \ T'(g,h)\phi=\phi\}. \ee By
continuity, $L_{\phi}\subset \mathcal G$ is a closed subgroup.
%Clearly, if $\phi=0 $, $L_{\phi}=\mathcal G$ and the orbit of
%$\phi=0$ is just the origin.
%The resulting {\scriptsize\rm IRS} of
%$\H\B_{c}$ are trivial in the subgroup $A(\N)$, and are really
%only {\scriptsize\rm IRS} of ${\mathcal G}$. These representations
%are unphysical, since they describe quantum systems with 'zero
%supermomentum' ; they are ignored here, as in earlier work for
%$B_{c}$ and $CB(\N)_{c}$. Henceforth, then, it is  assumed that
%$\phi \neq  0$.

Attention is confined to measures on
${L^{2}}^{'}(\mathcal P,\lambda ,R)$ which are concentrated on
single orbits of the $
%G^{2}
\mathcal G
-$action $T^{\prime }.$ These measures
give rise to IRs of $ \H \B $ which are induced in a sense
generalising \cite{mac6} Mackey's \cite{Wigner,Mackey,Mackey1,Simms,Isham,Mackey2}.
In fact in \cite{Mel} it was shown that \it all \normalfont
IRs of the $\H \B$
with the Hilbert topology are
%shown to be
derivable by the inducing construction.
The inducing construction
is realized as follows. Let $ \mathcal O \subset
{L^{2}}^{'}(\mathcal P,\lambda ,R)$ be any orbit of the dual action
${\;T}^{\prime }$ of
$
%G^{2}
\mathcal G
$ on ${L^{2}}^{'}(\mathcal P,\lambda ,R).$
There is a natural homomorphism $ \mathcal O = \mathcal G \phi_{\rm{o}} \simeq
%G^{2}
\mathcal G
/L_{\phi_{\rm{o}}} $ where $L_{\phi_{\rm{o}}} $ is the `little group' of the point
$\phi_{\rm{o}} \subset \mathcal O .$ Let $U$ be a continuous
irreducible unitary representation of $L_{\phi_{\rm{o}}}$ on a Hilbert space
$D.$
%Every orbit $ \mathcal O \subset L_{e}^{2}({\rm T}^{2})^{\prime }$
%of $G^{2}$ in $L_{e}^{2}({\rm T}^{2})^{\prime }$ can be identified
%with a coset space of $G^{2}.$ Indeed, let $\zeta_{o}$ be a base
%point of $ \mathcal O .$ There is a natural homomorphism $
%\mathcal O \simeq G^{2}/L_{o} $ where $L_{o}$ is the `little
%group' of the point $\zeta_{o} \subset \mathcal O .$ Let $U$ be an
%irreducible unitary representation of $L_{o}$ on a Hilbert space
%$\mathcal H_{o}.$
Every coset space  $\mathcal G / L_{\phi_{\rm{o}}}
%\mathcal O
$ can be  equipped with a unique
class of measures which are quasi$-$invariant under the  action
${\;T}$ of $G^{2}$. Let $ \mu $ be any one  of these. Let $
D_{\mu}= L^{2}(\mathcal G/ L_{\phi_{\rm{o}}} , \mu, D)$ be the
Hilbert space of functions $  f :  \mathcal G/ L_{\phi_{\rm{o}}}  \rightarrow D$
which are square integrable with respect to $\mu.$ From a
given $\phi_{\rm{o}} $ and any continuous irreducible unitary
representation $U$ of $L_{\phi_{\rm{o}}}$ on a Hilbert space $D$
a continuous irreducible unitary representation of $ \H \B \ $
on $D_{\mu}$ can be constructed. The representation is said to
be induced from $U$ and $\phi_{\rm{o}}$.
% and is given by
%\begin{eqnarray}
%(\varrho_{o}f)(\varrho) & = & f(\varrho_{o}^{-1}\varrho),
%\nonumber \\
%(\alpha f)(\varrho) & = & e^{i\left\langle \varrho \zeta_{o},
%\alpha \right\rangle } f (\varrho),
%\end{eqnarray}
%where $\varrho,\varrho_{o} \in G^{2}$ and $\left\langle \varrho
%\zeta_{o}, \alpha \right\rangle$ is the scalar product in $
%L_{e}^{2}({\rm T}^{2}).$
Different points of an orbit $\mathcal G \phi$
have conjugate little groups and give rise to equivalent
representations of $ \H \B $.

%In the product topology of $
%L_{e}^{2}({\rm T}^{2}) \  \times  \ G^{2}, \ $ $ \ B(2,2) \ $ is
%not locally compact and as a consequence the problem of
%determining IRs of $ \ B(2,2) \ $ arising from strictly ergodic
%actions ${\;T}^{\prime }$ of $G^{2}$ on $L_{e}^{2}({\rm
%T}^{2})^{\prime }$ is hopeless.
%In \cite{Mel} All unitary irreducible representations of the Bondi-Metzner-Sachs group with the
%Hilbert topology are shown to be derivable by the inducing construction,
To conclude,
%everyone of the IRs of $\H \B $ Thus
every irreducible representation of
$\H \B$ is obtained by the inducing construction for each $\phi_{\rm{o}} \in  {L^{2}}^{'}(\mathcal P,\lambda ,R)$
and each irreducible representation $U$ of $L_{\phi_{\rm{o}}}$.
%an irreducible representation of
%$\H \B$ is obtained for each $\phi_{\rm{o}} \in  {L^{2}}^{'}(\mathcal P,\lambda ,R)$
%and each irreducible representation $U$ of $L_{\phi_{\rm{o}}}$.
%{\it every}
%representation of $ \H \B $ determined uniquely (up to
%equivalence) via induction by (1) an orbit $\mathcal O \in
%L_{e}^{2}({\rm T}^{2})^{\prime },$ (2) a class of equivalent IRs
%of any little group $L_{o},$ is irreducible \cite{mac5}.
%%It is not
%%known if there are other IRs of $ \ B(2,2) \ $ emanating from
%%strictly ergodic actions.
All the little groups $L_{\phi_{\rm{o}}}$ of $ \H \B $
are compact and they up to conjugation subgroups of
$ \rm {S O(2)} \times  \rm { S O (2)} $. They
include groups
which are finite as well as groups which are
infinite, both  connected
%as well as groups which are infinite but
and not$-$connected \cite{macmel}.
Therefore the construction of the IRs of $\mathcal H \mathcal B$
involves at the first instance the classification of all the subgroups of
$ \rm {S O(2)} \times  \rm { S O (2)} $.

%It turns out \cite{macmel} that

The infinite not$-$connected
%compact
subgroups of $ \rm {S O(2)} \times  \rm { S O (2)} $
%little groups
%$L_{\phi_{\rm{o}}}$
%subgroups  of  $\; {\rm S}{\rm
%O}(2)\times {\rm S}{\rm O}(2) \;$
were given  in \cite{macmel} and
the infinite connected
subgroups of $ \rm {S O(2)} \times  \rm { S O (2)} $
%little groups
were given in  \cite{Mel1}.
The IRs  of $\H \B$ induced from all the infinite
little groups, both connected and non$-$connected, were
constructed in \cite{Mel1}.

%In this paper we concentrate on describing in detail all the
%infinite connected little groups
%%$L_{\phi_{\rm{o}}}$
%and on finding the Hilbert spaces of the invariant functions
%$\phi_{\rm{o}}$ of all the infinite little groups $L_{\phi_{\rm{o}}}$,
%both connected and non$-$connected.
%Finally we give explicitly the operators of the IRs  of
%$\H \B$ induced from all the infinite little groups.

%In particular
The problem of
constructing the IRs of $\mathcal H \mathcal B$
induced from {\it finite} little groups
%The IRs induced from {\it
%infinite} little groups were constructed in full detail. The
%purpose of this paper is to construct the corresponding IRs for
%the {\it finite} little groups.
%It turns out (\cite{macmel}) that this problem
reduces to a seemingly very simple task \cite{macmel}; that of classifying all
subgroups of the Cartesian product group $\; {\rm C}_{{\rm n}}
\times {\rm C}_{{\rm m}}, \;$ where $\;{\rm C}_{{\rm r}}\;$ is the
cyclic group of order $\; {\rm r}, \;$ $\; {\rm r}\;$ being
finite. Surprisingly, this task is less simple than it may appear
at first sight. It turns out \cite{mel} that the solution is constructed from
the `fundamental cases' $\; {\rm n}={\rm p}^{a}, \;$ $\; {\rm
m}={\rm p}^{\beta}, \;$ (n,m are powers of the same prime), via
the prime decomposition of m and n. The classification of all the
subgroups of $\; {\rm C}_{{\rm n}} \times {\rm C}_{{\rm m}}, \;$
was given in \cite{mel}. Classifying the subgroups of
$\; {\rm C}_{{\rm n}} \times {\rm C}_{{\rm m}}$ is one thing,
constructing the IRs of $\mathcal H \mathcal B$ induced by the
finite little groups is quite another.
In this paper we fill this gap by isolating the finite little groups and by
explicitly  constructing the IRs of $\mathcal H \mathcal B$ induced by them.

%The explicit construction of the IRs of
%$\mathcal H \mathcal B$ induced from finite little groups is not
%a trivial task and it will be given in a different paper \cite{Mel1}.

%The finite little groups are

%The finite
%(compact)
%and
%The infinite not$-$connected
%%compact
%little groups
%%$L_{\phi_{\rm{o}}}$
%%subgroups  of  $\; {\rm S}{\rm
%%O}(2)\times {\rm S}{\rm O}(2) \;$
%were given  in \cite{macmel}.
%In this paper we concentrate on describing in detail all the
%infinite connected little groups
%%$L_{\phi_{\rm{o}}}$
%and on finding the Hilbert spaces of the invariant functions
%$\phi_{\rm{o}}$ of all the infinite little groups $L_{\phi_{\rm{o}}}$,
%both connected and non$-$connected.
%Finally we give explicitly the operators of the IRs  of
%$\H \B$ induced from all the infinite little groups.

%Here, we study in more detail, the infinite connected
%(compact) subgroups of $\; {\rm S}{\rm O}(2)\times {\rm S}{\rm
%O}(2) \;$.

%. The little groups $L_{o}$ for $\;B(2,2)\;$ are  the
%closed subgroups of $\;K=SO(2)\times SO(2)\;$ which contain the
%element $(-I,-I)$. \ These are (A) $K$ itself, (B) \ a class of
%one dimensional not connected Lie groups  which are described in
%detail in \cite{macmel} , and (C) \ all finite subgroups
%containing $(-I,-I)$. The finite subgroups of $\;K=SO(2)\times
%SO(2)\;$ are precisely the subgroups of ${\rm C}_{\rm n} \times
%{\rm C}_{\rm m}$ where both n and m are finite. These subgroups
%are not given in \cite{macmel} and so we proceed now to construct
%them explicitly.

\indent
The paper is organised as follows: In Section \ref{s1}
we prove that the elementary domains for the actions
of $G={\rm C}_{{\rm n}} \times {\rm C}_{{\rm m}}$ and of any of its
subgroups $F$ on the torus $\mathcal P$ are related in a very simple fashion:
the elementary domain for the action of $F$ on $\mathcal P$ can be obtained
by letting a set of representatives of the left cosets of $F$
in the coset space $G/F$ act on an elementary domain E for the action of
$G$ on $\mathcal P$.
%an elementary domain for the action of
%any of the subgroups of $\; {\rm C}_{{\rm n}}
%\times {\rm C}_{{\rm m}} \;$ on  the torus $\mathcal P$ is obtained by le
%the infinite connected subgroups of $ \rm  { S O(2)} \times \rm {S O (2)}. $
%In Section \ref{terolokiju} it is , this spinor version is given.
In Section \ref{s2} we use this result to find  the Hilbert spaces $\mathcal H(F)$
of all invariant vectors for each of the finite subgroups $F$ of $G
%={\rm C}_{{\rm n}} \times {\rm C}_{{\rm m}}
$.
%the infinite potential little
%groups, both connected and non$-$connected.
In Section \ref{s3} we find all the finite potential  little groups.
%prove
%that \it all \normalfont the infinite potential little groups are actual.
In Section \ref{s4} we prove that all the finite potential  little groups are actual.
%construct explicitly the IRs of $\H \B$ induced
%from the infinite little groups.
In Section \ref{s5} we describe explicitly all the finite little groups.
%make some remarks about the  IRs of $\H \B$
%constructed in Section \ref{s4} by the inducing method.
In Section \ref{s6} we give the form of the IRs  of $\mathcal H \mathcal B$  induced
by the IRs of the finite little groups and the corresponding invariant characters .
% elementary regions associated
% with the action of the little groups on the torus $P_{1}(R) \times P_{1}(R).$
Finally,
in Section \ref{s8} we make some remarks about the  IRs of $\H \B$
constructed in Section \ref{s6} by the inducing method.
%construct
%$\;\mathcal G-{\rm quasi}-{\rm invariant}$ measures on the orbits $
%%\mathcal G
%%/L_{\phi} \approx
%\mathcal G \phi_{\rm{o}}$ of the action of the little groups $L_{\phi_{\rm{o}}} $
%on the topological dual  ${L^{2}}^{'}(\mathcal P,\lambda ,R)$, measures which are necessary
%%which are necessary
%in order
%to give the operators of the induced IRs of $\H \B$ in Section \ref{s4}.

\section{Elementary regions  for Finite Groups}
\label{s1}
\label{uyfvfvjytddxhtesxe}
%We need some results regarding group
%actions. It was shown in chapter \ref{lyglygyhgyuhglyughbyhug}
%that
The representation theory of $\;\mathcal H \mathcal B
%_{c}
\;$
is governed by the dual action  $T'$ of $\mathcal G
%=G\times
%G
$
%$\;G={\rm S}{\rm L}(2,R) \;$,
on
%$A'(\mathcal N)\simeq
%A(\mathcal N)$,
${L^{2}}^{'}(\mathcal P,\lambda ,R)$
given by Eq. (\ref{tzatziki})
$$ \label{ulbkyuvrdfcrfxcr}
(T'(g,h)\phi)(x,y)=k_{g}^{-3}(x)s_{g}(x)k_{h}^{-3}(y)s_{h}(y)\phi(xg,yh).
$$
%We need some results regarding group
%actions.
We have an action \be \label{jkbb.jkb.b}(x,y) \longmapsto (xg,yh)
\ee of $\; \mathcal G \;$  on the torus $\; {\rm T}
%^{2}
\simeq
P_{1}(R) \times P_{1}(R) \;$.
We need at this point to recall some results regarding group
actions.
Recall that if the set $\; \mathcal
X \;$ is a $\; {\rm  G}$$-$$ {\rm space} \;$ then there is a
natural bijective correspondence between the set of actions of $\;
{\rm  G } \;$ on  $\; \mathcal X \;$ and the set of homomorphisms
from $\; {\rm G} \;$ into $\; {\rm Aut}(\mathcal X), \;$ where $\;
{\rm Aut}(\mathcal X) \;$ is the group of automorphisms of the set
$\; \mathcal X . \;$ Indeed, let the map \be \label{cukilok}
\mathcal X  \times  {\rm G}       \mapsto \mathcal X, \ee from $\;
\mathcal X  \times  {\rm G} \;$ to $\; \mathcal X \;$ be a right
action of $\;{\rm G} \;$ on the set $\; \mathcal X, \;$ with the
image of $\;(x,g)\;$ being denoted by $\; xg. \;$ Then the map
(\ref{cukilok}) satisfies the following conditions:
\begin{itemize}
\item{$\; xe=x \;$ for every  $\; x  \in \mathcal X. \;$}
\item{$\; x(g_{1}g_{2})=(xg_{1})g_{2} \;$ for every $\;g_{1},g_{2}
\in {\rm G} \;$ and $\; x \in \mathcal X. \;$}

\end{itemize}
Now for each $\; g \in {\rm G} \;$ define the map
\begin{eqnarray}
s_{g} & : & \mathcal X  \rightarrow  \mathcal X \nonumber \\
s_{g}(x) & = & xg \quad {\rm for} \quad {x \in \mathcal X}.
\end{eqnarray}
One can easily show (see e.g \cite{alperin}, p.28) that
$\;s_{g}\;$ is a member of the group $\; {\rm Aut}(\mathcal X).
\;$ Furthermore, the second condition in the definition of  a
group action ensures that we have $\;
s_{g_{1}g_{2}}=s_{g_{1}}\circ s_{g_{2}} \;$ for any $\;
g_{1},g_{2} \in {\rm G}. \;$ Consequently there exists an
homomorphism $\; h \;$such that
\begin{eqnarray}
h & : & {\rm G} \rightarrow {\rm Aut}( \mathcal X) \nonumber \\
\label{vgdgtuyikjkknm} h(g) & = & s_{g}.
\end{eqnarray}
Conversely, suppose that $\; h : {\rm G} \rightarrow {\rm
Aut}(\mathcal X) \;$ is a homomorphism. We define a map from $\;
\mathcal X \times {\rm G} \rightarrow \mathcal X \;$ by sending
$\; (x,g) \;$ to $\; h(g)(x). \;$ One can easily check that this
map is an action of $\; {\rm G} \;$ on $\; \mathcal X \;$. It is
an easy exercise  now to show that the correspondence between the
actionss (\ref{cukilok}) and the homomorphisms
(\ref{vgdgtuyikjkknm}) is bijective.

%For ease of notation, we   write $\mathcal P$
%for the torus ${\rm T} \simeq P_{1}(R)\times P_{1}(R)$, $P_{1}(R)$

%There is a natural bijective
%correspondence between the set of actions of $\; \mathcal G \;$ on
%$\; {\rm T}^{2} \;$ and the set of homomorphisms from $\; \mathcal
%G \;$ into the group $\; {\rm Aut}({\rm T}^{2}) \;$ of
%automorphisms of $\; {\rm T}^{2} \;$.
In the problem under consideration $\; {\rm G} \equiv \mathcal G
\;$ and $\; \mathcal X \equiv \mathcal P \simeq P_{1}(R)\times P_{1}(R), \;$ and,
there is a
natural bijective correspondence between the set of actions of $\;
\mathcal G \;$ on $\; \mathcal P \;$ and the set of homomorphisms
from $\; \mathcal G \;$ into the group $\; {\rm Aut}(\mathcal P)
\;$ of automorphisms of $\; \mathcal P. \;$ Let $\; h_{m} \;$ be
the homomorphism which is associated with the specific action
(\ref{jkbb.jkb.b}). The following Proposition gives the kernel of
$\; h_{m}. \;$
\begin{prop}
 The action from the right of $\;
\mathcal
 G \;$ on $\; P_{1}(R) \times P_{1}(R)
\;$ \be ((x,y)(g,h)) \mapsto (xg,yh), \ee where $\;(g,h) \in
\mathcal G \;$,  $g \in G$ is $ \left[ \begin{array}{lr}
   a    & b \\
   c    & d
\end{array}  \right],
$ and $ \label{uigbgkyuvmkuyhv}  xg=
%\frac {({\bf x}g)_{1}}{({\bf
%x}g)_{2}}=
%\frac {x_{1}a + x_{2}c} {x_{1}b+x_{2}d}=
\frac {xa+c}{xb+d} $ is not effective. The kernel of $\; h_{m} \;$
is \be {\rm K}= \left \{ (g,h) \in \mathcal G \quad | \quad
g=\underline{+}{\rm I} \quad , \quad h=\underline{+}{\rm I} \right
\}, \ee where $\;{\rm I}\;$ is the identity element $ \left[
\begin{array}{lr}
   1    & 0 \\
   0    & 1
\end{array}  \right]
$ of $\;G.\;$
\end{prop}
\begin{prf}
\normalfont If $\;(g,h)\in {\rm K}, \;$ $\;xg=x \;$ and $\;yh=y
\;$ for all $\;(x,y).\;$ So, if $g= \left[ \begin{array}{lr}
   a    & b \\
   c    & d
\end{array}  \right],
$ $   xg=
%\frac {({\bf x}g)_{1}}{({\bfx
%x}g)_{2}}=
%\frac {x_{1}a + x_{2}c} {x_{1}b+x_{2}d}=
\frac {xa+c}{xb+d}= x$ for all x. Taking $\; x=0 \;$ gives $\;
\frac{c}{d}=0, \;$ and so, $\; c=0. \;$ Taking $\; x= \infty \;$
gives $\; \frac {a}{b}= \infty, \;$ and therefore, $\; b=0. \;$ So
we obtain that $\;g= \left[ \begin{array}{ll}
   a    & 0 \\
   0    & a^{-1}
\end{array}  \right]
\;$ (it is det$g$=1). Therefore we have $\;xg=\frac{ax +
0}{0x+a^{-1}}=a^{2}x=x. \;$ Setting $\; x=1 \;$ gives
$\;a^{2}=1,\;$ so $\; a=\underline {+}1 \;$ and $\; g=
\underline{+}{\rm I}. \;$ Similarly $\; h= \underline{+}{\rm I}.
\;$ This completes the proof.
\end{prf}
The group K is a normal subgroup of $\; \mathcal G. \;$ The not
effective action of $\; \mathcal G \;$ on $\; \mathcal P \simeq
P_{1}(R) \times P_{1}(R) \;$ passes naturally to an effective
action of the group $\; \mathcal G / {\rm K} \;$ on the torus $\; \mathcal P. \;$
This last action is given by \be
\label{bgdjutrfgsss}(x,y)((g,h){\rm K})=(xg,yh), \ee where $\;
(g,h){\rm K} \;$ denotes an element of the coset space $\;
\mathcal G / {\rm K}. \;$ We will prove now that the  action
(\ref{bgdjutrfgsss}) not only is it effective but also fixed point
free (f.p.f). To prove this we first need the following lemma.

\begin{lm}
\label{nuhyyyteooeleodki}
 When $\; \mathcal G \;$ is restricted to
its subgroup $ \; {\rm S}{\rm O}(2) \times {\rm S}{\rm O}(2) \;$
the action (\ref{bgdjutrfgsss}) is f.p.f.
\end{lm}
\normalfont

\begin{prf}

\normalfont

 To simplify notation we denote
$ \; \left(  \left( \begin{array}{cc}
\cos\theta & \sin\theta \\
\!\!\!\!-\sin\theta & \cos\theta
\end{array}
\right),
 \left( \begin{array}{cc}
\cos\varphi & \sin\varphi \\
\!\!\!\!\!-\sin\varphi & \cos\varphi
\end{array}
\right) \right) \; $ by $\; \left ( R \left ( \theta \right ) ,
 R  \left ( \varphi
  \right
) \right ). \; $
%It was shown in chapter
%\ref{lyglygyhgyuhglyughbyhug}, (Eq.(\ref{choma})), that
If $\;
(g,h)= \left ( R \left ( \theta_{{\rm o}} \right ) ,
 R  \left ( \varphi_{{\rm o}}
  \right
) \right ) \; $ is an element of $\; {\rm S}{\rm O}(2) \times {\rm
S}{\rm O}(2) \;$ then the action (\ref{jkbb.jkb.b}) reads
$$
(\rho, \sigma) \longmapsto (\rho + 2\theta_{{\rm o}}, \sigma +
2\varphi_{{\rm o}}),
$$
where $\; x= {\rm cot} \frac {\rho }{2}, \;  \; y= {\rm cot} \frac
{\sigma}{2}. \;$ Suppose that for some point $ \; (\rho_{{\rm
o}},{\sigma}_{{\rm o}}) \;$ of the torus $ \; P_{1}(R) \times
P_{1}(R) \;$ we have \be \label {vbbcnncncmjok}(\rho_{{\rm o}} +
2\theta_{{\rm o}}, \sigma_{{\rm o}} + 2\varphi_{{\rm
o}})=(\rho_{{\rm o}},\sigma_{{\rm o}}), \ee for some element $\;
\left ( R \left ( \theta_{{\rm o}} \right ) ,
 R  \left ( \varphi_{{\rm o}}
  \right
) \right ) \; $ of $\; {\rm S}{\rm O}(2) \times {\rm S}{\rm O}(2).
\;$ Since, $\; \rho_{{\rm o}} \;$ and $\; \sigma_{{\rm o}} \;$
are defined only mod $\; 2\pi \;$ Eq. (\ref {vbbcnncncmjok}) can
be satisfied for some point $ \; (\rho_{{\rm o}},{\sigma}_{{\rm
o}}) \;$ of the torus $ \; P_{1}(R) \times P_{1}(R) \;$ if and
only if $\; \left ( R \left ( \theta_{{\rm o}} \right ) ,
 R  \left ( \varphi_{{\rm o}}
  \right
) \right ) \; $ is an element of the kernel K of the homomorphism
$\; h_{m}. \;$ This completes the proof.
\end{prf}

%\vspace{0.2cm}

\noindent
 Now we are ready to prove the following

\begin{lm} The action
(\ref{bgdjutrfgsss}) is f.p.f.
\end{lm}
\normalfont

\begin{prf}
\normalfont
 Suppose that for a point $\; (x_{{\rm o}},y_{{\rm o}})
\;$ of $\; P_{1}(R) \times P_{1}(R) \;$ and an element $\;
(g_{{\rm o}},h_{{\rm o}}) \;$ of $\; \mathcal G \;$ we have \be
\label{kolokijnmuhhyg} (x_{{\rm o}},y_{{\rm o}})(g_{{\rm
o}},h_{{\rm o}})= (x_{{\rm o}}g_{{\rm o}}, y_{{\rm o}}h_{{\rm
o}})=(x_{{\rm o}},y_{{\rm o}}). \ee The element $ \; g_{{\rm o}}=
\left[
\begin{array}{lr}
   a    & b \\
   c    & d
\end{array}  \right], \;
$ of $\; {\rm S}{\rm L}(2,R) \;$ can always be written in the form
\be g_{{\rm o}}=u_{{\rm o}}\delta_{{\rm o}} w_{{\rm o}},\quad
\delta_{{\rm o}} =\delta_{{\rm o}}(t)= \left[
\begin{array}{lr}
e^{t_{{\rm o}}/2} & 0 \\
   0    & e^{-t_{{\rm o}}/2}
\end{array}
\right ], \ee where $\;u_{{\rm o}},w_{{\rm o}} \in {\rm S}{\rm
O}(2)\;$ and $\;t_{{\rm o}}\;$ is real. We have \be x_{{\rm
o}}\delta_{{\rm o}}= \frac{x_{{\rm o}}e^{t_{{\rm
o}}/2}}{e^{-t_{{\rm o}}/2}}=x_{{\rm o}}e^{t_{{\rm o}}}. \ee When
$$
x_{{\rm o}}e^{t_{{\rm o}}}=x_{{\rm o}}
$$
and $\; x_{{\rm o}} \neq 0 \;$ we obtain
$$
e^{t_{{\rm o}}}=1\Leftrightarrow t_{{\rm o}}=0.
$$
From the last Equation and Lemma \ref{nuhyyyteooeleodki} we
conclude that for almost all $\;(x_{{\rm o}},y_{{\rm o}}) \;$ when
Eq. (\ref{kolokijnmuhhyg}) is satisfied then
$$
g_{{\rm o}}={\rm I} \quad {\rm or} \quad g_{{\rm o}}=-{\rm I}.
$$
Similarly,
% A similar argument applies to $\; h_{{\rm o}} \;$.
$$
h_{{\rm o}}={\rm I} \quad {\rm or} \quad h_{{\rm o}}=-{\rm I}.
$$
We conclude that for almost all $\;(x_{{\rm o}},y_{{\rm o}}) \;$
when Eq. (\ref{kolokijnmuhhyg}) is satisfied then
$$
(g_{{\rm o}},h_{{\rm o}}) \in {\rm K},
$$
where K is the kernel of the homomorphism $\; h_{m}. \;$ This
completes the proof.
\end{prf}

It will prove convenient to recall at this point one definition
and two propositions from elementary group theory ( see e.g. \cite
{alperin} pages 4, 6 and 11).
The symbol $\; \leq \;$ denotes subgroup,  the symbol $\;
\underline {\triangleleft} \;$ denotes normal subgroup,  whereas
the symbol $\; \cong \;$ denotes isomorphism.
%\vspace{0.4cm}

%\vspace{0.3cm}
%\noindent {\bf Definition}
\begin{de}
If $\; X \;$ and $\; Y \;$ are subgroups of a group $\; G, \;$
then we define the product  of $\; X \;$ and $\; Y \;$ in $\; G
\;$ to be $\; XY=\{{\rm x}{\rm y} \; | \; {\rm x} \in X \; , \;
{\rm y} \in Y \} \subseteq G. \; $
\end{de}

\begin{prop}
\label{byhjujjkiiklomnd}
 Let H and K  be subgroups of a group $\;
G. \;$ If $\; {\rm K} \;
 \underline {\triangleleft} \; G, \;$ then $\; {\rm H}{\rm K} \leq
G \;$ and $\; {\rm H} \cap {\rm K} \;  \underline{\triangleleft}
\; {\rm H}; \;$ if also $\; {\rm H} \; \underline{\triangleleft}
\; G, \;$ then $ \; {\rm H}{\rm K} \; \underline{\triangleleft} \;
G \; $and $\; {\rm H} \cap {\rm K} \; \underline{\triangleleft} \;
G. \;$
\end{prop}

\begin{prop}
\label{vyjkkilopllkio} Let $ \; G \;$ be a group. If $\; {\rm K}
\; \underline{\triangleleft} \; G \;$ and $ \; {\rm H} \leq G, \;$
then $\; {\rm H}{\rm K}/{\rm K} \cong {\rm H}/{\rm H} \cap {\rm
K}. \;$
\end{prop}
%The symbol $\; \leq \;$ denotes subgroup , the symbol $\;
%\underline {\triangleleft} \;$ denotes normal subgroup , whereas
%the symbol $\; \cong \;$ denotes isomorphism. \vspace{0.4cm}

Henceforth we will assume that the group $\; G \;$ is finite. Let
$\; |G| \;$ denote the order of a group G. We can prove now the
following

\begin{prop}
Let $ \; G \;$ be a finite group. If $\; {\rm K} \;
\underline{\triangleleft} \; G \;$ and $ \; {\rm H} \leq G, \;$
then the number $\; |G|/|{\rm H}| \;$ is divisible by the number
$\;|{\rm K}|/|{\rm K} \cap {\rm H}|. \;$
\end{prop}

\begin{prf}
\normalfont The number $\; |G|/|{\rm H}| \;$ is divisible by the
number $\;|{\rm K}|/|{\rm K} \cap {\rm H}| \;$ if and only if the
number $\; |G|/|{\rm K}| \;$ is divisible by the number $\;|{\rm
H}|/|{\rm K} \cap {\rm H}|. \;$ According to Proposition
\ref{byhjujjkiiklomnd} the subset $\; {\rm H}{\rm K} \;$ of $\; G
\;$ is a subgroup of $\; G. \;$ The coset spaces $\; G/{\rm K} \;$
and $\; {\rm H}{\rm K}/{\rm K} \;$ can be equipped with a group
structure since by assumption $\; {\rm
K}\;\underline{\triangleleft} \; G. \;$ Since $\; {\rm H}{\rm K}
\leq G, \;$ the group $\; {\rm H}{\rm K}/{\rm K} \;$ is a subgroup
of $\; G/{\rm K}. \;$ Therefore the number $\; |G|/|{\rm H}| \;$
is divisible by the number $\; |{\rm H}{\rm K}|/|{\rm K}|. \;$ But
according to Proposition \ref{vyjkkilopllkio}
$$
 |{\rm H}{\rm K}|/|{\rm K}| = |{\rm
H}|/|{\rm K} \cap {\rm H}|.
$$
This completes the proof.
\end{prf}

There are $\; |{\rm H}|^{\frac{|G|}{|{\rm H}|}} \;$ different ways
to choose coset representatives from the coset space $\; G/{\rm H}.
\;$ Let $\; \mathcal S \;$ be the set
$$
\mathcal S=\left \{ S_{1},S_{2},...,S_{|{\rm H}|^{\frac{|G|}{|{\rm
H}|}}} \right \}
$$
whose elements are the different choices of coset representatives
from the coset space $\; G/{\rm H}. \;$ For our needs we need to
select specific elements of $\; \mathcal S. \;$ We proceed now to
give  an explicit description of these elements.

Consider the coset space $\; {\rm K}/{\rm K} \cap {\rm H}. \;$
There are $\; |{\rm K}\cap {\rm H}|^{\frac{|K|}{|{\rm K} \cap {\rm
H}|}} \;$ different ways to choose coset representatives from the
coset space $\; {\rm K}/{\rm K} \cap {\rm H}. \;$ Let $\; \Im \;$
be the set
$$
\Im =\left \{ \sigma_{1},\sigma_{2},...,\sigma_{|{\rm K}\cap {\rm
H}|^{\frac{|K|}{|{\rm K} \cap {\rm H}|}}} \right \}
$$
whose elements are the different choices of coset representatives
from the coset space $\; {\rm K}/{\rm K} \cap {\rm H}. \;$ Let us
denote by
$$
g_{i}\sigma_{j}
$$
the elements of $\; G \;$ which are obtained by multiplying the
element $\; g_{i} \;$ of $\; G \;$ with the specific choice of
coset representatives $\; \sigma_{j} \;$($\;j \in \left \{
1,2,...,|{\rm K}\cap {\rm H}|^{\frac{|K|}{|{\rm K} \cap {\rm H}|}}
\right \} \;$). Therefore $\; g_{i} \sigma_{j} \;$ denotes
collectively $ \; \frac{|K|}{|{\rm K} \cap {\rm H}|} \;$ elements
of $\; G \;$. We write them explicitly as follows
$$
g_{i}\sigma_{j1},g_{i}\sigma_{j2},...,g_{i}\sigma_{j \omega },
$$
where $\; \omega \equiv  \frac{|K|}{|{\rm K} \cap {\rm H}|} \;$.
The following 5 steps describe in algorithmic fashion the way we
make our specific choices of elements of $\; \mathcal S. \;$
\begin{enumerate}
\item{Pick up an element $\; g_{i_{1}} \;$ of $\; G \;$ and an
element $\; \sigma_{i_{1}} \;$ of $\; \Im \;$ and construct the
cosets \be \label{mncccvffder} g_{i_{1}}\sigma_{i_{1}1}{\rm
H},g_{i_{1}}\sigma_{i_{1}2}{\rm
H},...,g_{i_{1}}\sigma_{i_{1}\omega}{\rm H} \ee of the coset space
$\; G / {\rm H}. \;$} \item{ Choose an element $\; g_{i_{2}} \;$
of $\; G \;$ which {\it does not} \normalfont belong to the
previous cosets (\ref{mncccvffder}) and an element $\;
\sigma_{i_{2}} \;$ of $\; \Im. \;$ The element $\; \sigma_{i_{2}}
\;$ might be identical to $\; \sigma_{i_{1}} \;$ or different from
it. Construct the cosets \be \label{dsgghhhnbvc}
g_{i_{2}}\sigma_{i_{2}1}{\rm H},g_{i_{2}}\sigma_{i_{2}2}{\rm
H},...,g_{i_{2}}\sigma_{i_{2}\omega}{\rm H} \ee of the coset space
$\; G / {\rm H}. \;$} \item{ Choose an element $\; g_{i_{3}} \;$
of $\; G \;$ which {\it does not} \normalfont belong to the
previous cosets (\ref{mncccvffder}) and (\ref{dsgghhhnbvc})
 and an element $\;
\sigma_{i_{3}} \;$ of $\; \Im. \;$ The element $\; \sigma_{i_{3}}
\;$ might be identical either to $\; \sigma_{i_{1}} \;$ or to $\;
\sigma_{i_{2}}, \;$  or, it could be  different from both of them.
Construct the cosets \be \label{nmmmcmcmcmcmmc}
g_{i_{3}}\sigma_{i_{3}1}{\rm H},g_{i_{3}}\sigma_{i_{3}2}{\rm
H},...,g_{i_{3}}\sigma_{i_{3}\omega}{\rm H} \ee of the coset space
$\; G / {\rm H}. \;$} \item{Repeat the same procedure $\;
\frac{\frac{|G|}{|{\rm H}|}}{\frac{|{\rm K}|}{|{\rm K} \cap {\rm
H}|}}= \frac { |G| |{\rm H} \cap {\rm K}|}{|{\rm H}|{|{\rm K}|}}
\;$ times and obtain finally the following set of cosets \be
\label{mmmmkkkklolol} g_{i_{\theta}}\sigma_{i_{\theta}1}{\rm
H},g_{i_{\theta}}\sigma_{i_{\theta}2}{\rm
H},...,g_{i_{\theta}}\sigma_{i_{\theta}\omega}{\rm H} \ee of the
coset space $\; G / {\rm H}, \;$ where $\; \theta \equiv
\frac{\frac{|G|}{|{\rm H}|}}{\frac{|{\rm K}|}{|{\rm K} \cap {\rm
H}|}}= \frac { |G| |{\rm H} \cap {\rm K}|}{|{\rm H}|{|{\rm
K}|}}.\;$} \item{ Consider the following set $\; S_{i} \;$ of
elements of $\; G \;$ \begin{eqnarray} S_{i} & = & \left \{
g_{i_{1}}\sigma_{i_{1}1},g_{i_{1}}\sigma_{i_{1}2},...,g_{i_{1}}\sigma_{i_{1}\omega},
g_{i_{2}}\sigma_{i_{2}1},g_{i_{2}}\sigma_{i_{2}2},...,g_{i_{2}}\sigma_{i_{2}\omega},
\right .  \nonumber \\ && \label{nujjjhyghbn}\left .
g_{i_{3}}\sigma_{i_{3}1},g_{i_{3}}\sigma_{i_{3}2},...,g_{i_{3}}\sigma_{i_{3}\omega},...,
 g_{i_{\theta}}\sigma_{i_{\theta}1},g_{i_{\theta}}\sigma_{i_{\theta}2},...,g_{i_{\theta}}\sigma_{i_{\theta}\omega}
\right \}. \end{eqnarray} }
\end{enumerate}
This completes the construction.

\vspace{0.5cm}

\noindent Let $\; S_{i'} \;$ be a set of elements of $\; G, \;$
which if differs from $\; S_{i}, \;$ it differs in the choice of
the elements $\; \sigma_{\xi},\;$ $\; \xi \in \{1,2,...,\theta \},
\;$ chosen from the set $\; \Im. \;$ We will write \be
\label{vbcncnmx} S_{i}\sim S_{i'}. \ee

\noindent Let $\; \mathcal C_{1} \;$ and $\; \mathcal C_{2} \;$ be
the following sets of cosets

\be \mathcal C_{1} = \left \{ g_{i}\sigma_{i_{1}1}{\rm
H},g_{i}\sigma_{i_{1}2}{\rm H},...,g_{i}\sigma_{i_{1}\omega}{\rm
H} \right \} \ee and \be \mathcal C_{2} = \left \{
g_{i}\sigma_{i_{2}1}{\rm H},g_{i}\sigma_{i_{2}2}{\rm
H},...,g_{i}\sigma_{i_{2}\omega}{\rm H} \right \}, \ee where $\;
\sigma_{1} \;$ and $\; \sigma_{2} \;$ are {\it distinct}
\normalfont  elements of $\; \Im. \;$ Then we have the following.

\begin{prop}\label{nujkilokmnjhyu}The following are true
\begin{enumerate}
\item{ Every element of $\; \mathcal C_{1} \;$ belongs to $\;
\mathcal C_{2} \;$ (and vice versa).} \item{ The set $\; S_{i} \;$
is an element of $\; \mathcal S, \;$ i.e., the elements  of $\;
S_{i} \;$ are coset representatives of the coset space $\; G /
{\rm H}. \;$} \item{ The relation (\ref{vbcncnmx}) is an
equivalence relation.}
\end{enumerate}
\end{prop}

\begin{prf}
\normalfont \hspace{0.5cm} 1. \hspace{0.5cm} Choose an element $
\; g_{i}\sigma_{i_{1}j}{\rm H} \;$ of $\; \mathcal C_{1} \;$ (so
$\; j \in \{1,2,...,\omega. \}\; $) This element belongs to $\;
\mathcal C_{2} \;$ if and only if for given $\; \sigma_{i_{1}j}
\;$ there always exists an element $ \; g_{i}\sigma_{i_{2}j'}{\rm
H} \;$ of $\; \mathcal C_{2} \;$ ( $\; j' \in \{1,2,...,\omega
\}\;$) which satisfies \be \label{nhygtfgtnh}
 g_{i}\sigma_{i_{1}j}=g_{i}\sigma_{i_{2}j'}h \Leftrightarrow \sigma_{i_{1}j}=\sigma_{i_{2}j'}h
\ee for some $\; \sigma_{i_{2}j'} \;$ and some $ \; h \in {\rm H
}.\;$ Now the element $\;\sigma_{i_{1}j}\;$ of $\; {\rm K} \;$
{\it always} \normalfont belongs to some (a unique) coset of the
coset space $\; {\rm K} / {\rm K} \cap {\rm H}. \;$ Therefore,
always there exists an element $\; \sigma_{i_{2}j'} \;$ of $\;
\sigma_{i_{2}} \;$  such that $\; \sigma_{i_{1}j} \in
\sigma_{i_{2}j'}({\rm K} \cap {\rm H}). \;$ Consequently there
always exist (unique) elements $\; \sigma_{i_{2}j'} \;$ of $\;
\sigma_{i_{2}} \;$ and $\; h \in {\rm K} \cap {\rm H} \;$ which
satisfy Eq. (\ref{nhygtfgtnh}). \vspace{0.5cm}

\hspace{1.3cm} 2. \hspace{0.5cm} Consider the cosets \be
\label{bsretatttadf} g_{i_{\nu}}\sigma_{i_{\nu}1}{\rm
H},g_{i_{\nu}}\sigma_{i_{\nu}2}{\rm
H},...,g_{i_{\nu}}\sigma_{i_{\nu}\omega}{\rm H},
 \ee
where, $\; \nu \in \{1,2,...,\theta \}. \;$  Any two $\;
g_{i_{\nu}}\sigma_{i_{\nu}\mu_{1}}{\rm
H},\;$ $\;g_{i_{\nu}}\sigma_{i_{\nu}\mu_{2}}{\rm H},\;$
$\;\mu_{1},\mu_{2} \in \{1,2,...,\omega \}\;$ of them are
distinct. Indeed, suppose they are not. Then they coincide, and
for every $\; h_{1} \in {\rm H}, \;$  there always exists a unique
$\; h_{2} \in {\rm H}, \;$  such that \be \label{xaaaqwerfvcdsxa}
g_{i_{\nu}}\sigma_{i_{\nu}\mu_{1}}h_{1}=
g_{i_{\nu}}\sigma_{i_{\nu}\mu_{2}}h_{2} \Leftrightarrow
\sigma_{i_{\nu}\mu_{2}}^{-1}\sigma_{i_{\nu}\mu_{1}}=h_{2}h_{1}^{-1}.
\ee
 Since
$\; \sigma_{i_{\nu}\mu_{2}}^{-1}\sigma_{i_{\nu}\mu_{1}}=\sigma \in
{\rm K} \;$ we have that $\;h_{2}h_{1}^{-1}=\sigma \in {\rm K}
\cap {\rm H}. \;$ From Eq. (\ref{xaaaqwerfvcdsxa}) we obtain
$$
\sigma_{i_{\nu}\mu_{1}}= \sigma_{i_{\nu}\mu_{2}}\sigma, \qquad
\qquad \sigma \in {\rm K} \cap {\rm H}.
$$
Therefore, $\; \sigma_{i_{\nu}\mu_{1}}, \sigma_{i_{\nu}\mu_{2}}
\;$ belong to the same coset of the coset space $\; {\rm K}/{\rm
K} \cap {\rm H}. \;$ Contradiction. We conclude that any two $\;
g_{i_{\nu}}\sigma_{i_{\nu}\mu_{1}}{\rm
H},\;$ $\;g_{i_{\nu}}\sigma_{i_{\nu}\mu_{2}}{\rm H}\;$,
$\;\mu_{1},\mu_{2} \in \{1,2,...,\omega \},\;$ of the cosets
(\ref{bsretatttadf}) are different from one another.

Choose now a coset $\; g_{i_{\tau}}\sigma_{i_{\tau}j}{\rm H} \;$
and without loss of generality assume that $\; \tau > \nu. \;$ The
coset $\; g_{i_{\tau}}\sigma_{i_{\tau}j}{\rm H} \;$ is different
from the cosets (\ref{bsretatttadf}). Indeed, suppose that it
coincides with one of them, say
$\;g_{i_{\nu}}\sigma_{i_{\nu}k}{\rm H}\;$ (where $\; k \in
\{1,2,...,\omega \}).\;$ Then for every $\;h_{1} \;$ there always
exists a unique $\; h_{2} \in {\rm H} \;$ such that \be
g_{i_{\tau}}\sigma_{i_{\tau}j}h_{1}=
g_{i_{\nu}}\sigma_{i_{\nu}k}h_{2} \Leftrightarrow
g_{i_{\tau}}\sigma_{i_{\tau}j}=
g_{i_{\nu}}\sigma_{i_{\nu}k}h_{2}h_{1}.^{-1}
 \ee
By writing $\;h_{2}h_{1}^{-1}=h \in {\rm  H} \;$ the last equation
gives \be \label{zassssdcxsda} g_{i_{\tau}}=
g_{i_{\nu}}\sigma_{i_{\nu}k}(h \sigma_{i_{\tau}j}^{-1}h^{-1})h.
\ee Since, $\;\sigma_{i_{\tau}j}^{-1}\; \in {\rm K}  \;  $, and
since, $\;{\rm K} \; \underline{\triangleleft} \; G\;$we have that
$\; h \sigma_{i_{\tau}j}^{-1}h^{-1} \in {\rm K}, \;$ so say, $\; h
\sigma_{i_{\tau}j}^{-1}h^{-1}= \sigma' \in {\rm K}. \;$ By writing
$\;\sigma_{i_{\nu}k} \sigma'= \overline{\sigma} \;$ Eq.
(\ref{zassssdcxsda}) gives
$$
g_{i_{\tau}}= g_{i_{\nu}}\overline{\sigma}h
$$
for some $\; h \in {\rm H}. \;$ This implies that $\;g_{i_{\tau}}
\in g_{i_{\nu}}\overline{\sigma}{\rm H}.\;$ In 1 it was shown that
the coset $\;g_{i_{\nu}}\overline{\sigma}{\rm H} \;$ coincides
with one of the cosets (\ref{bsretatttadf}). Therefore,
$\;g_{i_{\tau}}\;$ belongs to one of the cosets
(\ref{bsretatttadf}). Contradiction, since by assumption $\; \tau
> \nu, \;$ and therefore, the element $\;g_{i_{\tau}}\;$ does not
belong to any of the cosets (\ref{bsretatttadf}). \vspace{0.5cm}

\hspace{1.3cm} 3. \hspace{0.5cm} Let $\; S_{1}, \;$ $\;S_{2} \;$
and $\;S_{3} \;$ be three sets of elements of $\;G\;$ of the form
(\ref{nujjjhyghbn}). Assume that $\; S_{1}  \sim S_{2} \;$ and
also that $\; S_{2}  \sim S_{3}. \;$ The relation defined in Eq.
(\ref{vbcncnmx}) is reflexive ($\; S_{1} \sim S_{1} \;$),
symmetric ($\; S_{1} \sim S_{2} \;$ implies $\; S_{2} \sim S_{1}
\;$), and transitive ($\; S_{1} \sim S_{2} \;$ and  $\; S_{2} \sim
S_{3} \;$ together, imply $\; S_{1} \sim S_{3} \;$). This
completes the proof.
\end{prf}

\noindent Let $\; \overline{\mathcal S} \;$ be the subset of $\;
\mathcal S \;$ which has only members of the form
(\ref{nujjjhyghbn}). The relation defined in Eq. (\ref{vbcncnmx})
is an equivalence relation on $\; \overline{\mathcal S}. \;$ If
$\; S \in \overline{\mathcal S},  \;$  then $\; \tilde {S} \;$
denotes the equivalence class of $\; S. \;$

%\vspace {0.5cm}

%The group K has three subgroups of order 2, namely the following,
%${\rm K}_{1}=\; \left \{ ({\rm I},{\rm I}),(-{\rm I},-{\rm I})
%\right \},\;$ ${\rm K}_{2}=\; \left \{ ({\rm I},{\rm I}),({\rm
%I},-{\rm I}) \right \},\;$ and, ${\rm K}_{3}=\; \left \{ ({\rm
%I},{\rm I}),(-{\rm I},{\rm I}) \right \}.\;$ If $\; S\;$ is a
%finite subgroup of $\; K=SO(2) \times SO(2) \;$ then $\; S\;$ is a
%subgroup of $\;{\rm C}_{\rm n} \times {\rm C}_{\rm m} \;$. If K or
%$\;{\rm K}_{\rm i} \ , \ {\rm i}=1,2,3\;$ is a subgroup of
%$\;S\;$, then it is a normal subgroup of $\; S \; $.
%\ Here we consider the remaining cases; the closed subgroups $H$
%of $K$ with dimension $0$. \ As shown in Section 6, every such
%subgroup is a finite subgroup $ F$ of $C_{n}\times C_{m}$:
%\begin{equation}
%H=F\subset C_{n}\times C_{m}.
%end{equation}
%In order to find the invariant subspaces $\mathcal{H}(F)$, it is
%convenient first to prove some general results about group
%actions. \ It is also convenient, in this Section, to change
%notation; instead of $G$ denoting $ {\rm S}{\rm L}(2,R)$, $G$ will
%here denote a finite group.

In our study the finite group $\;G\;$ is the group $\; C_{{\rm n}}
\times C_{{\rm  m}} \;$ which acts on the torus $\;
%{\rm T}^{2}
\mathcal P
\simeq P_{1}(R) \times P_{1}(R).\;$  With a view to apply the
previous theory to the problem under consideration we leave the
group $\; G \;$ to act on any manifold ( the theorems proven here
are slightly more general than it is strictly needed, in fact we
allow $\; G \;$ to act on any topological space) and we establish
a few more facts. The details are as follows.

%\vspace{0.5cm}

 Let $M$ be any topological space, and let $G$ be any
finite group which acts on $M$ from the right. \ That is to say,
we are given a map $M\times G\rightarrow M$, denoted
$(x,g)\longmapsto xg$, with the following properties. \ For each
$g\in G$, the map $x\longmapsto xg$ is a homeomorphism of $M$ onto
itself. \ If $g=e$ (the identity element), $xe=x$ for every $x\in
M$. \ For any $x\in M$ and $g_{1},g_{2}\in G$, $
(xg_{1})g_{2}=x(g_{1}g_{2})$. There is
 an homomorphism $\; h_{G} \;$ from $\; G \;$ into the group of
automorphisms $\; {\rm Aut}(M) \;$ of $\; M \;$ which is naturally
associated with this action. Let $\; {\rm K} \;$ be the kernel of
$\; h_{G} \;$. The action of $\;G\;$ on  $\;M\;$ passes naturally
to an action of $\;G/{\rm K}\;$ on $\;M.\;$ Henceforth, we will
assume that this action
$$
x(g{\rm K})=xg,
$$
where, $\;x \in M, \;$ and $\;g \in G \;$ is fixed point free.
 %We assume that $G$ acts
%effectively on $M$; that is, $xg=x$ for all $x\in M$ implies
%$g=e$.
\ An \it{elementary domain} \normalfont for the given action is an
open subset ${\rm E}\subset M$ such that the following conditions
are satisfied:
%\underset
\begin{equation}
\begin{array}{l}
{\rm (A) \ For\; any \;}g_{1},g_{2}\in G,\; \; {\rm with }\;
g_{1}\neq g_{2}k,\; \;k \in {\rm K},\; \;{\rm E}g_{1}\cap
\overline{{\rm E}}g_{2}=\varnothing , \\
\label{mkkkliojhgt} {\rm (B) \;} {{\large \bigcup }}
_{_{_{_{\!\!\, \;\!\!\!\!\!\!\!\!\!\large g \in
G}}}}\!\!\overline{{\rm E}}g=M.
%\quad ({\rm disjoint \quad union}).
%({\rm disjoint \quad union}).
%{\rm (B)} \quad {{\large g\in G}}{{\Large \cup }}\overline{E}g=M.
\end{array}
\end{equation}
Here bar means topological closure, and $\varnothing $ means the empty
set.  Now let $\;{\rm  H}\subset G\;$ be any subgroup of $\;G,\;$
and let $\;S \subset G \;$ be a set of representatives of the left
cosets of $\;{\rm H}\;$ which is a member of $\;
\overline{\mathcal S}. \;$ Since $\; S \subset G \;$ is a set of
representatives of the left cosets of $\; {\rm H} \;$ in the coset
space $\; G/ {\rm H} \;$ we have
%\underset
\begin{equation}
\begin{array}{l}
{\rm (C) \; }G=S{\rm H}, \\
{\rm (D) \; }S{\rm H}=~
%\underset
\!\!^{{\large \bigcup }} _{_{_{\!\ \!\!\!\large s \in
S}}}\!\!s{\rm H}\quad ({\rm disjoint \quad union}).
\end{array}
\end{equation}
Then we have the following result relating elementary domains for
$G$ and $ {\rm H}:$
%\footnote{By using the results of sections one could try and find `` nice-
%looking '' connected elementary regions for the finite actual little groups .
%}:
%\textbf{8.1}:

\begin{prop}
\label{nncnncncnmjdk} ${\rm E}S$  is an elementary domain for $\;
{\rm H}. \;$
\end{prop}
\begin{prf}
\normalfont \ First note that condition (B) can be written
$\overline{{\rm E}} G=M$. \ To verify condition (B) for ${\rm H}$,
we calculate, using (C), as follows. \ $(\overline{{\rm E}}S){\rm
H}= \overline{{\rm E}}(S{\rm H})=\overline{{\rm E}}G=M,$ as
required. \ To verify condition (A) for ${\rm H}$, we must prove
that, for any $ h_{1},h_{2}\in {\rm H}$ with $h_{1}\neq
h_{2}k_{{\rm H}}$,$\;k_{{\rm H}} \in {\rm K} \cap {\rm H}, \;$
$({\rm E}S)h_{1}$ and $(\overline{{\rm E}S})h_{2}=(\overline{{\rm
E}}S)h_{2}$ are disjoint. \ Assume that they are not disjoint. \
Then there exists an element $x\in ({\rm E}S)h_{1}\cap
(\overline{{\rm E}}S)h_{2}$. \ So, by definition, $
x=z_{1}s_{1}h_{1}=z_{2}s_{2}h_{2}$ for some $z_{1}\in {\rm E}$,
$z_{2}\in \overline{{\rm E}}$ and $s_{1},s_{2}\in S$. \ Write
$g_{1}=s_{1}h_{1}$ and $ g_{2}=s_{2}h_{2}$. \ Then
$z_{1}g_{1}=z_{2}g_{2}$, so ${\rm E}g_{1}\cap \overline{{\rm
E}}g_{2}\neq \varnothing $. \ So, by property (A), $\;g_{1}=g_{2}k,\;$
$\;k \in {\rm K}.\;$
 \ That is,
 \be \label{mnbghjkiu} s_{1}h_{1}=s_{2}h_{2}k. \ee
We distinguish two cases.

\vspace{0.3cm}

\hspace{0.3cm} 1. \hspace{0.3cm} Assume that $\; {\rm H} \cap {\rm
K}={\rm K}. \;$ Then Eq. (\ref{mnbghjkiu}) gives
$$
s_{1}{\rm H}=s_{2}{\rm H}\Longrightarrow s_{1}=s_{2}.
$$
Substituting back into Eq. (\ref{mnbghjkiu}) we obtain $\;
h_{1}=h_{2}k, \;$ where $\; k \in {\rm H} \cap {\rm K}={\rm K}.
\;$ \ This contradicts $\;h_{1}\neq h_{2}k\;$, and so $({\rm
E}S)h_{1}$ and $(\overline{{\rm E}}S)h_{2}$ are actually disjoint;
condition (A) is satisfied for ${\rm H}$. Henceforth we can assume
that $\; {\rm H} \cap {\rm K} < {\rm K}, \;$(the symbol $\; < \;$
indicates proper subgroup), and that the coset space $\;{\rm
K}/{\rm H} \cap {\rm K} \;$ has at least two elements.

\vspace{0.3cm}

\hspace{0.3cm} 2. \hspace{0.3cm}Assume that $\; {\rm H} \cap {\rm
K} < {\rm K}. \;$ Eq. (\ref{mnbghjkiu}) implies that \be
\label{mcmcmmcdse} s_{1}h_{1}k^{-1}=s_{2}h_{2}.\ee Now $\;k^{-1}
\in {\rm K} \;$ and $\;{\rm K} \;$ is a normal subgroup of $\;G.
\;$ Therefore, $\;h_{1}k^{-1}=k'h_{1}, \;$ for some $\; k' \in
{\rm K}. \;$ Substituting into Eq. (\ref{mcmcmmcdse}) we obtain
\be \label{ncdgshe} s_{1}k'h_{1}=s_{2}h_{2}. \ee We distinguish
three cases. \vspace{0.15cm}

\hspace{0.7cm} 2a. \hspace{0.3cm} Assume now that $\; s_{1}k' \in
S, \;$ where $\; S \;$ is the fixed  set of coset representatives
we chose at the beginning. Then Eq. (\ref{ncdgshe}) gives \be
s_{1}k'{\rm H}=s_{2}{\rm H} \Longrightarrow s_{1}k'=s_{2}. \ee
Substituting into Eq. (\ref{ncdgshe})  we obtain $\;h_{1}=h_{2}
\;$. This contradicts $\;h_{1}\neq h_{2}k,\;$ and so henceforth we
can assume that $\; s_{1}k' \; / \!\!\!\!\! \in S.\;$
\vspace{0.15cm}

\hspace{0.7cm} 2b. \hspace{0.3cm} Assume now that $\; s_{1}k' \; /
\!\!\!\!\!\! \in S \;$ and also that $\; s_{1}k' \in s_{1}{\rm H}.
\;$ Then there exists $\; h \in {\rm H} \;$ such that
$\;s_{1}k'=s_{1}h. \;$Therefore, $\;k'=h \;$ and $\;k' \in {\rm K}
\cap {\rm H}. \;$According to Proposition \ref{byhjujjkiiklomnd}
the group $\;{\rm K} \cap {\rm H} \;$ is normal in $\; {\rm H},
\;$ and therefore, $\;k'h_{1}=h_{1} \overline {k}, \;$ for some
$\; \overline {k} \;$ in $\;{\rm K} \cap {\rm H}. \;$ Substituting
back into Eq. (\ref{ncdgshe}) we obtain \be \label{mkjuhygtfr}
s_{1}h_{1}=s_{2}h_{2}{\overline{k}}^{\;\;-1}. \ee From the last
Equation we obtain
$$ s_{1}{\rm H}=s_{2}{\rm
H} \Longrightarrow s_{1}=s_{2}. $$ Substituting back into Eq.
(\ref{mkjuhygtfr}) we obtain $$
h_{1}=h_{2}{\overline{k}}^{\;\;-1}. $$ This contradicts
$\;h_{1}\neq h_{2}k,\;$ and so henceforth we can assume that $\;
s_{1}k' \; / \!\!\!\!\! \in S,\;$ and that $\; s_{1}k' \; /
\!\!\!\!\! \in s_{1}{\rm H}.\;$ \vspace{0.15cm}

\hspace{0.9cm} 2c. \hspace{0.3cm} Assume now that $\; s_{1}k' \; /
\!\!\!\!\!\! \in S \;$ and also that $\; s_{1}k' \; / \!\!\!\!\!\!
\in s_{1}{\rm H}. \;$  From Eq. (\ref{ncdgshe}) we obtain \be
\label{kolllkkklkio} s_{1}k' \in s_{2}{\rm H}. \ee From
Proposition \ref{nujkilokmnjhyu} (statement 1) we have that \be
\label{njuhyghjkmnb} s_{1}k' \in s_{1}k_{{\rm o}}{\rm H}, \ee for
some (unique) $\;s_{1}k_{{\rm o}} \in S. \;$ From Equations
(\ref{kolllkkklkio}) and (\ref{njuhyghjkmnb}) we have $\;
s_{2}{\rm H}=s_{1}k_{{\rm o}}{\rm H} \;$ and therefore we obtain
$\; s_{2}=s_{1}k_{{\rm o}}. \;$ Substituting back into Eq.
(\ref{ncdgshe}) we obtain
$$
 k'h_{1}=k_{{\rm o}}h_{2}.
$$
The last Equation gives $\; k_{{\rm o}}^{-1}k'=h_{2}h_{1}^{-1}.
\;$ Therefore we have $\;k_{{\rm o}}^{-1}k'=k'' \in {\rm K} \cap
{\rm H}. \;$ So we have $\; k''h_{1}=h_{2}. \;$ The group $\; {\rm
K} \cap {\rm H} \;$ is normal in $\; {\rm H} \;$ and therefore $\;
k''h_{1}=h_{1}\overline{\overline{k}}, \;$ for some
$\;\overline{\overline{k}} \in {\rm K} \cap {\rm H}. \;$ So
finally we obtain
$$
h_{1}=h_{2} \overline{\overline{k}}^{\;\;-1}.
$$
\ This contradicts $\;h_{1}\neq h_{2}k\;$, $\; k \in {\rm K} \cap
{\rm H}. \;$ So $({\rm E}S)h_{1}$ and $(\overline{{\rm E}}S)h_{2}$
are actually disjoint; condition (A) is satisfied for $\;{\rm H}
\;$. This completes the proof.
\end{prf}

% \ So $s_{1}H\cap s_{2}H\neq
%\varphi $. \ So, by (D), $s_{1}=s_{2}=s$, say. \ This gives
%$sh_{1}=sh_{2}$ and so, cancelling the $s$ , $h_{1}=h_{2}$. \ This
%contradicts $h_{1}\neq h_{2}$, and so $({\rm E}S)h_{1}$ and
%$(\overline{{\rm E}}S)h_{2}$ are actually disjoint; condition (A)
%is satisfied for $H$.
\section{Invariant Subspaces for Finite Groups}
\label{s2}
For reasons which will become clear later, in this section we take
the positive integers $\;n \;  {\rm and} \; m \;$ to be even.
Consider the particular (non$-$effective) action $\; S^{1}\times
C_{n}\rightarrow S^{1}\;$ given by \be\label{njukiloha}(\rho
,g_{r})\longmapsto \rho g_{r}=\rho +\frac{4\pi }{n}r,\ee where
$\rho $ is the usual angular coordinate on the circle, taken mod
$2\pi $, and $g_{r}\in C_{n}$ are the elements of $C_{n}$, where $
\; 0\leq r\leq (n-1) \;$. Let $\; h' \;$ be the homomorphism
associated with the action (\ref{njukiloha}).
%The homomorphism $\;
%h' \;$ has domain the group $\; C_{n} \;$ and range the group of
%automorphisms of the circle $\; S^{1}. \;$
One can easily show that the kernel $\;{\rm K}_{C_{n}} \;$ of this
homomorphism is $\;{\rm K}_{C_{n}} = \{ {\rm I}, -{\rm I} \},\;$
where $\; {\rm I} \;$ is the identity element of the group $\;
C_{n}. \;$ \ Now we  show that
\begin{prop}

The open set
\begin{equation}
{\rm E}_{n}=\{\rho \in S^{1}~\left| ~0<\rho <4\pi /n\right. \}
\end{equation}
is elementary for the action (\ref{njukiloha}).
\end{prop}
\begin{prf}
\normalfont Let $\;g_{1} \;$ and $\; g_{2} \;$ be two elements of
$\; C_{n}. \;$ Then we have \begin{eqnarray} \label{nbhdhhhdhdbcs}
{\rm E}_{n}g_{1} & = & \left \{\rho \in S^{1}~\left| \rho=\rho'+
\frac{4\pi}{n}r_{1} \; , \; \rho' \in {\rm E}_{n}
% ~\frac{4\pi}{n}r_{1}<\rho
%+ \frac{4\pi}{n}r_{1}
%<\frac {4\pi} {n} + \frac{4\pi}{n}r_{1}
\right. \right  \}
%\nonumber
\\
%{\rm E}_{n}g_{1}
%\label{nbhdhhhdhdbcs}
& = & \left \{\rho \in S^{1}~\left| ~\frac{4\pi}{n}r_{1}<\rho
%+ \frac{4\pi}{n}r_{1}
<\frac {4\pi} {n} + \frac{4\pi}{n}r_{1}\right. \right  \}
\nonumber
\end{eqnarray} and
 \begin{eqnarray} \label{nncnncsaer}
\overline {{\rm E}}_{n}g_{2} & = & \left \{\rho \in S^{1}~\left|
\rho=\rho'+ \frac{4\pi}{n}r_{2} \; , \; \rho' \in \overline {{\rm
E}}_{n}
%~\frac{4\pi}{n}r_{2} \leq \rho
%+ \frac{4\pi}{n}r_{2}
%\leq \frac {4\pi} {n} + \frac{4\pi}{n}r_{2}
\right. \right  \}
\\
% \overline {{\rm
%E}}_{n}g_{2}
& = &\left \{\rho \in S^{1}~\left| ~\frac{4\pi}{n}r_{2} \leq \rho
%+ \frac{4\pi}{n}r_{2}
\leq \frac {4\pi} {n} + \frac{4\pi}{n}r_{2}\right. \right  \}.
\nonumber
\end{eqnarray} Assume now that $\; {\rm E}_{n}g_{1} \cap \overline {{\rm
E}}_{n}g_{2} \neq \varnothing. \;$ Then from Eqs.
(\ref{nbhdhhhdhdbcs}) and (\ref{nncnncsaer}) we conclude that
there must exist some
%$\; \rho, \;$ say
$\; \rho_{{\rm o}} \in {\rm E}_{n} \;$such that \be
\label{zbbabaghjjhaa} \rho_{{\rm o}} + \frac{4\pi}{n}r_{1}=
\rho_{{\rm o}} + \frac{4\pi}{n}r_{2} + a2\pi,\ee for some integer
$\; a, \;$  since $\; \rho \;$ is only defined mod $\; 2\pi.
\;$Eq. (\ref{zbbabaghjjhaa}) gives \be r_{1}-r_{2}=a \frac{n}{2}.
\ee The last Equation is equivalent to \be g_{1}=g_{2}k, \ee where
$\; k \in {\rm K}_{C_{n}}, \;$ and $\;{\rm K}_{C_{n}}\;$ is the
kernel of the homomorphism $\; h'. \;$ Therefore the condition $\;
{\rm (A) } \;$  of Eq. (\ref{mkkkliojhgt}) is satisfied. We note
that that {\it every} \normalfont $\; \rho_{{\rm o}} \in S^{1} \;$
belongs to some interval of the form
 $$  \overline {{\rm E}}_{n}g_{r}=\left \{\rho \in S^{1}~\left|
~\frac{4\pi}{n}r \leq \rho
%+ \frac{4\pi}{n}r
\leq \frac {4\pi} {n} + \frac{4\pi}{n}r\right. \right  \}, $$
where, $\; r \in \{0,1,2,...,\frac{n}{2}-1 \}. \;$ We conclude
that \be \label{nnnmcmmckddsl}
%{\rm (B) \;}
{{\large \bigcup }} _{_{_{_{\!\!\, \;\!\!\!\!\!\!\!\!\!\large
g_{r} \in C_{n} }}}}\!\!\overline{{\rm E}}_{n}g_{r}=S^{1}. \ee (In
fact Eq. (\ref{nnnmcmmckddsl}) is satisfied even when $\; r \;$
runs only through the set of values $\; \{ 0,1,2,...,\frac{n}{2}-1
\}). \;$ Therefore the condition $\; {\rm (B) } \;$  of Eq.
(\ref{mkkkliojhgt}) is satisfied. This completes the proof.
\end{prf}
 \ Next, consider
the (non$-$effective) action $
\mathcal P
%{\rm T}^{2}
\times (C_{n}\times
C_{m})\rightarrow
\mathcal P
%{\rm T}^{2}
$ given by
%\tag*{(8.1)}\begin{equation}
\begin{equation}
\label{nnnnmnbbxxxcz} ((\rho ,\sigma ),(g_{i},g_{j}))\longmapsto
(\rho g_{i},\sigma g_{j})
\end{equation}
where $0\leq i\leq (n-1)$, $0\leq j\leq (m-1)$, $g_{i}\in C_{n}$,
$g_{j}\in C_{m}$ and \be \label{bhnjmkloiulkjuhygtfr} \rho
g_{i}=\rho +\frac{4\pi }{n}i, \quad \quad \sigma g_{j}=\sigma
+\frac{ 4\pi }{m}j. \ee \ Define the set ${\rm F}_{nm}\subset
\mathcal P
%{\rm
%T}^{2}
$ by the formula
\begin{equation}
{\rm F}_{nm}={\rm E}_{n}\times {\rm E}_{m}\subset
\mathcal P.
%{\rm T}^{2}.
\end{equation}
We prove now that

\begin{prop}
\label{mnjhkloiujk} ${\rm F}_{nm}$ is an elementary domain for the
action (\ref{nnnnmnbbxxxcz}).
\end{prop}
\begin{prf}
\normalfont \ Since $\;{\rm E}_{n}$ and ${\rm E}_{m}\;$ are open,
so is $F_{nm}$.  If $\;(g_{i},g_{j})\neq (g_{i^{\prime
}},g_{j^{\prime }})k,\;$ where $\;k \in {\rm K}=\{({\rm I},{\rm
I}), (-{\rm I},-{\rm I}), ({\rm I},-{\rm I}),(-{\rm I},{\rm I})\},
\;$ then either $\; g_i \neq g_{i^{\prime }} \lambda \;$ or $\;
g_j \neq g_{j^{\prime }} \lambda \;$ or both, where $\; \lambda
\in {\rm K}_{C_{n}}= \{{\rm I},-{\rm I} \}. \;$ \ We have
\begin{equation}
{\rm F}_{nm}(g_{i},g_{j})=({\rm E}_{n}g_{i})\times ({\rm
E}_{m}g_{j}),
\end{equation}
\begin{equation}
\overline{{\rm F}}_{nm}(g_{i^{\prime }},g_{j^{\prime
}})=(\overline{{\rm E}} _{n}g_{i^{\prime }})\times (\overline{{\rm
E}}_{m}g_{j^{\prime }}).
\end{equation}
Taking the intersection gives
\begin{equation}
({\rm E}_{n}g_{i}\times {\rm E}_{m}g_{j})\cap (\overline{{\rm
E}}_{n}g_{i^{\prime }}\times \overline{{\rm E}}_{m}g_{j^{\prime
}})=({\rm E}_{n}g_{i}\cap \overline{{\rm E}} _{n}g_{i^{\prime
}})\times ({\rm E}_{m}g_{j}\cap \overline{{\rm E}}_{m}g_{j^{\prime
}}).
\end{equation}
Since ${\rm E}_{n}$ is elementary for the original action of
$C_{n}$ on $S^{1}$, at least one of the factors on the RHS is
empty. \ So condition (A) is satisfied. \ For (B), let $(\rho
,\sigma )\in {\rm T}^{2}$ be given. \ Then, because ${\rm E}_{n}$
is elementary for the original action, $\rho \in \overline{ {\rm
E}}_{n}g_{i}$ for some $i=0,1,2,...,(n-1)$, and $\sigma \in
\overline{{\rm E}} _{m}g_{j}$ for some $j=0,1,2,...,(m-1)$. \ So
\begin{equation}
(\rho ,\sigma )\in \overline{{\rm E}}_{n}g_{i}\times
\overline{{\rm E}}_{m}g_{j}= \overline{{\rm F}}_{nm}(g_{i},g_{j}).
\end{equation}
That is, every $(\rho ,\sigma )$ belongs to some $\overline{{\rm
F}} _{nm}(g_{i},g_{j})$ and so condition (B) is satisfied. This
completes the proof.
\end{prf}

 Now let $F\subset C_{n}\times C_{m}$ be any
subgroup of $C_{n}\times C_{m}$.  Let $\; S  \; $ be a selection
of representatives of left cosets of $\; F \; $ which belongs to
$\; \overline{S}. \; $  Then, by Proposition \ref{nncnncncnmjdk},
the set
\begin{equation}
{\rm E}={\rm F}_{nm}S\subset \mathcal P
%{\rm T}^{2}
\end{equation}
is an elementary domain for the subgroup $\; F. \; $ Let $\;
\mathcal{H}(
\mathcal P
%{\rm T}^{2}
)\; $ denote the Hilbert space of functions
$f:
%{\rm T}^{2}
\mathcal P
\rightarrow R$ which are square integrable with
respect to the usual Lesbegue measure $d\theta \wedge d\phi $ on
the torus $\;
\mathcal P
%{\rm T}^{2}
\simeq P_{1}(R) \times P_{1}(R). \;$ \ We
now want to find the Hilbert space $\mathcal{H}(F)$ of all
invariant vectors for each of the finite subgroups $F$ of
$C_{n}\times C_{m}$. Thus, $\mathcal{H}(F)$ is the space
%\ So definition
%\ref{dskjnbpotfcnkm} becomes, for
\begin{equation}
\mathcal{H(}F)=\left\{ \widetilde{\zeta }\in \mathcal{H}(
\mathcal P
%{\rm T}^{2}
)\left | ~T^{\prime }(h)\widetilde{\zeta }=\widetilde{\zeta
}\quad {\rm for} \quad {\rm for} \; {\rm all}\;h \in F\right.
\right\}. \end{equation} $\mathcal{H}(F)$ is a closed subspace of
 $\; \mathcal{H}(
\mathcal P
% {\rm T}^{2}
 ).\; $
 First note that ${\rm E}={\rm F}_{nm}S$ is a  union of open
rectangles in $
\mathcal P
%{\rm T}^{2}
$, so inherits the Lesbegue measure
$d\theta \wedge d\phi $ from $
\mathcal P
%{\rm T}^{2}
$. \ Let $\mathcal{H}(\rm
E)$ denote the Hilbert space of square integrable functions
$f:{\rm E}\rightarrow R$, and $\mathcal{H}_{\rm E}(
\mathcal P
%{\rm T}^{2}
)$
the Hilbert subspace of $\mathcal{H}(
\mathcal P
%{\rm T}^{2}
)$ consisting of
functions which vanish outside ${\rm E}$;
\begin{equation}
\mathcal{H}_{\rm E} (
\mathcal P
%{\rm T}^{2}
)=\left\{ f\in \mathcal{H}(
\mathcal P
%{\rm
%T}^{2}
)~\left| ~f(x)=0,\quad {\rm all } \;\;x\notin E\right.
\right\} .
\end{equation}
 There is a bijection between $\mathcal{H}({\rm E})$ and
$\mathcal{H}_{\rm E}(
\mathcal P
%{\rm T}^{2}
)$, obtained as follows.  Define
maps $\alpha: \mathcal{H}({\rm E})\rightarrow\mathcal{H}_{\rm
E}(
\mathcal P
%{\rm T}^{2}
)$ and $\beta:\mathcal{H}_{\rm E}(
\mathcal P
%{\rm T}^{2}
)
\rightarrow \mathcal{H}(\rm E)$ by ``extending to $0$ outside
${\rm E}$'' and ``restricting to ${\rm E}"$ respectively:
%0,\text{ \ \ \ \ \ }x\notin E
\begin{equation}
(\alpha (l))(x)=\left\{
\begin{array}{l}
l(x),\quad x \in {\rm E}\\
0,\quad \quad \, x \notin {\rm E}
\end{array}
\right\} ,
\end{equation}
\begin{equation}
(\beta (f))(x)=f(x),\quad x\in {\rm E}.
\end{equation}
Here typical elements are denoted by $l\in \mathcal{H}({\rm E})$
and $f\in \mathcal{H}_{\rm E} (
\mathcal P
%{\rm T}^{2}
)$.  Then routine
checks show that $(\beta \alpha )(l)=l$ for all $l$, and $(\alpha
\beta )(f)=f$ for all $f$; $\beta \alpha $ is the identity map on
$\mathcal{H}({\rm E})$, $\alpha \beta $ is the identity map on
$\mathcal{H}_{\rm E} (
\mathcal P
%{\rm T}^{2}
)$. Thus there is an induced
bijection
%\tag*{(8.2)}
\begin{equation}
\label{jnmnhbgsdnmjk} \mathcal{H}({\rm E})\leftrightarrow
\mathcal{H}_{\rm E} (
\mathcal P
%{\rm T}^{2}
).
\end{equation}
%M=\underset{{\large g\in G}}{{\Large \cup }}\overline{E}g
%%%%%%%%%%%%%%%%%%%%%%%%%%%%%%%%%%%%%%%%%%%%%%%%%%%%%%%%%%%%%%%%%%%%%%%%%%%%%%%
Before our next result, note first that, from the conditions (A)
and (B) defining an elementary region, the union
%\tag*{(8.3)}
\begin{equation}
\label{bvcdfghnmklcdxsza}
%{\rm (B) \;} {{\large \bigcup }}
%_{_{_{_{\!\!\, \;\!\!\!\!\!\!\!\!\!\large g \in G}}}}\!\!
%\overline{{\rm E}}g=M
%\quad ({\rm disjoint \quad union}).
%M={{\large g\in G}}{{\Large \cup }}\overline{\rm E}g
M={{\large \bigcup }} _{_{_{_{\!\!\,
\;\!\!\!\!\!\!\!\!\!\!\!\!\!\!\large g{\rm K} \in G / {\rm
K}}}}}\!\!
%{{\large g\in G}}{{\Large \cup }}
\overline{\rm E}(g{\rm K})
\end{equation}
is ``almost disjoint'', in the sense that overlaps between the
sets $ \overline{\rm E}g_{1}$ and $\overline{\rm E}g_{2}$, for any
$\;g_{1} \neq g_{2}k, \quad k \in {\rm K}, \;$ can only involve
boundary points of ${\rm E}$. \ In our situation ($M=
\mathcal P
%{\rm T}^{2}
$,
$F\subset G=C_{n}\times C_{m} \quad {\rm and} \quad {\rm
K}=\{({\rm I},{\rm I}),({\rm I},-{\rm I}),(-{\rm I},{\rm
I}),(-{\rm I},-{\rm I}) \} )$, these sets (involving boundary
points) are, at most, one dimensional, and so certainly of measure
zero.
%\ We will write $ G$ instead of $F$ for convenience.
%\ So definition
%\ref{dskjnbpotfcnkm} becomes, for
%\begin{equation}
%\mathcal{H(}G)=\left\{ \widetilde{\zeta }\in \mathcal{H}({\rm
%T}^{2})\left | ~T^{\prime }(h)\widetilde{\zeta }=\widetilde{\zeta
%}\quad {\rm for} \quad {\rm for} \; {\rm all}\;g \in G\right.
%\right\} .
%\end{equation}

%In our situation
%It was shown in chapter
%\ref{lyglygyhgyuhglyughbyhug} that
The representation theory of $\;\mathcal H \mathcal B
%_{c}
\;$ is
governed by the dual action $T'$ of $\mathcal G={\rm S}{\rm
L}(2,R) \times {\rm S}{\rm L}(2,R) \;$
on
%$\phi_{\rm{o}} \in  {L^{2}}^{'}(\mathcal P,\lambda ,R)$
$
{L^{2}}^{'}(\mathcal P,\lambda ,R)
%A'(\mathcal N)
\simeq
{L^{2}}(\mathcal P,\lambda ,R)
%A(\mathcal N)
$ given by (Eq. (\ref{tzatziki}))
$$ \label{ulbkyuvrdfcrfxcr}
(T'(g,h)\phi)(x,y)=k_{g}^{-3}(x)s_{g}(x)k_{h}^{-3}(y)s_{h}(y)\phi(xg,yh).
$$
%With a view to apply the theory we develop in our situation
Subsequently, for notational simplicity,  we will write $\;
\label{mkiuyhhhnbgdfdffdr} (T'(g)f)(x)= \gamma(x,g)f(xg), \;$ and
it will be understood that
%in our situation
by  $\;\gamma(x,g)\;$ we denote the multiplier \be
\label{vdsermklokmnj}
\gamma(x,g)=k_{g}^{-3}(x)s_{g}(x)k_{h}^{-3}(y)s_{h}(y), \ee and
that when we write $\; (T'(g)f)(x)= \gamma(x,g)f(xg), \;$ by $\; x
\;$ we denote a point of the torus
$
\mathcal P,
%\; P_{1}(R) \times P_{1}(R)
\;$ and by $\; g \;$ an element of
$ \mathcal G.
%\; {\rm S}{\rm L}(2,R) \times
%{\rm S}{\rm L}(2,R)
\;$

Now we
define two maps $\; \rho \;$ and $\; \sigma \;$ (we use the same
symbols to denote the usual angular coordinates on the torus
$
\mathcal P;
%P_{1}(R) \times P_{1}(R)
$
hopefully the context will cause no
doubts regarding their meaning when they are encountered) as
follows
\begin{eqnarray}
\sigma & : & \mathcal{H}_{{\rm E}}(
\mathcal P
%{\rm T}^{2}
)\rightarrow
\mathcal{H}(F) \nonumber \\
\label{mnjzzzzxcdsefs} \sigma (f)& = & \frac {1}{|{\rm
K}_{F}|}\sum_{{\large g\in F}}T^{\prime }(g)f,
\end{eqnarray}
where $\; f \;$ is an element of $\; \mathcal{H}_{{\rm E}}(
\mathcal P
%{\rm
%T}^{2}
), \;$ and $\;|{\rm K}_{F}|\;$ denotes the order of the
group $\;{\rm K}_{F}=F  \cap {\rm K}. \;$ The map $\; \rho \;$ is
defined as follows
\begin{eqnarray}
\rho & : & \mathcal{H}(F)\rightarrow \mathcal{H} _{{\rm E}}(
\mathcal P
%{\rm
%T}^{2}
) \nonumber \\
\label{bbxnnxxnhhskkaju} (\rho (f))(x) & = & \chi _{{\rm
E}}(x)f(x)=(\chi _{{\rm E}}\cdot f)(x),
\end{eqnarray}
 where $\; x \in
 \mathcal P,
% {\rm T}^{2}
  \;$ $\; f \in \mathcal{H}(F), \;$
 and $\chi _{{\rm E}}$
denotes the characteristic function of ${\rm E}$, equal to $1$
inside $\;{\rm  E}, \;$ and $\;0\;$ outside.  The dot denotes
pointwise multiplication. We now prove the following theorem

 %have theF%\textbf{8.1}:
\begin{thrm}
\label{zaaasdhjjnmkijk} Let $\; G= C_{n} \times C_{m} \;$ and let
$\; F \;$ be a subgroup of $\; G. $ Let $\; h_{G} \; : \; G
\mapsto {\rm Aut}(
\mathcal P),
%{\rm T}^{2}),
\;$ $\;{\rm Aut}(
\mathcal P
%{\rm T}^{2}
)\;$ is
the group of automorphisms of $\;
\mathcal P,
%{\rm T}^{2},
\;$ be the
homomorphism which is associated with the action
(\ref{nnnnmnbbxxxcz}). Let $\; {\rm K} \;$ be the kernel of the
homomorphism $\; h_{G}. \;$ Let $\; {\rm K}_{F} \;$ denote the
intersection $\; {\rm K} \cap F. \;$ Then we have the following
\begin{enumerate}
\item{ The map $\; \rho \sigma \;$, computed on $\;
\mathcal{H}_{{\rm E}}(
\mathcal P
%{\rm T}^{2}
), \;$ gives \be
\label{jilokokiiiju} (\rho (\sigma (f)))(x)= \frac {1}{|{\rm
K}_{F}|} \left [ \gamma(x,e) + \gamma(x,k_{2}) + \gamma(x,k_{3}) +
... + \gamma(x,k_{|{\rm K}_{F}|}) \right ]f(x),
 \ee
where, $\; x \in
\mathcal P,
%{\rm T}^{2}
\;$ $\; f \in \mathcal{H}_{{\rm
E}}(
\mathcal P
%{\rm T}^{2}
), \;$ $\;|{\rm K}_{F}|\;$ denotes the order of
the group $\; {\rm K} \cap F \;$ and $\; {\rm K} \cap F= \left \{
e, k_{2}, k_{3}, ... , k_{|{\rm K}_{F}|} \right \} \;$.} \item{
The map $\; \sigma \rho, \;$ computed on $\; \mathcal{H}(F), \;$
%_{{\rm E}}({\rm T}^{2} ), \;$
gives \be  \label{nmnmnnnmmmbbbbvabcd}(\sigma (\rho (f)))(x)=
f(x),
 \ee
where, $\; x \in
\mathcal P
%{\rm T}^{2}
\;$ and  $\; f \in \mathcal{H}(F).
\;$ }
\end{enumerate}

%$\; {\rm K} \;$ be the kernel of the homomorphism

%\ Next, consider the (non-effective) action $ {\rm T}^{2}\times
%(C_{n}\times C_{m})\rightarrow {\rm T}^{2}$ given by
%\tag*{(8.1)}\begin{equation}%
%\begin{equation}
%\label{nnnnmnbbxxxcz} ((\rho ,\sigma ),(g_{i},g_{j}))\longmapsto
%(\rho g_{i},\sigma g_{j})
%\end{equation}
%where $0\leq i\leq (n-1)$, $0\leq j\leq (m-1)$, $g_{i}\in C_{n}$,
%$g_{j}\in C_{m}$ and $\rho g_{i}=\rho +\frac{4\pi }{n}i$, $\sigma
%g_{j}=\sigma +\frac{ 4\pi }{m}j$.

%There is a bijection $\mathcal{H}_{\rm E}({\rm T}^{2})
%\leftrightarrow \mathcal{H}(G).$
\end{thrm}

\begin{prf}
\normalfont
%\ We must find two maps $\sigma :\mathcal{H}_{{\rm
%E}}({\rm T}^{2} )\rightarrow \mathcal{H}(G)$ and $\rho
%:\mathcal{H}(G)\rightarrow \mathcal{H} _{{\rm E}}({\rm T}^{2})$
%such that $\rho \sigma $ is the identity map on $\mathcal{H
%}_{{\rm E}}({\rm T}^{2})$, and $\sigma \rho $ is the identity map
%on $\mathcal{H}(G)$. \ Define $\sigma $ by
%\begin{equation}
%\sigma (f)=\sum_{{\large g\in G}}T^{\prime }(g)f.
%\end{equation}
Note that $\sigma $ is well defined since, for any $g_{0}\in F$,
\begin{eqnarray}
T^{\prime }(g_{0})[\sigma (f)] & = & T^{\prime }(g_{0}) \left
\{\frac{1}{|{\rm K}_{F}|}\sum_{{\large g\in F} }T^{\prime }(g)f
\right \}=\frac{1}{|{\rm K}_{F}|}\sum_{{\large g\in F}}T^{\prime
}(g_{0})T^{\prime }(g)f= \nonumber \\
\label{bhjooolkoolkiiio} && \frac{1}{|{\rm K}_{F}|}\sum_{{\large
g\in F}}T^{\prime }(g_{0}g)f.
\end{eqnarray}
But, as $g$ runs over $F$, so does $g_{0}g$, so the final sum is
just $ \sigma (f)$. \ So $\sigma (f)\in \mathcal{H}(F)$, as
required.
%\ Define $ \rho $ by
%\begin{equation}
%(\rho (f))(x)=\chi _{{\rm E}}(x)f(x)=(\chi _{{\rm E}}\cdot f)(x).
%\end{equation}
%Here $\chi _{{\rm E}}$ is the characteristic function of ${\rm
%E}$, equal to $1$ inside $ E$, and $0$ outside. \ The dot denotes
%pointwise multiplication.
We first compute $\rho \sigma $. \ We have
\begin{eqnarray}
(\sigma (f))(x) & = & \frac{1}{|{\rm K}_{F}|}\sum_{{\large g\in
F}}(T^{\prime }(g)f)(x)=\frac{1}{|{\rm K}_{F}|} \sum_{{\large g\in
F}}\gamma(x,g)f(xg) \nonumber \\
& =& \frac{1}{|{\rm K}_{F}|} \left \{
\gamma(x,e)f(xe)+\gamma(x,g_{2})f(xg_{2})
+\gamma(x,g_{3})f(xg_{3})+ \right . \nonumber \\
&&\left . + ... +\gamma(x,g_{n})f(xg_{n}) \right \},
\end{eqnarray}
where $e=g_{1},g_{2},g_{3},...,g_{n}$ are the elements of $F.$ \
Also $\rho (\sigma (f))=\chi _{E}.(\sigma (f))$, so, evaluating at
$x$,
\begin{eqnarray}
(\rho (\sigma (f)))(x) & = & \chi _{E}(x)(\sigma (f))(x)  = \nonumber \\
 &&\frac{1}{|{\rm K}_{F}|}\chi _{{\rm
E}}(x) \left [\gamma(x,e)f(xe)+\gamma(x,g_{2})f(xg_{2})+\gamma(x,g_{3})f(xg_{3})+ \right. \nonumber \\
\label{nmcmckdseraw} && +...+\gamma(x,g_{n})f(xg_{n})].
\end{eqnarray}
If $x\in {\rm E}$, $xg_{i}\notin {\rm E}$ when  $\; g_{i} \notin
{\rm K}_{F},\;$ $\;i=2,3,...,n,\;$ because ${\rm E}$ is
elementary. Therefore, the RHS of Eq. (\ref{nmcmckdseraw}) equals
to
\begin{eqnarray}
(\rho (\sigma (f)))(x) & = & \chi _{E}(x)(\sigma (f))(x)  = \nonumber \\
 &&\frac{1}{|{\rm K}_{F}|}\chi _{{\rm
E}}(x) \left [\gamma(x,e)f(xe)+\gamma(x,k_{2})f(xk_{2})+\gamma(x,k_{3})f(xk_{3})+ \right. \nonumber \\
\label{mkilollpkjiu} && +...+\gamma(x,k_{|{\rm
K}_{F}|})f(xk_{|{\rm K}_{F}|})],
\end{eqnarray}
where $\; {\rm K} \cap F= \left \{ e, k_{2}, k_{3}, ... , k_{|{\rm
K}_{F}|} \right \}. \;$ \ If $x\notin {\rm E}$, the RHS of Eq.
(\ref{nmcmckdseraw}) is $0$.  Since $f\in \mathcal{H}_{{\rm
E}}(
\mathcal P),
%{\rm T}^{2}
$ this means that for all such $f$, \be
\label{mkloioloplko}  (\rho (\sigma (f)))(x)= \frac {1}{|{\rm
K}_{F}|} \left [ \gamma(x,e) + \gamma(x,k_{2}) + \gamma(x,k_{3}) +
... + \gamma(x,k_{|{\rm K}_{F}|}) \right ]f(x), \ee where, $\; x
\in
\mathcal P,
%{\rm T}^{2}
\;$ $\; f \in \mathcal{H}_{{\rm E}}(
\mathcal P),
%{\rm T}^{2}
\;$ $\;|{\rm K}_{F}|\;$ denotes the order of the group $\; {\rm
K} \cap F \;$ and $\; {\rm K} \cap F= \left \{ e, k_{2}, k_{3},
... , k_{|{\rm K}_{F}|} \right \}, \;$ as required.
%$\rho \sigma (f)=f$,
%as required.

\vspace{0.3cm} \noindent Next, we compute $\sigma \rho $. \ We
have
\begin{equation}
\sigma (\rho (f))=\frac{1}{|{\rm K}_{F}|}\sum_{{\large g\in
F}}T^{\prime }(g)(\rho (f))=\frac{1}{|{\rm K}_{F}|}\sum_{ {\large
g\in F}}T^{\prime }(g)(\chi _{{\rm E}}\cdot f).
\end{equation}
Evaluating at $x$,
\begin{eqnarray}
%\be
(\sigma (\rho (f)))(x) & = & \frac{1}{|{\rm K}_{F}|}\sum_{{\large
g\in F}}(T^{\prime }(g) (\chi _{{\rm E}}\cdot f))(x)=
\frac{1}{|{\rm K}_{F}|}\sum_{{\large g\in F}}\gamma(x,g)(\chi
_{{\rm E}}\cdot f)(xg)=
\nonumber \\
&& \frac{1}{|{\rm K}_{F}|}\sum_{{\large g\in F} }\gamma(x,g)\chi
_{{\rm E}}(xg)f(xg).
%\ee
\end{eqnarray}
But $f$ is invariant, so $\gamma(x,g)f(xg)=f(x)$ for all $x\in
\mathcal P,
%{\rm T}^{2}
$ so the sum on the RHS becomes
\begin{equation}
\frac{1}{|{\rm K}_{F}|}f(x)\sum_{{\large g\in F}}\chi _{{\rm
E}}(xg)= \frac{1}{|{\rm K}_{F}|} f(x) [\chi _{{\rm E}}(x)+\chi
_{{\rm E}}(xg_{2})+...+\chi _{{\rm E}}(xg_{n})].
\end{equation}
By (\ref{mkkkliojhgt}), almost every $x \in
\mathcal P
%{\rm T}^{2}
$ belongs
to $\; |{\rm K}_{F}| \;$ open sets of the form
 ${\rm E}g_{i}^{-1}$
$(i=1,2,...,n \ {\rm and} \ n=|F|)$. \ Thus $xg_{i}\in {\rm E}$
for  $\; |{\rm K}_{F}| \;$ values of $i$, for almost all $x$. \
The sum in square brackets thus, almost always, contributes $\;
|{\rm K}_{F}| \;$, ($1$ for each $\chi _{{\rm E}}(xg_{i})$, $i$
being one of these $\; |{\rm K}_{F}| \;$ values,
 and $0$ for
all the other terms). \ So, for almost all $x$, the sum is $\;
|{\rm K}_{F}|, \;$ and so the RHS is $f(x)$. \ Hence \be
\label{dnnmdjksmki}(\sigma (\rho (f)))(x)=f(x) \ee for almost all
$x$, and so, as Hilbert space maps, $\sigma \rho (f)=f$ for all
$f\in \mathcal{H}(F)$. \ This completes the proof.
\end{prf}

%\vspace{0.5cm}

\section{Finite Potential Little Groups}
\label{s3}

\noindent Define the groups $\; A, \;$ $\; B, \;$ $\; C, \;$ and
$\; D, \;$ as follows $\; A = \{ ({\rm I},{\rm I}),(-{\rm I},-{\rm
I}) \}, \;$ $\;B = \{ ({\rm I},{\rm I}),(-{\rm I},{\rm I}) \}, \;$
$\; C = \{  ({\rm I},{\rm I}),({\rm I},-{\rm I}) \}, \;$ and,  $\;
D = \{ ({\rm I},{\rm I}),  (-{\rm I},{\rm I}),    ({\rm I},-{\rm
I}), (-{\rm I},-{\rm I}) \}. \;$ The subgroups of $\; C_{n} \times
C_{m} \;$  fall into one of the following disjoint classes:

\begin{enumerate}
\item{ \bf {Class O}. \normalfont \ \ \ This class contains those
subgroups of $\; C_{n} \times C_{m} \;$ which are of odd order.}

\item{ $\bf {Class \; E_{1}}. $ \normalfont \ \ This class
contains those subgroups of $\; C_{n} \times C_{m} \;$ which are
of even order and contain the group $\; A \;$ but they do not
contain the group $\; D. \;$ }

\item{ $\bf {Class \; E_{2}}. $ \normalfont \ \ This class
contains those subgroups of $\; C_{n} \times C_{m} \;$ which are
of even order and contain the group $\; B \;$ but they do not
contain the group $\; D. \;$ }

\item{ $\bf {Class \; E_{3}}. $ \normalfont \ \ This class
contains those subgroups of $\; C_{n} \times C_{m} \;$ which are
of even order and contain the group $\; C \;$ but they do not
contain the group $\; D. \;$ }

\item{ $\bf {Class \; E_{4}}. $ \normalfont \ \ This class
contains those subgroups of $\; C_{n} \times C_{m} \;$ which are
of even order and contain the group  $\; D. \;$ }

\end{enumerate}
From Eq. (\ref{tzatziki}) we have that \be \label{bykilokmnj}
(T'(-{\rm I},-{\rm I})\phi)(x,y)=\phi(x,y). \ee We also have \be
\label{ncccdretjkilo} (T'(-{\rm I},{\rm I})\phi)(x,y)=-\phi(x,y),
\ee and \be \label{ncdsawernm} (T'({\rm I},-{\rm
I})\phi)(x,y)=-\phi(x,y). \ee Eq. (\ref{bykilokmnj}) implies that
every potential little group contains the element $\;(-{\rm
I},-{\rm I}).\;$ Eq. (\ref{ncccdretjkilo}) implies that the
invariant subspace of the groups which belong to the class $\bf
{Class \; E_{2}} $ consists just of one element, namely, the zero
function. Eq. (\ref{ncdsawernm}) implies that the invariant
subspace of the groups which belong to the class $\bf {Class \;
E_{3}} $ consists just of one element, namely, the zero function.
Finally, both   Eq. (\ref{ncccdretjkilo}) and Eq.
(\ref{ncdsawernm}) imply that the invariant subspace of the groups
which belong to the class $\bf {Class \; E_{4}} $ consists just of
one element, namely, the zero function. Now we examine in turn the
invariant subspaces of the groups which belong to the $\bf {Class
\; E_{1}} $ and to the class $\bf {Class \; O}. $

\vspace{0.3cm}

\noindent $\bf {Class \; E_{1}.} $ \ \ \ \ \ \ If $\; F \;$ is a
group which belongs to this class then $\; {\rm K}_{F} \equiv A =
\{ ({\rm I},{\rm I}),(-{\rm I},-{\rm I}) \}, \;$ and $\; |{\rm
K}_{F}|=2. \;$ Moreover, from Eq. (\ref{vdsermklokmnj}) we have
that \be \label{dserkmnjuiol} \gamma(x,e)= \gamma(x,k_{2})=1, \ee
where $\; k_{2} = (-{\rm I},-{\rm I}). \;$ Substituting the result
of Eq. (\ref{dserkmnjuiol}) into Eq. (\ref{mkloioloplko}) and
taking into account that $\; |{\rm K}_{F}|=2 \;$ we obtain \be
\label{lmnnbbgdhhjshdyyeu}  (\rho (\sigma (f)))(x)= f(x), \ee
where, $\; x \in
\mathcal P,
%{\rm T}^{2}
\;$ and $\; f \in \mathcal{H}_{{\rm
E}}(
\mathcal P
%{\rm T}^{2}
). \;$ Moreover, Eq. (\ref{dnnmdjksmki}) reads
$$
(\sigma (\rho (f)))(x)=f(x)
$$
for all $ \; f\in \mathcal{H}(F). \;$ From equations
(\ref{lmnnbbgdhhjshdyyeu}) and (\ref{dnnmdjksmki}) we conclude
that there are two maps, namely,  $\sigma :\mathcal{H}_{{\rm
E}}(
\mathcal P
%{\rm T}^{2}
)\rightarrow \mathcal{H}(F)$ and $\rho
:\mathcal{H}(F)\rightarrow \mathcal{H} _{{\rm E}}(
\mathcal P
%{\rm T}^{2}
)$
such that $\rho \sigma $ is the identity map on $\mathcal{H
}_{{\rm E}}(
\mathcal P
%{\rm T}^{2}
),$ and $\sigma \rho $ is the identity map
on $\mathcal{H}(F)$. Therefore, we have the bijection
$$
%\mathcal{H}({\rm E})\leftrightarrow
\mathcal{H}_{{\rm E}} (
\mathcal P
%{\rm T}^{2}
)\leftrightarrow \mathcal{H}(F).
$$
Combining this bijection with (\ref{jnmnhbgsdnmjk}),we have the
following bijections:
\begin{equation}
\mathcal{H}({\rm E})\leftrightarrow \mathcal{H}_{{\rm E}} (
\mathcal P
%{\rm
%T}^{2}
)\leftrightarrow \mathcal{H}(F),
\end{equation}
which completely describes all invariant functions for $F$; they
are just given by ``arbitrary'' functions on the elementary region
${\rm E}$ for $F$.

%\ Define $\sigma $ by
%\begin{equation}
%\sigma (f)=\sum_{{\large g\in G}}T^{\prime }(g)f.
%\end{equation}
%\label{zxcvbvcxzkmj}
%\ Hence
%\be \label{dnnmdjksmki}(\sigma (\rho (f)))(x)=f(x) \ee for
%almost all $x$, and so, as Hilbert space maps, $\sigma \rho (f)=f$
%for all $f\in \mathcal{H}(G)$. \ This completes the proof.
%$$ \label{ulbkyuvrdfcrfxcr}
%(T'(g,h)\phi)(x,y)=k_{g}^{-3}(x)s_{g}(x)k_{h}^{-3}(y)s_{h}(y)\phi(xg,yh).
%$$
%(Eq. (\ref{tzatziki})),

\vspace{0.3cm}

\noindent $\bf {Class \; O.} $ \ \ \ \ \ \ If $\; F \;$ is a group
which belongs to this class then $\; {\rm K}_{F}  = \{ ({\rm
I},{\rm I}) \}, \;$ and $\; |{\rm K}_{F}|=1. \;$ Moreover, from
Eq. (\ref{vdsermklokmnj}) we have that \be \label{nbbbxvdfdgstsbg}
\gamma(x,e)=1, \ee where $\; e = ({\rm I},{\rm I}). \;$
Substituting the result of Eq. (\ref{nbbbxvdfdgstsbg}) into Eq.
(\ref{mkloioloplko}) and taking into account that $\; |{\rm
K}_{F}|=1 \;$ we obtain \be \label{nmjkmkjjjjjjsssss}  (\rho
(\sigma (f)))(x)= f(x), \ee where, $\; x \in
\mathcal P,
%{\rm T}^{2}
\;$ and
$\; f \in \mathcal{H}_{{\rm E}}(
\mathcal P
%{\rm T}^{2}
). \;$ Moreover, Eq.
(\ref{dnnmdjksmki}) reads
$$
(\sigma (\rho (f)))(x)=f(x)
$$
for all $ \; f\in \mathcal{H}(F). \;$ From equations
(\ref{nmjkmkjjjjjjsssss}) and (\ref{dnnmdjksmki}) we conclude that
there are two maps, namely,  $\sigma :\mathcal{H}_{{\rm E}}(
\mathcal P
%{\rm
%T}^{2}
)\rightarrow \mathcal{H}(F)$ and $\rho
:\mathcal{H}(F)\rightarrow \mathcal{H} _{{\rm E}}(
\mathcal P
%{\rm T}^{2}
)$
such that $\rho \sigma $ is the identity map on $\mathcal{H
}_{{\rm E}}(
\mathcal P
%{\rm T}^{2}
),$ and $\sigma \rho $ is the identity map
on $\mathcal{H}(F)$. Therefore, we have the bijection
$$
%\mathcal{H}({\rm E})\leftrightarrow
\mathcal{H}_{{\rm E}} (
\mathcal P
%{\rm T}^{2}
)\leftrightarrow \mathcal{H}(F).
$$
Combining this bijection with (\ref{jnmnhbgsdnmjk}), we have the
following bijections:
\begin{equation}
\mathcal{H}({\rm E})\leftrightarrow \mathcal{H}_{{\rm E}} (
\mathcal P
%{\rm
%T}^{2}
)\leftrightarrow \mathcal{H}(F),
\end{equation}
which completely describes all invariant functions for $F$; they
are just given by ``arbitrary'' functions on the elementary region
${\rm E}$ for $F$.

%\ Define $\sigma $ by
%\begin{equation}
%\sigma (f)=\sum_{{\large g\in G}}T^{\prime }(g)f.
%\end{equation}
%\label{zxcvbvcxzkmj}
%\ Hence
%\be \label{dnnmdjksmki}(\sigma (\rho (f)))(x)=f(x) \ee for
%almost all $x$, and so, as Hilbert space maps, $\sigma \rho (f)=f$
%for all $f\in \mathcal{H}(G)$. \ This completes the proof.
%$$ \label{ulbkyuvrdfcrfxcr}
%(T'(g,h)\phi)(x,y)=k_{g}^{-3}(x)s_{g}(x)k_{h}^{-3}(y)s_{h}(y)\phi(xg,yh).
%$$
%(Eq. (\ref{tzatziki})),

%Combining this bijection with (\ref{jnmnhbgsdnmjk}), and returning
%to our original notation (writing $F$ for $G$), we have the
%following bijections:
%\begin{equation}
%\mathcal{H}({\rm E})\leftrightarrow \mathcal{H}_{{\rm E}} ({\rm
%T}^{2})\leftrightarrow \mathcal{H}(F),
%\end{equation}
%which completely describes all invariant functions for $F$; they
%are just given by ``arbitrary'' functions on the elementary region
%${\rm E}$ for $F$.
By collecting the previous results we have the following theorem:

\begin{thrm}
\label{kilopjuhynbgderts} The invariant subspace of the groups
which belong to the classes $\; {\rm E}_{2}, \;$ $\; {\rm E}_{3},
\;$ and $\; {\rm E}_{4} \;$ consists just of one element, namely,
the zero function. The invariant subspace of a group $\; F \;$
which either belongs to the class $\; {\rm O} \;$ or to the class
$\; {\rm E}_{1} \;$ consists of functions which are ``arbitrary''
on the elementary region ${\rm E}$ for $F.$ In particular, the
following bijections hold
$$
\mathcal{H}({\rm E})\leftrightarrow \mathcal{H}_{{\rm E}} (
\mathcal P
%{\rm
%T}^{2}
)\leftrightarrow \mathcal{H}(F).
$$
\end{thrm}

A subgroup of $\; C_{n} \times C_{m} \;$ of odd order  can never
be a potential little group. Indeed, let $\; O \;$ be such a
group. Then, there always exists a bigger group, namely $\; O
\times A, \;$  where, $\;A = \{ ({\rm I},{\rm I}), (-{\rm I},-{\rm
I}) \} \simeq  Z_{2} \;$ which has the same invariant vectors as
$\; O. \;$ So, by taking into account theorem
\ref{kilopjuhynbgderts} we arrive at the following

\begin{prop}
The only finite potential little groups are those subgroups of $\;
C_{n} \times C_{m} \;$ which fall into the class $\; {\rm E}_{1}.
\;$
\end{prop}

%\newpage

\section{All the Finite Potential Little Groups are Actual}
\label{s4}
%A subgroup of $\; C_{n} \times C_{m} \;$ of odd order  can never be a potential little group. Indeed, let $\; O \;$
%be such a group. Then, there always exists a bigger group , namely $\; O \times A \;$ , where,
%$\;A = \{ ({\rm I},{\rm I}), (-{\rm I},-{\rm I}) \} \simeq  Z_{2} \;$ which has the same invariant vectors
%as $\; O \;$. So, by taking into account theorem \ref{kilopjuhynbgderts} we arrive at the following

%\begin{prop}
%The only finite potential little groups are those subgroups of $\; C_{n} \times C_{m} \;$
%which fall into the class $\; {\rm E}_{1} \;$.
%\end{prop}

\ So far, for all potential little groups $\; F, \;$ we  have
given an explicit description of the invariant subspaces
$\mathcal{H}(F)\subset \mathcal{H} (
\mathcal P
%{\rm T}^{2}
)$; $\;
\mathcal{H}(F) \;$ is in bijective correspondence with $\;
\mathcal{H}({\rm E}). \;$ \ For most  of the vectors $\;f \in
\mathcal{H}(F),\;$ the little group will be precisely $\; F. \;$ \
However, to be certain that any given potential little group
really does occur as a little group, we must exclude the
possibility that all vectors $\;f \in \mathcal{H}(F)\;$ are also
invariant under a bigger group $\; \Theta \supset F. \;$ So, in
order to prove that $\; F \;$ does occur as a little group it
suffices to prove that there exists  one $\;f' \in
\mathcal{H}(F)\;$ which has no higher symmetry;  there is no group
$\; \Theta \supset F \;$ which leaves $\; f' \;$ invariant. The
proof will be constructive. For each $\; F \;$ we will give an $\;
f_{F} \in \mathcal{H}(F)\;$ which is invariant only under $\; F.
\;$
%give , for each $\; F \;$, a specific example of
%$f \in \mathcal{H}(F)$
%with no ``higher symmetry''; no element $k\notin F$ with $
%T^{\prime }(k)f =f $.

The subgroups $\; F \;$ which fall into the class $\; {\rm E}_{1}
\; $ may be taken as subgroups $F\subset C_{n}\times C_{m}$, where
$\pi _{1}(F)=C_{n}$ and $\pi _{2}(F)=C_{m}$. Since the groups $\;
F \;$ contain the element  $\; (-{\rm I},-{\rm I}) \;$ both the
numbers $\;n \;$ and $\; m \;$ are even. \ Let $F_{nm}$ be the
open rectangle
\begin{equation}
{\rm F}_{nm}={\rm E}_{n}\times {\rm E}_{m}\subset \mathcal P
%{\rm T}^{2}
\end{equation}
%defined in
%Section \ref{alskjdfhgvbcn},
where ${\rm E}_{n}=\{\theta \in S^{1}$ \ $\left| 0<\theta <4\pi
/n\right. \}$. \ Let $R$ be the open rectangle with center $(2\pi
/n,2\pi /m)$ and side lengths $\pi /n$ and $\pi /m$. \ That is,
$R$ is the set
\begin{equation}
R=(\frac{2\pi }{n}-\frac{\pi }{2n},\frac{2\pi }{n}+\frac{\pi
}{2n})\times ( \frac{2\pi }{m}-\frac{\pi }{2m},\frac{2\pi
}{m}+\frac{\pi }{2m}).
\end{equation}

Now consider the non$-$effective action $
%{\rm T}^{2}
\mathcal P
\times (C_{n}
\times C_{m}) \rightarrow
\mathcal P
%{\rm T}^{2}
$ defined by equation
(\ref{bhnjmkloiulkjuhygtfr}).\ Let $R_{ij}$ the translate of $R$
by the element $(g_{i},g_{j})$ of $C_{n}\times C_{m}$:
\begin{equation}
R_{ij}=R(g_{i},g_{j}).
\end{equation}
%These sets are all disjoint.
More generally, let $R_{\omega \xi }$ denote the translate of $R$
by any element $(g(\omega ),g(\xi ))\in K$; $\; K=SO(2) \times
SO(2). \;$
\begin{equation}
R_{\omega \xi }=R(g(\omega ),g(\xi )).
\end{equation}
Then we have the following Propositions (\ref{drftgvcxseawqax}), (\ref{jkimnxxxsdhgjuk}), and  (\ref{nmvcvvcvvcfdr}), regarding $R_{\omega \xi }$. 
The proofs of these Propositions follow closely the proofs of the Propositions 4, 5, and 6 
given in \cite{macmel}, and as such are omitted.
\begin{prop}
\label{drftgvcxseawqax} Let $\;R_{\omega \xi }\;$ be any translate
of $\;R.\;$ Then $\;R_{\omega \xi }\;$ cannot intersect two
distinct rectangles of the form $R_{ij}$.
\end{prop}

%\begin{prf}
%The proof is similar to the proof of  Proposition \ref
%{drftgvcxseawqax}.
%\end{prf}
%\begin{prop}
%\label{drftgvcxseawqax} Let $\;R_{\omega \chi }\;$ be any
%translate of $\;R\;$. Then $\;R_{\omega \chi }\;$ cannot intersect
%two distinct rectangles of the form $R_{ij}$.
%\end{prop}
%\begin{prf}
%\normalfont
%\ Let $p\in R_{ij}$ and $q\in R_{kl}$ be two points in the
%distinct rectangles $R_{ij}$ and $R_{kl}$. \ Then the minimum difference
%between the $\theta $ coordinates of $p$ and $q$ is given by this minimum
%quantity for the closest pair of rectangles, for example $R=R_{00}$ and $
%R_{01}$. \ This minimum is
%\begin{equation}
%\frac{2\pi }{n}-\frac{\pi }{2n}=\frac{3\pi }{2n}.
%\end{equation}
%Now let $r$ and $s$ be two points of $R_{\omega \chi }$. \ Then the minimum
%diference between the $\theta $ coordinates of $r$ and $s$ is the side
%length $\pi /2n$ of $R_{\omega \chi }$. \ Similarly for the $\phi $
%coordinates, with $n$ replaced by $m$. \ Since $3\pi /2n>\pi /2n$ and $3\pi
%/2m>\pi /2m$, $R_{\omega \chi }$ cannot intersect both rectangles $R_{ij}$
%and $R_{kl}$. \ This completes the proof.
%\end{prf}

Now let $F$ be any finite subgroup of $K$ which falls into the
class of subgroups $\; {\rm E}_{1},$  with $\pi _{1}(F)=C_{n}$
and $\pi _{2}(F)=C_{m}.$ \ So $F\subset C_{n}\times C_{m}.$ \
Write, for convenience, $G=F,$ and denote the elements of $G$ as
follows:
\begin{equation}
G=\{g_{1},g_{2},g_{3},...,g_{l}\}
\end{equation}
where $g_{1}=e$ and $l$ is the order of $G$. \ Let $A$ denote the
union
% \tag*{(10.5)}
\begin{equation}
\label{njbhgyjukolopolkimnj}
%\label{wsljytbvcdfhga}
A=Rg_{1}\cup Rg_{2}\cup ...\cup Rg_{l}\equiv RG,
\end{equation}
and let $B$ denote the  union
%\tag*{(10.6)}
\begin{equation}
\label{mnjbhvgcfxd}
%\label{plkjhnbvcfdresxwa}
B=R_{\omega \xi }g_{1}\cup R_{\omega \xi }g_{2}\cup ...\cup
R_{\omega \xi }g_{l}=R_{\omega \xi }G.
\end{equation}
Then we have the following:

%\textbf{10.2}
\begin{prop}
\label{jkimnxxxsdhgjuk} Let $\;A\;$ and $\;B\;$ be the sets
defined by (\ref{njbhgyjukolopolkimnj}) and (\ref{mnjbhvgcfxd}). \
Then either no set of the form $\;R_{\omega \xi }g_{i}\;$
intersects any set of the form $\;Rg_{j},\;$  or every set of the
form $\;R_{\omega \xi }g_{i}\;$ intersects exactly two identical
sets of the form $\;Rg_{j}.\;$  \ In the latter case, we have, for
some $\;k,\;$ the
%disjoint
union
%\tag*{(10.7)}
\begin{equation}
\label{ffdsswqqyijj} A\cap B={\large \cup }_{j=1}^{l}(R_{\omega
\xi }\cap Rg_{k})g_{j}.
\end{equation}
\end{prop}
%\begin{prf}
%The proof is similar to the proof of Proposition
%\ref{cddxxsszzeerrffgyt}
%\end{prf}
%\begin{prf}
%\normalfont
%\ Suppose first that $R_{\omega \chi }$ intersects none of
%the sets $Rg_{j}$; $R_{\omega \chi }\cap Rg_{j}=\varphi $ for all $j$. \
%Then ($R_{\omega \chi }\cap Rg_{j})g_{i}=R_{\omega \chi }g_{i}\cap
%Rg_{j}g_{i}=\varphi $ for all $i$ and $j$; this is the first possibility. \
%Suppose next that $R_{\omega \chi }$intersects at least one of the sets $
%Rg_{j}$. \ Then, by Proposition \ref{drftgvcxseawqax},
%$R_{\omega \chi }$intersects exactly
%one of the sets $Rg_{j}$, say $Rg_{k}$. \ It follows that each set of the
%form $R_{\omega \chi }g_{j}$ intersects exactly one of the $Rg_{i}$, namely $
%Rg_{k}g_{j}$. \ So, taking the intersection of $A$ and $B$
%in (\ref{wsljytbvcdfhga}) and
%(\ref{plkjhnbvcfdresxwa}) gives the formula displayed
%in the Proposition. \ This completes the
%proof.
%\end{prf}

Now, for any measurable subset $S$ of $
\mathcal P
%{\rm T}^{2}
$, let $\alpha
(S)$ denote the area of $S$. \ Note that, for any $(\omega ,\xi
)$, we have $\alpha (R_{\omega \xi })=\alpha (R)$; the area is
invariant under all translations. \ Now we have the following
result.
%extbf{10.3}:
\begin{prop}
\label{nmvcvvcvvcfdr} Let $\;A=RG\;$ and $\; B=R_{\omega \xi }G\;$
be the sets given by (\ref{njbhgyjukolopolkimnj}) and
(\ref{mnjbhvgcfxd}), and suppose that
% \tag*{(10.8)}
\begin{equation}
\label{pollkkmmmmncd} \alpha (A)=\alpha (A\cap B).
\end{equation}
Then $\;R_{\omega \xi }=Rg_{k},\;$  the sets $\;A\;$ and $\;B\;$
coincide, and either $\;(g(\omega ),g(\xi ))=g_{k}\;$ or
$\;(g(\omega ),g(\xi ))=g_{k}({\rm I},{\rm I})\;$.
\end{prop}
%\begin{prf}
%The proof is similar to the proof of Proposition
%\ref{ytfsdtresdnn}.
%\end{prf}

%\bigskip
We are now ready to construct a function $\zeta _{3}\in
\mathcal{H}(G)$ with little group $G.$ \ Let $\chi _{R}$ be the
characteristic function of $R,$ and define $\zeta _{3}$ by
\begin{equation}
\zeta _{3}=\frac{1}{|{\rm K}_{F}|}\sum_{{\large g\in G}}T^{\prime
}(g^{-1})\chi _{R}.
\end{equation}
By construction, (as in the proof of Theorem
\ref{zaaasdhjjnmkijk}, (Eq. (\ref{bhjooolkoolkiiio})), $\zeta
_{3}\in \mathcal{H}(G).$ Recall that we denote by $\; \gamma(x,g)
\;$ the multiplier $\; k_{g}^{-3}(x)s_{g}(x)k_{h}^{-3}(y)s_{h}(y)
\;$ in the representation $ \; (T'(g)f)(x)= \gamma(x,g)
%k_{g}^{-3}(x)s_{g}(x)k_{h}^{-3}(y)s_{h}(y)
f(xg), \;$ where $\; x \;$ denotes a point  of the torus and $\; g
\;$ an element of $\; {\rm S}{\rm L}(2,R) \times {\rm S}{\rm
L}(2,R). \;$ Henceforth, for notational convenience we will write
$\; \gamma_{g}(x) \;$ instead of $\; \gamma(g,x). \;$

%\begin{eqnarray}
%\rho & : & \mathcal{H}(F)\rightarrow \mathcal{H} _{{\rm E}}({\rm
%T}^{2}) \nonumber \\
%\label{bbxnnxxnhhskkaju}
%(\rho (f))(x) & = & \chi _{{\rm E}}(x)f(x)=(\chi _{{\rm E}}\cdot
%f)(x),
%\end{eqnarray}
% where $\; x \in {\rm T}^{2}, \;$ $\; f \in \mathcal{H}(F) \;$,
% and $\chi _{{\rm E}}$
%denotes the characteristic function of ${\rm E}$, equal to $1$
%inside $\;{\rm  E} \;$, and $\;0\;$ outside. \ The dot denotes
%pointwise multiplication. We now prove the following theorem

%Subsequently we will write $$ \label{mkiuyhhhnbgdfdffdr}
%(T'(g)f)(x)= \gamma(x,g)f(xg), $$ and it will be understood that
%in our situation
%the multiplier $\;\gamma(x,g)\;$ equals to
%\be
%\label{vdsermklokmnj}
%\gamma(x,g)=k_{g}^{-3}(x)s_{g}(x)k_{h}^{-3}(y)s_{h}(y).
%\ee

Since for $\; g \in G \;$ we have $\; T^{\prime }(g^{-1})\chi
_{R}=\gamma_{g^{-1}} \cdot \chi _{Rg},\; $ we obtain
%, and since the
%sets $Rg$ are disjoint,
\begin{equation}
\zeta _{3}=\frac{1}{|{\rm K}_{F}|}\sum_{{\large g\in
G}}\gamma_{g^{-1}} \cdot \chi _{Rg}.  \qquad
\end{equation}
%In Eq. (\ref{bxxmjknuiolkolkijkm})
The $\; \cdot \;$ in $\; \gamma_{g^{-1}} \cdot \chi _{Rg} \;$
denotes pointwise multiplication. Now note that if
$\;h^{-1}=(g(\omega )^{-1},g(\xi )^{-1}) \;$ is an element of $\;
SO(2) \times SO(2), \;$ then
%\newline
$ \; T^{\prime }(g(\omega )^{-1},g(\xi )^{-1})\zeta _{3} \;$
evaluated at $\; x= ( \rho, \sigma ) \;$ gives
\begin{eqnarray}
\label{bxxmjknuiolkolkijkm} \left ( T^{\prime }(g(\omega
)^{-1},g(\xi )^{-1})\zeta _{3} \right ) (\rho , \sigma) &  =  &
\frac{1}{|{\rm K}_{F}|} \sum_{{\large g\in G}}\left ( T^{\prime
}(h^{-1})\left ( \gamma_{g^{-1}} \cdot \chi _{Rg} \right )
\right )( \rho, \sigma )  \nonumber \\
%\label{bhyioklopoiuhjy}
& = & \frac{1}{|{\rm K}_{F}|} \sum_{{\large g\in G}} \left \{
\gamma_{h^{-1}}(\rho , \sigma )
\gamma_{g^{-1}}(\rho-2\omega,\sigma-2\xi) \right . \nonumber \\
\label{bhyioklopoiuhjy} && \left . \gamma_{Rgh}(\rho,\sigma)
\right \}
%\chi _{A(g(\omega
%),g(\chi ))}=\chi _{R_{\omega \chi }G}=\chi _{B}.
\end{eqnarray}
%In Eq. (\ref{bxxmjknuiolkolkijkm}) the $\; \cdot \;$ in
%$\; \gamma_{g^{-1}} \cdot \chi _{Rg} \;$
%denotes pointwise multiplication.
%In Eq. (\ref{bhyioklopoiuhjy})
%we used the fact that $\; SO(2) \times SO(2) \;$ is abelian and therefore $\; Rgh=Rhg \;$( $\; h \in
%SO(2) \times SO(2) \quad {\rm and} \quad g \in G \subset SO(2) \times SO(2) \;$).
A simple calculation shows
$$
\gamma_{g^{-1}}(\rho-2\omega,\sigma-2\xi)=
\gamma_{g^{-1}h^{-1}}(\rho, \sigma) \gamma_{h^{-1}}(\rho, \sigma).
\quad
$$
Substituting back into Eq. (\ref{bhyioklopoiuhjy}) gives \be
T^{\prime }(h^{-1})\zeta _{3}  = \frac{1}{|{\rm
K}_{F}|}\sum_{{\large g\in G}} \gamma_{g^{-1}h^{-1}} \cdot \chi
_{Rgh}.  \qquad  \ee
%where
%the $\; \cdot \;$ in
%$\; \gamma_{g^{-1}h^{-1}} \cdot \chi _{Rgh} \;$
%denotes pointwise multiplication.

\noindent
We are now ready to prove the following:

%{Theorem 10.4}:
\begin{thrm}
The little group of $\;\zeta _{3}\;$ is precisely $\;G{\rm ;}$
$\;L(\zeta _{3})=G.\;$
\end{thrm}
\begin{prf}
\normalfont \ We must prove that, if $T^{\prime }(g(\omega
)^{-1},g(\xi )^{-1})\zeta _{3}=\zeta _{3}$, then
\begin{equation}
(g(\omega )^{-1},g(\xi )^{-1})\in G.
\end{equation}
The equation holds in the Hilbert space sense, and is just
$$
\sum_{{\large g\in G}} \gamma_{g^{-1}h^{-1}} \cdot \chi _{Rgh} =
\sum_{{\large g\in G}}\gamma_{g^{-1}} \cdot \chi _{Rg}.   \qquad
$$
%$\chi _{A}=\chi
%_{B}$.
So we have
% \tag*{(10.9)}
\begin{eqnarray}
\label{jkcygfgfqbkxjql}
0 & = & \left\|  \zeta _{3} - T^{\prime }(h^{-1})\zeta _{3} \right\| ^{2}  =  \nonumber \\
\label{mxxxsdfhjyu}
%\left
&&\int_{{\rm T}^{2}} \left [ \sum_{{\large g\in G}}\left (
\gamma_{g^{-1}} \cdot \chi _{Rg} \right ) (x)  - \sum_{{\large
g\in G}} \left (  \gamma_{g^{-1}h^{-1}} \cdot \chi _{Rgh} \right )
(x)
%\chi_{A}(x)-\chi _{B}(x)
\right ]^{2}d\mu (x).
%\right .
\end{eqnarray}
By using the expansion formula for
$\;(a_{1}+a_{2}+a_{3}+a_{4}+...+a_{\nu})^{2},\;$ where $\;
a_{1},a_{2},...,a_{\nu} \;$ are any real numbers, and the
equalities $ \;
%\label{kolkijjjhnju}
\chi_{{\rm E}}^{2}(x)=\chi_{{\rm E}}(x) \quad {\rm and} \quad
\gamma_{g}^{2}=1, \quad $ where $\; g \in G \subset C_{n} \times
C_{m}, \;$  Eq. (\ref{mxxxsdfhjyu}) gives

%Recall that for any real numbers $\; a_{1},a_{2},...,a_{\nu} \;$ we have
%\begin{eqnarray}
%(a_{1}+a_{2}+a_{3}+a_{4}+...+a_{\nu})^{2} & = & a_{1}^{2}+a_{2}^{2}+a_{3}^{2}+a_{4}^{2}+...
%+a_{\nu}^{2}+ \nonumber \\
%&&2a_{1}(a_{2}+a_{3}+a_{4}+...+a_{\nu})+2a_{2}(a_{3}+a_{4}+...+a_{\nu})+ \nonumber \\
%\label{njiklomnju}
%&&2a_{3}(a_{4}+...+a_{\nu})+...+2a_{\nu-1}a_{\nu} \quad .
%\end{eqnarray}
%We also have
%\be
%\label{kolkijjjhnju}
%\chi_{{\rm E}}^{2}(x)=\chi_{{\rm E}}(x) \qquad {\rm and} \qquad \gamma_{g}^{2}=1 \quad,
%\ee
%where $\; g \in G \subset C_{n} \times C_{m} \;$.
%By using Eqs. (\ref{njiklomnju}) and (\ref{kolkijjjhnju}) Eq. (\ref{mxxxsdfhjyu}) gives

\be \label{mikkllloooikjuhygt} 0  = \int_{{\rm T}^{2}} \left \{
\sum_{{\large g\in G}} \chi _{Rg}  (x) +\sum_{{\large g\in G}}
 \chi _{Rgh}  (x) + \mathcal{B} + \mathcal{C} - \Lambda
\right \}{\rm d}\mu (x), \qquad  \ee where, denoting by $\; l \;$
the order of the group $\; G, \;$  we have
\begin{eqnarray}
\mathcal B & = & 2(\gamma_{g_{1}^{-1}} \cdot \chi_{Rg_{1}})(x)
 \sum_{i \geq 2}^{l}(\gamma_{g_{i}^{-1}} \cdot \chi_{Rg_{i}})(x)+  \\
&&2(\gamma_{g_{2}^{-1}} \cdot \chi_{Rg_{2}})(x) \sum_{i \geq
3}^{l}(\gamma_{g_{i}^{-1}} \cdot \chi_{Rg_{i}})(x)+...+
2(\gamma_{g_{l-1}^{-1}} \cdot
\chi_{Rg_{l-1}})(x)(\gamma_{g_{l}^{-1}} \cdot \chi_{Rg_{l}})(x),
\quad  \nonumber
\end{eqnarray}

\begin{eqnarray}
\mathcal C & = & 2(\gamma_{g_{1}^{-1}h^{-1}} \cdot
\chi_{Rg_{1}h})(x)
 \sum_{i \geq 2}^{l}(\gamma_{g_{i}^{-1}h^{-1}} \cdot \chi_{Rg_{i}h})(x)+  \\
&&2(\gamma_{g_{2}^{-1}h^{-1}} \cdot \chi_{Rg_{2}h})(x)
\sum_{i \geq 3}^{l}(\gamma_{g_{i}^{-1}h^{-1}} \cdot \chi_{Rg_{i}h})(x)+...+ \nonumber \\
&&2(\gamma_{g_{l-1}^{-1}h^{-1}} \cdot
\chi_{Rg_{l-1}h})(x)(\gamma_{g_{l}^{-1}h^{-1}} \cdot
\chi_{Rg_{l}h})(x), \quad  \nonumber
\end{eqnarray}

\be \Lambda  =  2 \sum_{g_{i} \in G}     \left
\{(\gamma_{g_{i}^{-1}} \cdot \chi_{Rg_{i}})(x) \left \{
\sum_{g_{j} \in G}(\gamma_{g_{j}^{-1}h^{-1}} \cdot
\chi_{Rg_{j}h})(x)\right \} \right \}.
%(\gamma_{g_{1}^{-1}h^{-1}} \cdot \chi_{Rg_{1}h})(x)
% \sum_{i \geq 2}^{l}(\gamma_{g_{i}^{-1}h^{-1}} \cdot \chi_{Rg_{i}h})(x)+  \\
%&&2(\gamma_{g_{2}^{-1}h^{-1}} \cdot \chi_{Rg_{2}h})(x)
%\sum_{i \geq 3}^{l}(\gamma_{g_{i}^{-1}h^{-1}} \cdot \chi_{Rg_{i}h})(x)+...+ \nonumber \\
%&&2(\gamma_{g_{l-1}^{-1}h^{-1}} \cdot \chi_{Rg_{l-1}})(x)(\gamma_{g_{l}^{-1}h^{-1}}
%\cdot \chi_{Rg_{l}})(x) \quad . \nonumber
\ee For convenience, call $\; \imath=(-{\rm I},-{\rm I}), \;$ and
enumerate the elements of $\; G \;$ as follows
$$
g_{1}=e, \  \ g_{2}, \  \ g_{3}, \, \ ... \ , \
g_{\frac{l}{2}+1}=\imath, \  \ g_{\frac{l}{2}+2}=g_{2}\imath, \  \
g_{\frac{l}{2}+3}=g_{3}\imath,  \  \ ... \ , \
g_{l}=g_{\frac{l}{2}}\imath.
$$
One then can show that \be \label{nhjjjfkfujhgydert} \mathcal B  =
 2\sum_{i =1}^{l/2}\chi_{Rg_{i}}(x), \quad  \quad
\mathcal C  =
 2\sum_{i =1}^{l/2}\chi_{Rg_{i}h}(x). \quad  \quad
\ee Combining Eq. (\ref{nhjjjfkfujhgydert}) with
$$
\sum_{{\large g\in G}} \chi _{Rg}  (x)= 2\sum_{i
=1}^{l/2}\chi_{Rg_{i}}(x), \quad  \quad \sum_{{\large g\in G}}
\chi _{Rgh}  (x) =  2\sum_{i =1}^{l/2}\chi_{Rg_{i}h}(x), \quad
$$
and substituting back into Eq. (\ref{mikkllloooikjuhygt}) we
obtain \be \label{kilopolokilokijuiujki} 0=\int_{{\rm T}^{2}}
\left (4\sum_{i =1}^{l/2}\chi_{Rg_{i}}(x)+4\sum_{i
=1}^{l/2}\chi_{Rg_{i}h}(x) - \Lambda \right ) {\rm d} \mu (x).
\quad  \ee Noting that $\;\sum_{i
=1}^{l/2}\chi_{Rg_{i}}(x)=\chi_{A}(x) \;$ and that $\;\sum_{i
=1}^{l/2}\chi_{Rg_{i}h}(x)= \chi_{B}(x), \;$  where $\; A \;$ and
$\; B \;$ are the areas which are displayed in equations
(\ref{njbhgyjukolopolkimnj}) and (\ref{mnjbhvgcfxd}) respectively,
Eq. (\ref{kilopolokilokijuiujki}) gives \be
\label{kiloplokiollkoikli} 0=4\int_{{\rm T}^{2}} \left
(\chi_{A}(x)+\chi_{B}(x) \right ){\rm d} \mu (x)   - \int_{{\rm
T}^{2}}\Lambda {\rm d} \mu (x). \quad  \ee
%Note that the characteristic functions cancel on the overlap $(A\cap B)$. \
%So $\chi _{A}(x)-\chi _{B}(x)=\chi _{A^{\prime }}(x)-\chi _{B^{\prime }}(x)$
Define  $A^{\prime }$ and $B^{\prime }$ to be  the disjoint sets
$A^{\prime }=A-(A\cap B)$ and $B^{\prime }=B-(A\cap B)$. By noting
that $\; \chi_{A}=\chi_{A^{\prime }} + \chi_{A\cap B} \;$ and that
$\; \chi_{B}=\chi_{B^{\prime }} + \chi_{A\cap B} \;$ Eq.
(\ref{kiloplokiollkoikli}) gives \be \label{njuihjhgtyfrddef}
0=4\alpha(A^{\prime }) + 4\alpha(B^{\prime }) + 8 \alpha (A\cap B)
- \int_{{\rm T}^{2}}\Lambda {\rm d} \mu (x), \quad  \ee where $\;
\alpha (A) \;$ denotes the area of the region $\; A. \;$ We
distinguish now the following 2 cases
%\vspace{0.5cm}

\noindent 1. \hspace{0.3cm}  $\;\alpha(A \cap B)=0 .\;$ In this
case we have $\;\int_{{\rm T}^{2}}\Lambda {\rm d} \mu (x)=0 \;$
and Eq. (\ref{njuihjhgtyfrddef}) gives \be 0=\alpha (A^{\prime })
+ \alpha(B^{\prime }). \quad  \ee
%\ So equation (\ref{jkcygfgfqbkxjql}) becomes
%\begin{eqnarray}
%0 & = & \left\| \chi _{A}-\chi _{B}\right\| ^{2}=\int_{{\rm T}^{2}}[\chi
%_{A^{\prime }}(x)-\chi _{B^{\prime }}(x)]^{2}d\mu (x) \\
%%\end{equnarray}
%\begin{equation}
%& = & \int_{{\rm T}^{2}}[\chi _{A^{\prime }}^{2}(x)+\chi _{B^{\prime
%}}^{2}(x)]d\mu (x)=\alpha (A^{\prime })+\alpha (B^{\prime }).
%\end{eqnarray}
Since areas are always non negative, it follows that $\alpha
(A^{\prime })=\alpha (B^{\prime })=0$. \ Since $A=A^{\prime }\cup
(A\cap B)$ (disjoint union), we have $\alpha (A)=\alpha (A^{\prime
})+\alpha (A\cap B)$ and therefore we obtain $\; \alpha(A)=0. \;$
But $\;\alpha(A)=\frac{l}{2}\alpha(R)>0. \;$  Contradiction. So
we cannot have $\; \alpha(A \cap B)=0.\;$
%\vspace{0.5cm}
\noindent
When $\alpha(A \cap B)> 0 \;$ from Proposition \ref
{jkimnxxxsdhgjuk} we have that $ \alpha(R_{h} \cap Rg_{k}) > 0,$
for some $k \in \{1,2,3,...,\frac{l}{2} \}$ ($R_{h}
\equiv R_{\omega \xi}$).
When $\; \alpha(R_{h} \cap Rg_{k}) >
0  \;$ either ($R_{h} \cap Rg_{k} \neq \varnothing \quad {\rm and}
\quad  R_{h} \neq Rg_{k}$) or $\;R_{h} = Rg_{k}\;$.
%($\;\O \;$
%denotes the empty set ) .
Now we show that the first alternative
cannot happen.

%\vspace{0.5cm}

\noindent 2. \hspace{0.3cm}  $\; \alpha(A \cap B)>0 \;$  and
($\;R_{h} \cap Rg_{k} \neq \varnothing,  \quad  \quad  R_{h} \neq
Rg_{k}$). In this case we have $\;\int_{{\rm T}^{2}}\Lambda {\rm
d} \mu (x)<8\alpha(A \cap B) \;$ and Eq. (\ref{njuihjhgtyfrddef})
gives \be \alpha(A^{\prime }) + \alpha(B^{\prime })<0. \quad  \ee
Contradiction, since areas are always non negative. So we cannot
have $\; \alpha(A \cap B)>0 \;$  and   ($\;R_{h} \cap Rg_{k} \neq
\varnothing,  \quad  \quad  R_{h} \neq Rg_{k}$).
%\vspace{0.5cm}
\noindent The only remaining possibility is $\; \alpha(A \cap B)>0
\;$  and   $\;R_{h} = Rg_{k}.\;$ Then $\; A=B \;$ and either
$h=\;(g(\omega ),g(\xi ))=g_{k}\;$ or $\;h=(g(\omega ),g(\xi
))=g_{k}(-{\rm I},-{\rm I}).\;$ This completes the proof.

%\vspace{0.5cm}

\noindent It is worth pointing out that that Equations $\;R_{h} =
Rg_{k}\;$ and $ 0=4\alpha(A^{\prime }) + 4\alpha(B^{\prime }) + 8
\alpha (A\cap B) - \int_{{\rm T}^{2}}\Lambda {\rm d} \mu (x) $ are
mutually consistent. Indeed, when $\;R_{h} = Rg_{k},\;$
$\;\int_{{\rm T}^{2}}\Lambda {\rm d} \mu (x)=8\alpha(A \cap B).
$ Substituting back into Eq. (\ref{njuihjhgtyfrddef}) we obtain
$\alpha (A^{\prime })=\alpha (B^{\prime })=0.$ The last Equation
combined with $\alpha (A)=\alpha (A^{\prime })+\alpha (A\cap B)$
gives $\;\alpha (A)=\alpha (A\cap B).  \;$ Then according
Proposition \ref{nmvcvvcvvcfdr} $\;R_{h }=Rg_{k}.$
%either $h=\;(g(\omega ),g(\xi ))=g_{k}\;$ or
%$\;h=(g(\omega ),g(\xi ))=g_{k}({\rm I},{\rm I})\;$.
% \vspace{0.5cm}
%\noindent
%3. \hspace{0.3cm}  $\; \alpha(A \cap B)>0 \;$  and   $\;R_{h} = Rg_{k}\;$.
%In this case we have $\;\int_{{\rm T}^{2}}\Lambda {\rm d} \mu (x)=8\alpha(A \cap B) \;$
%and Eq. (\ref{njuihjhgtyfrddef}) gives
%\be
%\alpha(A^{\prime }) + \alpha(B^{\prime })=0 \quad .
%\ee
%It follows that $\alpha (A^{\prime
%})=\alpha (B^{\prime })=0$. \ From  $\alpha (A)=\alpha (A^{\prime })+\alpha (A\cap B)$ we obtain
%$\;\alpha (A)=\alpha (A\cap B)  \;$. Then according Proposition \ref{nmvcvvcvvcfdr}
%$\;R_{h }=Rg_{k}\;$ and
%either $h=\;(g(\omega ),g(\xi ))=g_{k}\;$ or
%$\;h=(g(\omega ),g(\xi ))=g_{k}({\rm I},{\rm I})\;$.
%. \ So $
%\alpha (A)=\alpha (A\cap B)$.
%\ So, by Proposition \ref{ytfsdtresdnn}, $(g(\omega ),g(\chi
%))=g_{k}$ for some $g_{k}\in G$. \ This completes the proof.
\end{prf}
%We may collect together the results on little groups as follows:
%\textbf{10.5}:
\noindent
The previous results on little groups can be summarized
as  follows
%\begin{thrm}
%\label{xsdczzammnnbblk}
%The potential little groups of $B(2,2)$
%are as follows. \ There is exactly one two dimensional such
%group, namely $K$ itself. \ The one dimensional such groups are as
%follows. \ If $p$ and $q$ are both odd, $H(N,p,q)$
% is always such a group. \ If $p$ and $q$ have
%opposite parity, $H(N,p,q)$ is a potential little group if and
%only if, in the expression $p/N=p^{\prime }/N^{\prime }$ in lowest
%terms, the integer $N^{\prime }$ is even. \ The groups
%$C_{N}\times S^{1}$,
%$S^{1}\times C_{N}$ are potential
%little groups if and only if $N$
%is even. \ The zero dimensional
%potential little groups have been (partially) described in the paragraph
%preceding this theorem.
%\end{thrm}
\begin{thrm}
\label{nhmjknmnhbgvfcd} The actual little groups for $\; \H \B
%_c
\; $ are either one dimensional or zero dimensional. \ The infinite one
dimensional such groups are the groups $\; H(N,{\rm p},{\rm q}), \;
$  where, $ N, $ $  {\rm p}, $ and $  {\rm q}, $
are all odd, and, $  {\rm p}  $ and $  {\rm q}  $ are
relatively prime.  The finite zero dimensional actual little groups are
the groups which fall into the class $ {\rm E}_{1}. \;$
\end{thrm}

\section{Description of the Finite Zero Dimensional Actual Little Groups
} \label{bcnnnmcmce} \normalfont
\label{s5}

\noindent
%As it was explained in section
%\ref{lyigbhjbukygfi7k6tykj}
The finite zero dimensional subgroups of $\rm{\;SO(2)\times
SO(2)}\;$ have either one or two generators and have been given in detail in \cite{mel}.
We will use this information to describe in detail the finite zero
dimensional actual little groups.

\subsection{Cyclic Actual Little Groups}
\noindent
Firstly we find the actual cyclic little groups.
The  result is given in Proposition \ref{yhyhyuhvbhygygku}.
%Firstly, by using
%subsection \ref{uhulihnbjkhbguybv} (and in particular Proposition
%\ref{lyuigy7tfduujtrxduj65} and Theorem
%\ref{kyugfkhbvkyutgol87yhgkl})
%we give in the following
%Firstly we deal with the cyclic subgroups.
%In Proposition \ref{yhyhyuhvbhygygku} an explicit description of the cyclic subgroups of
%$
%\;K=
%SO(2)\times SO(2)\;$ which fall into the class $\;{\rm
%E}_{1}$ is given.
%Firstly, by using
%subsection \ref{uhulihnbjkhbguybv} (and in particular Proposition
For this purpose we recall some relevant results,
Proposition \ref{lyuigy7tfduujtrxduj65} and Theorem \ref{kyugfkhbvkyutgol87yhgkl},
proved in \cite{mel}.
%Furthermore, we need to recall another basic result from \cite{mel}:
%%which contain the element $(-I,-I)$.
%%An examination of
The element $(-I,-I)$ can only be contained in
the cyclic subgroups of $\;{\rm C}_{2^{a}}\times
{\rm C}_{2^{\beta}}$, and the actual cyclic little groups,
%$\;K=SO(2)\times SO(2)\;$
%which are displayed in subsection \ref{ikytfygbol87ytol78yhlohug}
%in Appendix \ref{lyh,hv,hvmjtdcnmc} shows that the cyclic groups
the  cyclic subgroups of $ \rm{C_{n}} \times  \rm{C_{m}} $
%$\;K=SO(2)\times SO(2)\;$
which fall into the class $\;{\rm E}_{1},$
are precisely those
which contain the element $(-I,-I);$ the group $D$ is not cyclic.
The key result in \cite{mel} we need to give an explicit description of
the finite cyclic actual little groups is the following:
%leads to this conclusion
%is that
The {\it only} cyclic subgroups
of $\;{\rm C}_{2^{a}}\times {\rm C}_{2^{\beta}}\;$ which contain
the element $(-I,-I)$ are the groups
%one of the following cyclic subgroups of
%$\;{\rm C}_{2^{a_{1}}} \times {\rm C}_{2^{\beta_{1}}}\;$
$ \left ( R \left ( \left ( \frac {2 \pi}{2^ {{\rm k}}} {\rm r}
\right )  i_{1} \right ) ,
 R  \left ( \frac {2 \pi}{2^{{\rm k}}}    i_{1} \right )
\right), $ where $\;1 \leq {\rm k} \leq {\rm min}(a,\beta),\;$
$\; {\rm r} \;$ parametrises the groups and takes values in the
set $\; \left \{ 1,2,...,2^{{\rm k}}-1 \right \}- \left \{2,2
\cdot 2,...,(2^{{\rm k}-1}-1)2 \right \}$. Now we are ready to
give an explicit description of the finite cyclic actual little
groups. Firstly we recall some useful results, Proposition
\ref{lyuigy7tfduujtrxduj65} and Theorem \ref{kyugfkhbvkyutgol87yhgkl},
proved in \cite{mel}.
%these groups is given
%in the following Proposition.

\begin{prop}
\label{lyuigy7tfduujtrxduj65}
Let \ ${\rm C}_{\rm n} \times {\rm C}_{\rm m}$ \ be the direct product of
the cyclic groups of finite order
\ ${\rm C}_{\rm n}$ \ and \ $ {\rm C}_{\rm m}$. Let \
${\rm n}={\rm p}_{1}^{a_{1}} \cdot $
${\rm p}_{2}^{a_{2}} \cdot \cdot \ \! \cdot $
${\rm p}_{\rm s}^{a_{\rm s}} $ \ and \
${\rm m}={\rm p}_{1}^{\beta_{1}} \cdot $
${\rm p}_{2}^{\beta_{2}} \cdot \cdot \ \! \cdot $
${\rm p}_{\rm s}^{\beta_{\rm s}} $ \ be the prime decomposition of the
integers \ ${\rm n}$ \ and \ ${\rm m}$, \, i.e.,\ $ {\rm p}_{\rm i},$
\  { \rm i=1,2,...,s,}  are distinct prime numbers and
${ a_{\rm i}, \ \beta_{\rm i}}$ \  are non$-$negative integers.
Then we have the following
\begin{enumerate}
\item
{A  group
\be
\label{xesedrtybbnnvvvfgy}
{\mathcal C}=
{\rm A}_{1} \times {\rm A}_{2} \times {\rm A}_{3}
\times ... \times
{\rm A}_{\rm s},
\ee
where ${\rm A}_{\rm i}$ is a cyclic subgroup,
not necessarily different from the identity element,
of
$\;{\rm C}_{{\rm p}_{\rm i}^{a_{\rm i}}}
\times {\rm C}_{{\rm p}_{\rm i}^{\beta_{\rm i}}},\;$ i=1,2,...,{\rm s},
is a cyclic subgroup of
$\;{\rm C}_{\rm n} \times {\rm C}_{\rm m}.\;$}
\item
{Every cyclic subgroup $\;{\mathcal C}\;$ of
$\;{\rm C}_{\rm n} \times {\rm C}_{\rm m}\;$ is of the form
$$
{\mathcal C}=
{\rm A}_{1} \times {\rm A}_{2} \times {\rm A}_{3}
\times ... \times
{\rm A}_{\rm s},
$$
where ${\rm A}_{\rm i}$ is a cyclic subgroup,
not necessarily different from the identity element,
of
$\;{\rm C}_{{\rm p}_{\rm i}^{a_{\rm i}}}
\times {\rm C}_{{\rm p}_{\rm i}^{\beta_{\rm i}}},\;$ i=1,2,...,{\rm s}.}
\item
{For every cyclic group $\;{\mathcal C}\;$
of $\;{\rm C}_{\rm n} \times {\rm C}_{\rm m}\;$ the expression
${\rm (\ref{xesedrtybbnnvvvfgy})}$  is unique.}
\end{enumerate}
A generator of the cyclic subgroup $\mathcal C$ is given by
\be
\label{ljihgiyutgoo87b}
(x^{\mathcal A},y^{\mathcal B}),
\ee
where, $\;x\;$ is a generator of $\;{\rm C}_{{\rm n}},\;$
$\;y\;$ is a generator of  $\;{\rm C}_{{\rm m}},\;$
\begin{eqnarray}
\mathcal A & = & \sum_{i=1}^{\nu} {\rm r}_{i}p_{i}^{{\rm
a}_{i}-{\rm k}_{i}} ({\rm n}/p_{i}^{{\rm a}_{i}}) +
%{\rm r}_{1}p_{1}^{{\rm a}_{1}-{\rm k}_{1}} ({\rm n}/p_{1}^{{\rm
%a}_{1}}) + {\rm r}_{2}p_{2}^{{\rm a}_{2}-{\rm k}_{2}} ({\rm
%n}/p_{2}^{{\rm a}_{2}}) +...+ {\rm r}_{\nu}p_{\nu}^{{\rm
%a}_{\nu}-{\rm k}_{\nu}}
%({\rm n}/p_{\nu}^{{\rm a}_{\nu}})+
%\nonumber \\
%&&
\sum_{i=\nu+1}^{\nu+\chi} p_{i}^{{\rm a}_{i}-{\rm k}_{i}} ({\rm
n}/p_{i}^{{\rm a}_{i}}) +
%p_{\nu+1}^{{\rm a}_{\nu+1}-{\rm k}_{\nu+1}}
%({\rm n}/p_{\nu+1}^{{\rm a}_{\nu+1}})+
%p_{\nu+2}^{{\rm a}_{\nu+2}-{\rm k}_{\nu+2}}
%({\rm n}/p_{\nu+2}^{{\rm a}_{\nu+2}})
%+...+
%p_{\nu+\chi}^{{\rm a}_{\nu+\chi}-{\rm k}_{\nu+\chi}}
%({\rm n}/p_{\nu+\chi}^{{\rm a}_{\nu+\chi}})
 %\nonumber \\
%&&
 \sum_{i=\nu+\chi+1}^{\nu+\chi+\tau} {\rm j}_{i}  ({\rm
n}/p_{i}^{{\rm a}_{i}}) +
%{\rm j}_{\nu+\chi+1} ({\rm n}/p_{\nu+\chi+1}^{{\rm
%a}_{\nu+\chi+1}}) +{\rm j}_{\nu+\chi+2} ({\rm
%n}/p_{\nu+\chi+2}^{{\rm a}_{\nu+\chi+2}}) +...+{\rm
%j}_{\nu+\chi+\tau} ({\rm n}/p_{\nu+\chi+\tau}^{{\rm
%a}_{\nu+\chi+\tau}})+
\nonumber \\
\label{bycoxenasha} &&
\sum_{i=\nu+\chi+\tau+1}^{\nu+\chi+\tau+\psi} p_{i}^{{\rm
a}_{i}-{\rm k}_{i}} ({\rm n}/p_{i}^{{\rm a}_{i}})
%p_{\nu+\chi+\tau+1}^{{\rm a}_{\nu+\chi+\tau+1}-{\rm k}_{\nu+\chi+\tau+1}}
%({\rm n}/p_{\nu+\chi+\tau+1}^{{\rm a}_{\nu+\chi+\tau+1}})+
%p_{\nu+\chi+\tau+2}^{{\rm a}_{\nu+\chi+\tau+2}-{\rm k}_{\nu+\chi+\tau+2}}
%({\rm n}/p_{\nu+\chi+\tau+2}^{{\rm a}_{\nu+\chi+\tau+2}})
%+...+
%\nonumber \\
%\label{bycoxenasha} &&
%\sum_{i=\nu+1}^{\nu+\chi} p_{i}^{{\rm
%a}_{i}-{\rm k}_{i}} ({\rm n}/p_{i}^{{\rm a}_{i}}) +
%p_{\nu+\chi+\tau+\psi}^{{\rm a}_{\nu+\chi+\tau+\psi} -{\rm
%k}_{\nu+\chi+\tau+\psi}} ({\rm n}/p_{\nu+\chi+\tau+\psi}^{{\rm
%a}_{\nu+\chi+\tau+\psi}})
, \qquad{\rm and}
\end{eqnarray}
\begin{eqnarray}
\mathcal B & = & \sum_{i=1}^{\nu} p_{i}^{{\rm b}_{i}-{\rm k}_{i}}
({\rm m}/p_{i}^{{\rm b}_{i}}) +
%p_{1}^{{\rm b}_{1}-{\rm k}_{1}}
%({\rm m}/p_{1}^{{\rm b}_{1}}) + p_{2}^{{\rm b}_{2}-{\rm k}_{2}}
%({\rm m}/p_{2}^{{\rm b}_{2}}) +...+ p_{\nu}^{{\rm b}_{\nu}-{\rm
%k}_{\nu}} ({\rm m}/p_{\nu}^{{\rm b}_{\nu}}) +
%\nonumber \\
%&&
\sum_{i=\nu+1}^{\nu+\chi} \rho_{i}p_{i}^{{\rm b}_{i}-{\rm
k}_{i}+1} ({\rm m}/p_{i}^{{\rm b}_{i}}) +
%\rho_{\nu+1}p_{\nu+1}^{{\rm b}_{\nu+1}-{\rm k}_{\nu+1}+1} ({\rm
%m}/p_{\nu+1}^{{\rm b}_{\nu+1}})+ \rho_{\nu+2}p_{\nu+2}^{{\rm
%b}_{\nu+2}-{\rm k}_{\nu+2}+2} ({\rm m}/p_{\nu+2}^{{\rm
%b}_{\nu+2}})
%+...+ \nonumber \\
%&&\rho_{\nu+\chi}p_{\nu+\chi}
%^{{\rm b}_{\nu+\chi}-{\rm k}_{\nu+\chi}+1}
%({\rm m}/p_{\nu+\chi}^{{\rm b}_{\nu+\chi}})
 %\nonumber \\
%&&
\sum_{i=\nu+\chi+1}^{\nu+\chi+\tau} p_{i}^{{\rm b}_{i}-{\rm
k}_{i}} ({\rm m}/p_{i}^{{\rm b}_{i}}) +
%p_{\nu+\chi+1}^{{\rm b}_{\nu+\chi+1}-{\rm k}_{\nu+\chi+1}} ({\rm
%m}/p_{\nu+\chi+1}^{{\rm b}_{\nu+\chi+1}})+ p_{\nu+\chi+2}^{{\rm
%b}_{\nu+\chi+2}-{\rm k}_{\nu+\chi+2}} ({\rm
%m}/p_{\nu+\chi+2}^{{\rm b}_{\nu+\chi+2}})
%+...+ \nonumber \\
%&&p_{\nu+\chi+\tau}^{{\rm b}_{\nu+\chi+\tau}-{\rm k}_{\nu+\chi+\tau}}
%({\rm m}/p_{\nu+\chi+\tau}^{{\rm b}_{\nu+\chi+\tau}})+
\nonumber \\
\label{bnkyuffdcikp;jp;0i9} &&
\sum_{i=\nu+\chi+\tau+1}^{\nu+\chi+\tau+\psi} {\rm j}_{i} ({\rm
m}/p_{i}^{{\rm b}_{i}}). \quad
%{\rm j}_{\nu+\chi+\tau+1} ({\rm m}/p_{\nu+\chi+\tau+1}^{{\rm
%b}_{\nu+\chi+\tau+1}})+ {\rm j}_{\nu+\chi+\tau+2} ({\rm
%m}/p_{\nu+\chi+\tau+2}^{{\rm b}_{\nu+\chi+\tau+2}})
%+...+ \nonumber \\
%&&{\rm j}_{\nu+\chi+\tau+\psi}
%({\rm m}/p_{\nu+\chi+\tau+\psi}^{{\rm b}_{\nu+\chi+\tau+\psi}}) \quad .
\end{eqnarray}
The order $\;|\mathcal C|\;$ of the group $\;\mathcal C\;$ is given by
\begin{eqnarray}
\label{lyuidfytde5rikyto9} |\mathcal C|&=&
\prod_{i=1}^{\nu+\chi+\tau+\psi}p_{i}^{{\rm k}_{i}}. \quad
%p_{1}^{{\rm k}_{1}}\cdot p_{2}^{{\rm k}_{2}} \cdot \cdot \cdot
%p_{\nu}^{{\rm k}_{\nu}} \cdot p_{\nu+1}^{{\rm
%k}_{\nu+1}}p_{\nu+2}^{{\rm k}_{\nu+2}}
%\cdot\cdot\cdot p_{\nu+\chi}^{{\rm k}_{\nu+\chi}}\cdot \nonumber \\
%\label{lyuidfytde5rikyto9}
%&&p_{\nu+\chi+1}^{{\rm k}_{\nu+\chi+1}}\cdot
%p_{\nu+\chi+2}^{{\rm k}_{\nu+\chi+2}}\cdot\cdot\cdot
%p_{\nu+\chi+\tau}^{{\rm k}_{\nu+\chi+\tau}}\cdot
%p_{\nu+\chi+\tau+1}^{{\rm k}_{\nu+\chi+\tau+1}}\cdot
%p_{\nu+\chi+\tau+2}^{{\rm k}_{\nu+\chi+\tau+2}}\cdot\cdot\cdot
%p_{\nu+\chi+\tau+\psi}^{{\rm k}_{\nu+\chi+\tau+\psi}}.
\end{eqnarray}
The non$-$negative integers
$\; \nu,\chi,\tau,\psi \;$ are such that $\; \nu+\chi+\tau+\psi \leq {\rm s}.
%\;
$
Moreover,
$\;
(p_{\rm i}^{{\rm a}_{\rm i}},p_{\rm i}^{{\rm b}_{{\rm i}}})=
\mathcal P({\rm p}_{{\rm i}}^{a_{{\rm i}}},{\rm p}_{{\rm i}}
^{\beta_{{\rm i}}}),\; \;{\rm i}=1,2,...,{\rm s},\;$
for some permutation  $\mathcal P$ of the $s$ pairs of numbers
$\;({\rm p}_{1}^{a_{1}},{\rm p}_{1}^{\beta_{1}}),
({\rm p}_{2}^{a_{2}},{\rm p}_{2}^{\beta_{2}}),...,
({\rm p}_{{\rm s}}^{a_{{\rm s}}},{\rm p}_{{\rm s}}^{\beta_{{\rm s}}}).\;$
Furthermore,
$\;{\rm r}_{\rm i} \in \left \{
0,1,2,...,p_{\rm i}^{{\rm k}_{\rm i}}-1 \right \},\; \;
{\rm i}\in \left \{1,2,...,\nu \right \},
\;$ and,
$\;{\rm j}_{\sigma} \in \left \{
0,1,2,...,p_{\sigma}^{{\rm a}_{\sigma}}-1 \right \},
\;
{\rm a}_{\sigma}<{\rm k}_{\sigma}\leq{\rm b}_{\sigma},\;
\;$
$
{\sigma}\in \left \{\nu+\chi+1,\nu+\chi+2,...,\right .$
\newline
$\left .
\nu+\chi+\tau \right \}.
\;$
Finally,
$\;{\rho}_{\rm q} ~\in \left \{
0,1,2,...,p_{\rm q}^{{\rm k}_{\rm q}-1}-1 \right \},\; \;
{\rm q}\in \left \{\nu+1,\nu+2,...,\nu+\chi \right \},
\;$ and,
$\;{\rm j}_{\theta} \in \left \{
0,1,2,...,p_{\theta}^{{\rm b}_{\theta}}-1 \right \},
\;
{\rm a}_{\theta} \geq {\rm k}_{\theta} >
{\rm b}_{\theta}, \; \;
{\theta} \in \left \{\nu+\chi+\tau+1,\nu+\chi+\tau+2,...,
\nu+\chi+\tau+ \right . $
\newline
$\left . \psi \right \}. \;$
\normalfont
\end{prop}

\noindent
One crucial feature of the cyclic subgroups $\;\mathcal C\;$ 
of $\;{\rm C}_{\rm n} \times {\rm C}_{\rm m}\;$  is that for 
\it every \normalfont subgroup$ \;\mathcal
C\;$ the expression (\ref{xesedrtybbnnvvvfgy}) is unique. This was proved in \cite{mel} with the use of 
Sylow's Second Theorem. 
It is noteworthy that one can use
instead of  Sylow's Second Theorem more
elementary group theory in order to prove the uniqueness of (\ref{xesedrtybbnnvvvfgy}).
%In fact in the case of cyclic subgroups $\;\mathcal C\;$
%of $\;{\rm C}_{\rm n} \times {\rm C}_{\rm m}\;$  one can use
%instead of  Sylow's Second Theorem more
%elementary
%and transparent
%group theory in order to prove that for every subgroup $\;\mathcal
%C\;$ the expression (\ref{xesedrtybbnnvvvfgy}) is unique. 
Let
$\;\gamma\;$ be a generator of the cyclic group $\;\mathcal C=
{\rm C}_{{\rm p}_{1}^{\rm k_{1}}} \times {\rm C}_{{\rm p}_{2}^{\rm
k_{2}}} \times {\rm C}_{{\rm p}_{3}^{\rm k_{3}}} \times \ ... \
\times {\rm C}_{{\rm p}_{\rm s}^{\rm k_{\rm s}}}. $ Then the
element $\; \omega=\gamma^{ {\rm p}_{2}^{\rm k_{2}} {\rm
p}_{3}^{\rm k_{3}}\ ... \ {\rm p}_{\rm s}^{\rm k_{\rm s}}}\;$
satisfies $\;\omega^{{\rm p}_{1}^{\rm k_{1}}}= (\gamma^{ {\rm
p}_{2}^{\rm k_{2}} {\rm p}_{3}^{\rm k_{3}}\ ... \ {\rm p}_{\rm
s}^{\rm k_{\rm s}}} ) ^{{\rm p}_{1}^{\rm k_{1}}}=1$. If the
element $\;\omega\;$ had order $\;{\rm p}_{1}^{\rm m_{1}},\;$
where $\;{\rm m}_{1}<{\rm k}_{1},\;$ then we would have $\gamma^{
{\rm p}_{2}^{\rm k_{2}} {\rm p}_{3}^{\rm k_{3}}\ ... \ {\rm
p}_{\rm s}^{\rm k_{\rm s}}  {\rm p}_{1}^{\rm m_{1}}}=1$. But the
order of $\; \gamma \;$ equals to the order of $\; \mathcal C. \;$
So we get a contradiction and therefore the order of $\; \omega
\;$ is $\;{\rm p}_{1}^{\rm k_{1}}.\;$ Thus $\; \omega \;$
generates $\;{\rm C}_{{\rm p}_{1}^{\rm k_{1}}}.$ Choose now any
element $\;x\;$ of order $\;{\rm p}_{1}^{{\rm k}_{1}}.\;$ Since
$\;\gamma \;$ generates the whole group $\;\mathcal C\;$ we have
$\;x=\gamma^{t}\;$ for some integer $\;t.\;$ Moreover, since
$\;x\;$ generates $\;{\rm C}_{{\rm p}_{1}^{\rm k_{1}}}\;$ we have
$\;\gamma^{t{\rm p}_{1}^{{\rm k}_{1}}}=1.\;$ Therefore, the order
of $\;\gamma\;$ divides $\;t  {\rm p}_{1} ^{{\rm k}_{1}}\;$
,i.e.,$\;{\rm p}_{1}^{{\rm k}_{1}}...{\rm p}_{{\rm s}}^{{\rm
k}_{{\rm s}}}\;$ divides $\;t{\rm p}_{1}^{{\rm k}_{1}}.\;$
Consequently, $\;{\rm p}_{2}^{{\rm k}_{2}}...{\rm p}_{{\rm
s}}^{{\rm k}_{{\rm s}}},\;$ divides $\;t.\;$ Thus $\;t=u{\rm
p}_{2}^{{\rm k}_{2}}...{\rm p}_{{\rm s}}^{{\rm k}_{{\rm s}}},\;$
for some integer u. We conclude that $x=\gamma^{t}=(\gamma^{{\rm
p}_{2}^{{\rm k}_{2}}... {\rm p}_{{\rm s}}^{{\rm k}_{{\rm
s}}}})^{u}=\omega^{u}.$ Therefore the element $\;x\;$ belongs to
the specific copy of $\;{\rm C}_{{\rm p}_{1}^{\rm k_{1}}}\;$ which
is generated by $\; \omega. \;$ We conclude that for every cyclic
subgroup $\;\mathcal C\;$ of $\;{\rm C}_{\rm n} \times {\rm
C}_{\rm m}\;$ the expression (\ref{xesedrtybbnnvvvfgy}) is unique.

\vspace{0.3cm}

\noindent
For the purposes of our study we also recall from \cite{mel}
%rewrite a 
that a cyclic
subgroup $\;\mathcal C\;$ of
$\;{\rm C}_{\rm n} \times {\rm C}_{\rm m}\;$ can be conveniently rewritten
as a subgroup of ${\rm S}{\rm O}(2) \times {\rm S}{\rm O}(2)\;$.

\begin{thrm}
\label{kyugfkhbvkyutgol87yhgkl}
Let $\;{\rm n}\;$ and $\;{\rm m}\;$ be {\it any} non-negative integers.
Then all the finite cyclic subgroups of
$\;{\rm S}{\rm O}(2) \times {\rm S}{\rm O}(2)\;$ are given by
\begin{eqnarray}
\mathcal{C} & = &
\left ( R \left ( \left ( \frac {2 \pi}{{\rm n}} {\mathcal A} \right )
 i
\right ) ,
 R  \left ( \left ( \left .
\frac {2 \pi}{{\rm m}}  {\mathcal B} \right ) i
\right ) \right . \right ),
\end{eqnarray}
where the expressions  $\; \mathcal A \;$ and $\; \mathcal B \;$ are given
in   (\ref{bycoxenasha}) and (\ref{bnkyuffdcikp;jp;0i9}) correspondingly.
The meaning and the ranges of the parameters
appearing in these
expressions  are displayed in Proposition
\ref{lyuigy7tfduujtrxduj65}.
For each specific subgroup
these parameters take specific values. Different values of the parameters
correspond to different subgroups and vice versa.
The order of the group $\; \mathcal C \;$ is given by
(\ref{lyuidfytde5rikyto9}).
\end{thrm}

\noindent 
We are ready now to use Proposition \ref{lyuigy7tfduujtrxduj65} 
and Theorem \ref{kyugfkhbvkyutgol87yhgkl}
to describe in detail the cyclic actual little groups. 
Their explicit description is given in the Proposition that follows.
In Theorem \ref{T7}  the cyclic actual little groups are rewritten as subgroups of     
${\rm S}{\rm O}(2) \times {\rm S}{\rm O}(2).\;$

\begin{prop}
\label{yhyhyuhvbhygygku}
 Let ${\rm n}$ and ${\rm m}$ be any
positive even numbers. Then
\begin{equation}
{\rm C}_{\rm n} \times {\rm C}_{\rm m}= ({\rm C}_{2^{a_{1}}}
\times {\rm C}_{2^{\beta_{1}}}) \times ({\rm C}_{{\rm
p}_{2}^{a_{2}}} \times {\rm C}_{{\rm p}_{2}^{\beta_{2}}}) \times
({\rm C}_{{\rm p}_{3}^{a_{3}}} \times {\rm C}_{{\rm
p}_{3}^{\beta_{3}}}) \times \ ... \ \times ({\rm C}_{{\rm p}_{\rm
s}^{a_{\rm s}}} \times {\rm C}_{{\rm p}_{\rm s}^{\beta_{\rm s}}})
\end{equation}
where $\;a_{1} \geq 1, \;$ $\; \beta_{1}  \geq 1, \;$ $\; {\rm
p}_{2},{\rm p}_{3},...,{\rm p}_{\rm s} \; $ are odd primes,
$\;a_{\rm i} \geq 0, \;$ and, $\; \beta_{\rm i} \geq 0,\;$
$\;{\rm i} \in \left \{ 2,3,...,{\rm s} \right \}. \;$ Every
cyclic subgroup of $\;K=SO(2)\times SO(2)\;$ which contains the
element $\;(-I,-I)\;$ is written uniquely in the form
$$
\mathcal C={\rm A}_{1} \times {\rm A}_{2} \times {\rm A}_{3}\times
... \times {\rm A}_{\rm s},
$$
where $\;{\rm A}_{{\rm i}}\;$ is a cyclic subgroup, not
necessarily different from the identity element, of $\; {\rm
C}_{{\rm p}_{{\rm i}}^{a_{\rm i}}} \times {\rm C}_{{\rm p}_{{\rm
i}}^{\beta_{\rm i}}}, \;$ $\;{\rm i} \in \left \{ 2,...,{\rm s}
\right \}, \;$ and $\;{\rm A}_{1},\;$ $\;({\rm p}_{1}=2),\;$ is
restricted to be one of the following cyclic subgroups of $\;{\rm
C}_{2^{a_{1}}} \times {\rm C}_{2^{\beta_{1}}}\;$
$$
\left ( R \left ( \left ( \frac {2 \pi}{2^ {{\rm k}_{1}}} {\rm r}
\right )  i_{1} \right ) ,
 R  \left ( \frac {2 \pi}{2^{{\rm k}_{1}}}    i_{1} \right )
\right ),
$$
where $\;1 \leq {\rm k}_{1} \leq {\rm min}(a_{1},\beta_{1}),\;$
$\; {\rm r} \;$ parametrises the groups and takes values in the
set $\; \left \{ 1,2,...,2^{{\rm k}_{1}}-1 \right \}- \left \{2,2
\cdot 2,...,(2^{{\rm k}_{1}-1}-1)2 \right \} \;$ and $\;i_{1}\;$
enumerates the elements of each group and takes values in the set
$\; \left \{0,1,2,...,2^{{\rm k}_{1}}-1 \right \}.$
A generator
of $\; \mathcal C \;$ is given by
(\ref{ljihgiyutgoo87b}), where, $\; (p_{\rm i}^{{\rm a}_{\rm
i}},p_{\rm i}^{{\rm b}_{{\rm i}}})= \mathcal P({\rm p}_{{\rm
i}}^{a_{{\rm i}}},{\rm p}_{{\rm i}} ^{\beta_{{\rm i}}})\;,\;{\rm
i} \in \left \{ 1,2,...,{\rm s} \right \},\;$ for some
permutation  $\mathcal P$ of the $s$ pairs of numbers $\;({\rm
p}_{1}^{a_{1}},{\rm p}_{1}^{\beta_{1}}), ({\rm p}_{2}^{a_{2}},{\rm
p}_{2}^{\beta_{2}}),..., ({\rm p}_{{\rm s}}^{a_{{\rm s}}},{\rm
p}_{{\rm s}}^{\beta_{{\rm s}}})\;$, and where, $\; {\rm
n},\;$ $\;{\rm m} \;$ are positive even numbers. In
(\ref{ljihgiyutgoo87b})  one of the primes
$\;p_{1},p_{2},...,p_{\nu}\;$ is the prime number $\rm 2$. If say,
$\;p_{t}=2,\;$ $\; t \in \left \{ 1,2,...,\nu \right \}, \;$ then
$\;{\rm r}_{t} \in \left \{1,2,...,2^{{\rm k}_{t}}-1 \right \}-
\left \{2,2 \cdot 2,...,(2^{{\rm k}_{t}-1}-1)2 \right \}. $ The
rest of the indices $\;{\rm r}_{d}\;$, $\;d \in \left \{
1,2,...,\nu \right \}- \left \{ t \right \}, \;$ take values in
the sets $\left \{0,1,2,...,p_{d}^{{\rm k}_{d}}-1 \right \}. $ The
other indices which appear in (\ref{ljihgiyutgoo87b}) take
values in the sets which are displayed in Proposition
(\ref{lyuigy7tfduujtrxduj65}). Some, or in fact all the exponents
$\;{\rm a}_{\rm i} \;$ and $\;{\rm b}_{\rm i} ,\;$ $\;{\rm i} \in
\left \{ 1,2,3,...,{\rm s} \right \} - \left \{ t \right \}, \;$
which appear in (\ref{ljihgiyutgoo87b}) can be equal to zero.
\end{prop}

\begin{thrm}
\label{T7}
The cyclic subgroup $\; \mathcal C \;$ can be written as
$$
\left ( R \left ( \left ( \frac {2 \pi}{{\rm n}} \mathcal A \right
)  i \right ) ,
 R  \left ( \left ( \frac {2 \pi}{{\rm m}} \mathcal B \right )   i \right )
\right ),
$$
where $\; \mathcal A \;$ and $\; \mathcal B \;$ are given
respectively by the expressions (\ref{bycoxenasha}) and
(\ref{bnkyuffdcikp;jp;0i9}). The ranges of the parameters 
which appear in $\; \mathcal A \;$ and $\; \mathcal B \;$
are specified in Proposition (\ref{yhyhyuhvbhygygku}).
%\label{yhyhyuhvbhygygku}
%In expressions (\ref{bycoxenasha})
%and (\ref{bnkyuffdcikp;jp;0i9}) $\; {\rm n} $ and $\;{\rm m} \;$ are
%positive even numbers and one of the primes
%$\;p_{1},p_{2},...,p_{\nu}\;$ is the prime number $ \rm 2$.
%If say,
%$\;p_{t}=2,\;$ $\; t \in \left \{ 1,2,...,\nu \right \}, \;$ then
%$\;{\rm r}_{t} \in \left \{1,2,...,2^{{\rm k}_{t}}-1 \right \}-
%\left \{2,2 \cdot 2,...,(2^{{\rm k}_{t}-1}-1)2 \right \}. $ The
%rest of the indices $\;{\rm r}_{d}\;$, $\;d \in \left \{
%1,2,...,\nu \right \}- \left \{ t \right \}, \;$ take values in
%the sets $\left \{0,1,2,...,p_{d}^{{\rm k}_{d}}-1 \right \}. $ The
%other indices which appear in (\ref{bycoxenasha}) and
%(\ref{bnkyuffdcikp;jp;0i9}) take values in the sets which are
%displayed in Proposition (\ref{lyuigy7tfduujtrxduj65}). In
%expressions (\ref{bycoxenasha}) and (\ref{bnkyuffdcikp;jp;0i9})
%some of the exponents
%%primes
%%In (\ref{nblyuitgo;hjlkgf7o}) some
%%of the exponents
%$\;{\rm a}_{\rm i} \;$ and $\;{\rm b}_{\rm i} ,\;$ $\;{\rm i} \in
%\left \{ 1,2,3,...,{\rm s} \right \} - \left \{ t \right \}, \;$
%or in fact all of them, can be equal to zero.
\end{thrm}

\subsection{Actual Little Groups with Two Generators }
%Now, by using subsection \ref{kiyufghgvjdfujtrdcyjtgcvj} (and in
%particular Proposition \ref{kittttttygy5ygyk5j55} and Theorem
%\ref{uyghkgvyjtrfitfvkyjtggfi7k6g}) we describe explicitly the
%finite actual little groups with two generators. 

We describe now explicitly in Proposition \ref{kloikloikloikloiklo} the
finite actual little groups with two generators. 
The key
observation here is that there is no subgroup of $\;{\rm
C}_{2^{a}}\times {\rm C}_{2^{\beta}}\;$ with two generators which
falls into the class $\;{\rm E}_{1}$ \cite{mel}. 
%This can be easily
%verified by examining these subgroups which are displayed in
%subsection \ref{mhhhdhdhhdgtter}  in Appendix
%\ref{lyh,hv,hvmjtdcnmc}. 
Therefore, as in the case of the cyclic
actual groups, the subgroups of $\;{\rm C}_{2^{a}}\times {\rm
C}_{2^{\beta}}\;$ are restricted to be one of the groups $ \left (
R \left ( \left ( \frac {2 \pi}{2^ {{\rm k}}} {\rm r} \right )
i_{1} \right ) ,
 R  \left ( \frac {2 \pi}{2^{{\rm k}}}    i_{1} \right )
\right ), $ where, $\;1 \leq {\rm k} \leq {\rm min}(a,\beta),\;$
$\; {\rm r} \;$ parametrises the groups and takes values in the
set $\; \left \{ 1,2,...,2^{{\rm k}}-1 \right \}- \left \{2,2
\cdot 2,...,(2^{{\rm k}-1}-1)2 \right \}$). 
We need to recall first some relevant results, Proposition 
\ref{kittttttygy5ygyk5j55} and Theorem \ref{uyghkgvyjtrfitfvkyjtggfi7k6g}, from \cite{mel}.

%The details are as
%follows:

%the subgroups of
%$\;K=SO(2)\times SO(2)\;$
%with two generators which contain the element $(-I,-I)$.
%in Appendix
%\ref{lyh,hv,hvmjtdcnmc}
%shows that the cyclic groups
%%cyclic subgroups of $\;K=SO(2)\times SO(2)\;$
%which fall into the class $\;{\rm E}_{1}\;$ are precisely those which contain the
%element $(-I,-I)$.
%The key observation which leads to this conclusion is that the {\it only} cyclic subgroups of
%$\;{\rm C}_{2^{a}}\times {\rm C}_{2^{\beta}}\;$
%\label{ik76tyfvnmgbfyhtrdtrcfhn}

\begin{prop}
\label{kittttttygy5ygyk5j55}
Let \ ${\rm C}_{\rm n} \times {\rm C}_{\rm m}$ \ be the direct product of
the cyclic groups of finite order
\ ${\rm C}_{\rm n}$ \ and \ $ {\rm C}_{\rm m}$. Let \
${\rm n}={\rm p}_{1}^{a_{1}} \cdot $
${\rm p}_{2}^{a_{2}} \cdot \cdot \ \! \cdot $
${\rm p}_{\rm s}^{a_{\rm s}} $ \ and \
${\rm m}={\rm p}_{1}^{\beta_{1}} \cdot $
${\rm p}_{2}^{\beta_{2}} \cdot \cdot \ \! \cdot $
${\rm p}_{\rm s}^{\beta_{\rm s}} $ \ be the prime decomposition of the
integers \ ${\rm n}$ \ and \ ${\rm m}$, \, i.e.,\ $ {\rm p}_{\rm i}$,
\  { \rm i=1,2,...,s}, are distinct prime numbers and
${ a_{\rm i}, \ \beta_{\rm i}}$ \  are non$-$negative integers.
Then we have the following
\begin{enumerate}
\item
{A  group
\be
\label{boyuifgity7fhjtrsdyu}
{\mathcal C}=
{\rm A}_{1} \times {\rm A}_{2} \times {\rm A}_{3}
\times ... \times
{\rm A}_{\rm s}
\ee
where ${\rm A}_{\rm i}$ is a subgroup,
not necessarily different from the identity element,
of
$\;{\rm C}_{{\rm p}_{\rm i}^{a_{\rm i}}}
\times {\rm C}_{{\rm p}_{\rm i}^{\beta_{\rm i}}}$, i=1,2,...,{\rm s},
is a  subgroup of
$\;{\rm C}_{\rm n} \times {\rm C}_{\rm m}\;$
with two generators if at least one of the
${\rm A}_{\rm i}$, { \rm i=1,2,...,s}, has two generators.
}
\item
{Every subgroup $\;{\mathcal C}\;$ of
$\;{\rm C}_{\rm n} \times {\rm C}_{\rm m}\;$ with two generators
is of the form
$$
{\mathcal C}=
{\rm A}_{1} \times {\rm A}_{2} \times {\rm A}_{3}
\times ... \times
{\rm A}_{\rm s},
$$
where ${\rm A}_{\rm i}$ is a subgroup,
not necessarily different from the identity element,
of
$\;{\rm C}_{{\rm p}_{\rm i}^{a_{\rm i}}}
\times {\rm C}_{{\rm p}_{\rm i}^{\beta_{\rm i}}}$, i=1,2,...,{\rm s},
and, where at least one of the
$\;{\rm A}_{\rm i}$, i=1,2,...,{\rm s}, has two generators.
}
\item
{For every subgroup $\;{\mathcal C}\;$
of $\;{\rm C}_{\rm n} \times {\rm C}_{\rm m}\;$ with two generators
the expression
${\rm (\ref{boyuifgity7fhjtrsdyu})}$  is unique.}
\end{enumerate}
Two  generators of a subgroup $\mathcal C$ of
$\;{\rm C}_{\rm n} \times {\rm C}_{\rm m}\;$ with two generators
are given by
\begin{enumerate}
\item
\be
\label{ku6tr967rfujr6xikt5fklo}
g_{1}=(x^{\mathcal A_{1}},y^{\mathcal B_{1}}),
\ee

where, $\;x\;$ is a generator of $\;{\rm C}_{{\rm n}}$,
$\;y\;$ is a generator of  $\;{\rm C}_{{\rm m}}$,
\begin{eqnarray}
\mathcal A_{1}  & = & \sum_{i=1}^{\nu} {\rm r}_{i}p_{1}^{{\rm
a}_{i}-{\rm k}_{i}} ({\rm n}/p_{i}^{{\rm a}_{i}}) +
%{\rm r}_{1}p_{1}^{{\rm a}_{1}-{\rm k}_{1}}
%({\rm n}/p_{1}^{{\rm a}_{1}}) +
%{\rm r}_{2}p_{2}^{{\rm a}_{2}-{\rm k}_{2}}
%({\rm n}/p_{2}^{{\rm a}_{2}})
%+...+
%{\rm r}_{\nu}p_{\nu}^{{\rm a}_{\nu}-{\rm k}_{\nu}}
%({\rm n}/p_{\nu}^{{\rm a}_{\nu}})+
%\nonumber \\
%&&
\sum_{i=\nu+1}^{\nu+\chi} p_{i}^{{\rm a}_{i}-{\rm k}_{i}} ({\rm
n}/p_{i}^{{\rm a}_{i}})+
%p_{\nu+1}^{{\rm a}_{\nu+1}-{\rm k}_{\nu+1}}
%({\rm n}/p_{\nu+1}^{{\rm a}_{\nu+1}})+
%p_{\nu+2}^{{\rm a}_{\nu+2}-{\rm k}_{\nu+2}}
%({\rm n}/p_{\nu+2}^{{\rm a}_{\nu+2}})
%+...+
%p_{\nu+\chi}^{{\rm a}_{\nu+\chi}-{\rm k}_{\nu+\chi}}
%({\rm n}/p_{\nu+\chi}^{{\rm a}_{\nu+\chi}}) +
 %\nonumber \\
%&&
\sum_{i=\nu+\chi+1}^{\nu+\chi+\tau} {\rm j}_{i} ({\rm
n}/p_{i}^{{\rm a}_{i}}) +
%{\rm j}_{\nu+\chi+1} ({\rm n}/p_{\nu+\chi+1}^{{\rm
%a}_{\nu+\chi+1}}) +{\rm j}_{\nu+\chi+2} ({\rm
%n}/p_{\nu+\chi+2}^{{\rm a}_{\nu+\chi+2}}) +...+{\rm
%j}_{\nu+\chi+\tau} ({\rm n}/p_{\nu+\chi+\tau}^{{\rm
%a}_{\nu+\chi+\tau}})+
\nonumber \\
&& \sum_{i=\nu+\chi+\tau+1}^{\nu+\chi+\tau+\psi} p_{i}^{{\rm
a}_{i}-{\rm k}_{i}} ({\rm n}/p_{i}^{{\rm a}_{i}})+
%p_{\nu+\chi+\tau+1}^{{\rm a}_{\nu+\chi+\tau+1}-{\rm k}_{\nu+\chi+\tau+1}}
%({\rm n}/p_{\nu+\chi+\tau+1}^{{\rm a}_{\nu+\chi+\tau+1}})+
%p_{\nu+\chi+\tau+2}^{{\rm a}_{\nu+\chi+\tau+2}-{\rm k}_{\nu+\chi+\tau+2}}
%({\rm n}/p_{\nu+\chi+\tau+2}^{{\rm a}_{\nu+\chi+\tau+2}})
%+...+
%\nonumber \\
%&&p_{\nu+\chi+\tau+\psi}^{{\rm a}_{\nu+\chi+\tau+\psi}
%-{\rm k}_{\nu+\chi+\tau+\psi}}
%({\rm n}/p_{\nu+\chi+\tau+\psi}^{{\rm a}_{\nu+\chi+\tau+\psi}})+
%\nonumber \\
%&&
\sum_{i=\nu+\chi+\tau+\psi+1}^{\nu+\chi+\tau+\psi+\sigma} {\rm
r}_{i} p_{i}^{{\rm a}_{i} -{\rm k}_{i}} ({\rm n}/p_{i} ^{{\rm
a}_{i}})+
%{\rm r}_{\nu+\chi+\tau+\psi+1} p_{\nu+\chi+\tau+\psi+1}^{{\rm
%a}_{\nu+\chi+\tau+\psi+1} -{\rm k}_{\nu+\chi+\tau+\psi+1}} ({\rm
%n}/p_{\nu+\chi+\tau+\psi+1}
%^{{\rm a}_{\nu+\chi+\tau+\psi+1}})+...+ \nonumber \\
%&&{\rm r}_{\nu+\chi+\tau+\psi+\sigma}
%p_{\nu+\chi+\tau+\psi+\sigma}^{{\rm a}_{\nu+\chi+\tau+\psi+\sigma}
%-{\rm k}_{\nu+\chi+\tau+\psi+\sigma}}
%({\rm n}/p_{\nu+\chi+\tau+\psi+\sigma}
%^{{\rm a}_{\nu+\chi+\tau+\psi+\sigma}})+
\nonumber \\
&&
\sum_{i=\nu+\chi+\tau+\psi+\sigma+1}^{\nu+\chi+\tau+\psi+\sigma+\theta}
p_{i} ^{{\rm a}_{i}-{\rm k}_ {i}} ({\rm n}/p_{i} ^{{\rm a}_{i}})+
%p_{\nu+\chi+\tau+\psi+\sigma+1} ^{{\rm
%a}_{\nu+\chi+\tau+\psi+\sigma+1}-{\rm k}_
%{\nu+\chi+\tau+\psi+\sigma+1}} ({\rm
%n}/p_{\nu+\chi+\tau+\psi+\sigma+1}
%^{{\rm a}_{\nu+\chi+\tau+\psi+\sigma+1}})+...+ \nonumber \\
%&&
%p_{\nu+\chi+\tau+\psi+\sigma+\theta}
%^{{\rm a}_{\nu+\chi+\tau+\psi+\sigma+\theta}-{\rm k}_
%{\nu+\chi+\tau+\psi+\sigma+\theta}}
%({\rm n}/p_{\nu+\chi+\tau+\psi+\sigma+\theta}
%^{{\rm a}_{\nu+\chi+\tau+\psi+\sigma+\theta}})+
%\nonumber \\
%&&
\sum_{i=\nu+\chi+\tau+\psi+\sigma+\theta+1}^{\nu+\chi+\tau+\psi+\sigma+\theta+\phi}
{\rm t}_{i} ({\rm n}/p_{i} ^{{\rm a}_{i}})+
%{\rm t}_{\nu+\chi+\tau+\psi+\sigma+\theta+1} ({\rm
%n}/p_{\nu+\chi+\tau+\psi+\sigma+\theta+1}
%^{{\rm a}_{\nu+\chi+\tau+\psi+\sigma+\theta+1}})+...+ \nonumber \\
%&&{\rm t}_{\nu+\chi+\tau+\psi+\sigma+\theta+\phi}
%({\rm n}/p_{\nu+\chi+\tau+\psi+\sigma+\theta+\phi}
%^{{\rm a}_{\nu+\chi+\tau+\psi+\sigma+\theta+\phi}})+
\nonumber \\
&& \label{vyitd6uir5er7i65rfiky}
\sum_{i=\nu+\chi+\tau+\psi+\sigma+\theta+\phi+1}^{\nu+\chi+\tau+\psi+\sigma+\theta+\phi+\xi}
p_{i}^ {{\rm a}_{i}- {\rm k}_{i}} ({\rm n}/p_{i} ^{{\rm a}_{i}})
%p_{\nu+\chi+\tau+\psi+\sigma+\theta+\phi+1}^
%{{\rm a}_{\nu+\chi+\tau+\psi+\sigma+\theta+\phi+1}-
%{\rm k}_{\nu+\chi+\tau+\psi+\sigma+\theta+\phi+1}}
%({\rm n}/p_{\nu+\chi+\tau+\psi+\sigma+\theta+\phi+1}
%^{{\rm a}_{\nu+\chi+\tau+\psi+\sigma+\theta+\phi+1}})+...+ \nonumber \\
%&&
%\label{vyitd6uir5er7i65rfiky}
%p_{\nu+\chi+\tau+\psi+\sigma+\theta+\phi+\xi}^
%{{\rm a}_{\nu+\chi+\tau+\psi+\sigma+\theta+\phi+\xi}-
%{\rm k}_{\nu+\chi+\tau+\psi+\sigma+\theta+\phi+\xi}}
%({\rm n}/p_{\nu+\chi+\tau+\psi+\sigma+\theta+\phi+\xi}
%^{{\rm a}_{\nu+\chi+\tau+\psi+\sigma+\theta+\phi+\xi}})
\\
&& \nonumber \\
 {\rm and}
\nonumber \\
&& \nonumber \\
\mathcal B_{1} & = & \sum_{i=1}^{\nu} p_{i}^{{\rm b}_{i}-{\rm
k}_{i}} ({\rm m}/p_{i}^{{\rm b}_{i}}) +
%p_{1}^{{\rm b}_{1}-{\rm k}_{1}} ({\rm m}/p_{1}^{{\rm b}_{1}}) +
%p_{2}^{{\rm b}_{2}-{\rm k}_{2}} ({\rm m}/p_{2}^{{\rm b}_{2}})
%+...+ p_{\nu}^{{\rm b}_{\nu}-{\rm k}_{\nu}} ({\rm m}/p_{\nu}^{{\rm
%b}_{\nu}})
%+
%\nonumber \\
%&&
\sum_{i=\nu+1}^{\nu+\chi} \rho_{i}p_{i}^{{\rm b}_{i}-{\rm
k}_{i}+1} ({\rm m}/p_{i}^{{\rm b}_{i}})+
%\rho_{\nu+1}p_{\nu+1}^{{\rm b}_{\nu+1}-{\rm k}_{\nu+1}+1} ({\rm
%m}/p_{\nu+1}^{{\rm b}_{\nu+1}})+ \rho_{\nu+2}p_{\nu+2}^{{\rm
%b}_{\nu+2}-{\rm k}_{\nu+2}+2} ({\rm m}/p_{\nu+2}^{{\rm
%b}_{\nu+2}})
%+...+ \nonumber \\
%&&\rho_{\nu+\chi}p_{\nu+\chi}
%^{{\rm b}_{\nu+\chi}-{\rm k}_{\nu+\chi}+1}
%({\rm m}/p_{\nu+\chi}^{{\rm b}_{\nu+\chi}}) +
%\nonumber \\
%&&
\sum_{i=\nu+\chi+1}^{\nu+\chi+\tau} p_{i}^{{\rm b}_{i}-{\rm
k}_{i}} ({\rm m}/p_{i}^{{\rm b}_{i}})+
%p_{\nu+\chi+1}^{{\rm b}_{\nu+\chi+1}-{\rm k}_{\nu+\chi+1}} ({\rm
%m}/p_{\nu+\chi+1}^{{\rm b}_{\nu+\chi+1}})+ p_{\nu+\chi+2}^{{\rm
%b}_{\nu+\chi+2}-{\rm k}_{\nu+\chi+2}} ({\rm
%m}/p_{\nu+\chi+2}^{{\rm b}_{\nu+\chi+2}})
%+...+ \nonumber \\
%&&
%p_{\nu+\chi+\tau}^{{\rm b}_{\nu+\chi+\tau}-{\rm k}_{\nu+\chi+\tau}}
%({\rm m}/p_{\nu+\chi+\tau}^{{\rm b}_{\nu+\chi+\tau}})+
\nonumber \\
&& \sum_{i=\nu+\chi+\tau+1}^{\nu+\chi+\tau+\psi} {\rm j}_{i} ({\rm
m}/p_{i}^{{\rm b}_{i}})+
%{\rm j}_{\nu+\chi+\tau+1} ({\rm m}/p_{\nu+\chi+\tau+1}^{{\rm
%b}_{\nu+\chi+\tau+1}})+ {\rm j}_{\nu+\chi+\tau+2} ({\rm
%m}/p_{\nu+\chi+\tau+2}^{{\rm b}_{\nu+\chi+\tau+2}})
%+...+ \nonumber \\
%&&{\rm j}_{\nu+\chi+\tau+\psi}
%({\rm m}/p_{\nu+\chi+\tau+\psi}^{{\rm b}_{\nu+\chi+\tau+\psi}})+
%\nonumber \\
%&&
\sum_{i=\nu+\chi+\tau+\psi+1}^{\nu+\chi+\tau+\psi+\sigma} p_{i}
^{{\rm b}_{i}-{\rm k}_{i}} ({\rm m}/p_{i}^{{\rm b}_{i}}) +
%p_{\nu+\chi+\tau+\psi+1} ^{{\rm b}_{\nu+\chi+\tau+\psi+1}-{\rm
%k}_{\nu+\chi+\tau+\psi+1}} ({\rm m}/p_{\nu+\chi+\tau+\psi+1}^{{\rm
%b}_{\nu+\chi+\tau+\psi+1}}) +...+
%\nonumber \\
%&&
%p_{\nu+\chi+\tau+\psi+\sigma}
%^{{\rm b}_{\nu+\chi+\tau+\psi+\sigma}-{\rm k}_{\nu+\chi+\tau+\psi+\sigma}}
%({\rm m}/p_{\nu+\chi+\tau+\psi+\sigma}^{{\rm b}_{\nu+\chi+\tau+\psi+\sigma}})+
\nonumber \\
&&
\sum_{i=\nu+\chi+\tau+\psi+\sigma+1}^{\nu+\chi+\tau+\psi+\sigma+\theta}
\rho_{i}p_{i}^ {{\rm b}_{i}- {\rm k}_{i}+1} ({\rm m}/p_{i}^ {{\rm
k}_{i}})+
%\rho_{\nu+\chi+\tau+\psi+\sigma+1}p_{\nu+\chi+\tau+\psi+\sigma+1}^
%{{\rm b}_{\nu+\chi+\tau+\psi+\sigma+1}-
%{\rm k}_{\nu+\chi+\tau+\psi+\sigma+1}+1}
%({\rm m}/p_{\nu+\chi+\tau+\psi+\sigma+1}^
%{{\rm k}_{\nu+\chi+\tau+\psi+\sigma+1}})+...+
%\nonumber \\
%&&
%\rho_{\nu+\chi+\tau+\psi+\sigma+\theta}p_{\nu+\chi+\tau+\psi+\sigma+\theta}^
%{{\rm b}_{\nu+\chi+\tau+\psi+\sigma+\theta}-
%{\rm k}_{\nu+\chi+\tau+\psi+\sigma+\theta}+1}
%({\rm m}/p_{\nu+\chi+\tau+\psi+\sigma+\theta}^
%{{\rm k}_{\nu+\chi+\tau+\psi+\sigma+\theta}})+
%\nonumber \\
%&&
\sum_{i=\nu+\chi+\tau+\psi+\sigma+\theta+1}^{\nu+\chi+\tau+\psi+\sigma+\theta+\phi}
p_{i}^ {{\rm b}_{i}- {\rm k}_{i}} ({\rm m}/p_{i} ^{{\rm b}_{i}})+
%p_{\nu+\chi+\tau+\psi+\sigma+\theta+1}^ {{\rm
%b}_{\nu+\chi+\tau+\psi+\sigma+\theta+1}- {\rm
%k}_{\nu+\chi+\tau+\psi+\sigma+\theta+1}} ({\rm
%m}/p_{\nu+\chi+\tau+\psi+\sigma+\theta+1}
%^{{\rm b}_{\nu+\chi+\tau+\psi+\sigma+\theta+1}})+...+ \nonumber \\
%&&
%\label{luhytoro678olygbol786t}
%p_{\nu+\chi+\tau+\psi+\sigma+\theta+\phi}^
%{{\rm b}_{\nu+\chi+\tau+\psi+\sigma+\theta+\phi}-
%{\rm k}_{\nu+\chi+\tau+\psi+\sigma+\theta+\phi}}
%({\rm m}/p_{\nu+\chi+\tau+\psi+\sigma+\theta+\phi}
%^{{\rm b}_{\nu+\chi+\tau+\psi+\sigma+\theta+\phi}})+
\nonumber \\
&& \label{luhytoro678olygbol786t}
\sum_{i=\nu+\chi+\tau+\psi+\sigma+\theta+\phi+1}^{\nu+\chi+\tau+\psi+\sigma+\theta+\phi+\xi}
{\rm t}_{i} ({\rm m}/p_{i} ^{{\rm b}_{i}})
%{\rm t}_{\nu+\chi+\tau+\psi+\sigma+\theta+\phi+1} ({\rm
%m}/p_{\nu+\chi+\tau+\psi+\sigma+\theta+\phi+1}
%^{{\rm b}_{\nu+\chi+\tau+\psi+\sigma+\theta+\phi+1}})+...+ \nonumber \\
%&&{\rm t}_{\nu+\chi+\tau+\psi+\sigma+\theta+\phi+\xi}
%({\rm m}/p_{\nu+\chi+\tau+\psi+\sigma+\theta+\phi+\xi}
%^{{\rm b}_{\nu+\chi+\tau+\psi+\sigma+\theta+\phi+\xi}}),
\end{eqnarray}
and by
\item
\be
\label{8oly7tg67ri65rf6ui5}
g_{2}=(x^{\mathcal A_{2}},y^{\mathcal B_{2}}),
\ee
where,
\begin{eqnarray}
\label{liyugh6yrfhjtxfktug} \mathcal A_{2}  & = &
\sum_{i=\nu+\chi+\tau+\psi+1}^{\nu+\chi+\tau+\psi+\sigma}
p_{i}^{{\rm a}_{i}- {\rm l}_{i}} ({\rm n}/p_{i}^{{\rm a}_{i}})+
%...+
%\nonumber \\
%&&
%p_{\nu+\chi+\tau+\psi+\sigma}^{{\rm a}_{\nu+\chi+\tau+\psi+\sigma}-
%{\rm l}_{\nu+\chi+\tau+\psi+\sigma}}
%({\rm n}/p_{\nu+\chi+\tau+\psi+\sigma}^{{\rm a}_{\nu+\chi+\tau+\psi+\sigma}})+
%\nonumber \\
%&&
\sum_{i=\nu+\chi+\tau+\psi+\sigma+\theta+1}^{\nu+\chi+\tau+\psi+\sigma+\theta+\phi}
p_{i} ^{{\rm a}_{i} -{\rm l}_{i}} ({\rm n}/p_{i} ^{{\rm a}_{i}})
\\
%p_{\nu+\chi+\tau+\psi+\sigma+\theta+1}
%^{{\rm a}_{\nu+\chi+\tau+\psi+\sigma+\theta+1}
%-{\rm l}_{\nu+\chi+\tau+\psi+\sigma+\theta+1}}
%({\rm n}/p_{\nu+\chi+\tau+\psi+\sigma+\theta+1}
%^{{\rm a}_{\nu+\chi+\tau+\psi+\sigma+\theta+1}})+...+ \nonumber \\
%&&
%\label{liyugh6yrfhjtxfktug}
%p_{\nu+\chi+\tau+\psi+\sigma+\theta+\phi}
%^{{\rm a}_{\nu+\chi+\tau+\psi+\sigma+\theta+\phi}
%-{\rm l}_{\nu+\chi+\tau+\psi+\sigma+\theta+\phi}
%}
%({\rm n}/p_{\nu+\chi+\tau+\psi+\sigma+\theta+\phi}
%^{{\rm a}_{\nu+\chi+\tau+\psi+\sigma+\theta+\phi}} )\\
&&  \nonumber \\
&&{\rm and} \nonumber \\
&& \nonumber \\
\label{kgd6u5rdf6ujtfk}
 \mathcal B_{2}  & = &
\sum_{i=\nu+\chi+\tau+\psi+\sigma+1}^{\nu+\chi+\tau+\psi+\sigma+\theta}
 p_ {i}^{{\rm
b}_{i}- {\rm l}_{i}} ({\rm m}/p_ {i}^{{\rm b}_{i}})+
% p_ {\nu+\chi+\tau+\psi+\sigma+1}^{{\rm
%b}_{\nu+\chi+\tau+\psi+\sigma+1}- {\rm
%l}_{\nu+\chi+\tau+\psi+\sigma+1}} ({\rm m}/p_
%{\nu+\chi+\tau+\psi+\sigma+1}^{{\rm
%b}_{\nu+\chi+\tau+\psi+\sigma+1}})+...+
%\nonumber \\
%&&
%p_
%{\nu+\chi+\tau+\psi+\sigma+\theta}
%{{\rm b}_{\nu+\chi+\tau+\psi+\sigma+\theta}-
%{\rm l}_{\nu+\chi+\tau+\psi+\sigma+\theta}}
%({\rm m}/p_
%{\nu+\chi+\tau+\psi+\sigma+\theta}
%^{{\rm b}_{\nu+\chi+\tau+\psi+\sigma+\theta}})+
%\nonumber \\
%&&
\!\! \! \! \! \! \! \!
\sum_{i=\nu+\chi+\tau+\psi+\sigma+\theta+\phi+1}^{\nu+\chi+\tau+\psi+\sigma+\theta+\phi+\xi}
p_{i} ^{{\rm b}_{i} -{\rm l}_{i}} ({\rm m}/p_{i} ^{{\rm b}_{i}}).
%p_{\nu+\chi+\tau+\psi+\sigma+\theta+\phi+1}
%^{{\rm b}_{\nu+\chi+\tau+\psi+\sigma+\theta+\phi+1}
%-{\rm l}_{\nu+\chi+\tau+\psi+\sigma+\theta+\phi+1}}
%({\rm m}/p_{\nu+\chi+\tau+\psi+\sigma+\theta+\phi+1}
%^{{\rm b}_{\nu+\chi+\tau+\psi+\sigma+\theta+\phi+1}})+...+
%\nonumber \\
%&&
%\label{kgd6u5rdf6ujtfk}
%p_{\nu+\chi+\tau+\psi+\sigma+\theta+\phi+\xi}
%^{{\rm b}_{\nu+\chi+\tau+\psi+\sigma+\theta+\phi+\xi}
%-{\rm l}_{\nu+\chi+\tau+\psi+\sigma+\theta+\phi+\xi}}
%({\rm m}/p_{\nu+\chi+\tau+\psi+\sigma+\theta+\phi+\xi}
%^{{\rm b}_{\nu+\chi+\tau+\psi+\sigma+\theta+\phi+\xi}}).
\end{eqnarray}
\end{enumerate}

The order $\;|\mathcal C|\;$ of the group $\;\mathcal C\;$ is given by
\begin{eqnarray}
\label{lygvukygvuklygoyg}
%\label{8oy7tgygvhggvyut67i5r7i65tf}
|\mathcal C|&=&
\prod_{i=1}^{\nu+\chi+\tau+\psi+\sigma+\theta+\phi+\xi}
p_{i}^{{\rm k}_{i}} \times
\prod_{i=\nu+\chi+\tau+\psi+1}^{\nu+\chi+\tau+\psi+\sigma+\theta+\phi}
p_{i}^{{\rm l}_{i}}. 
%\quad  .
%p_{\nu+\chi+\tau+\psi+1}^{{\rm l}_{\nu+\chi+\tau+\psi+1}}\cdot
%p_{\nu+\chi+\tau+\psi+2}^{{\rm l}_{\nu+\chi+\tau+\psi+2}}\cdot \cdot \cdot
%p_{\nu+\chi+\tau+\psi+\sigma}
%^{{\rm l}_{\nu+\chi+\tau+\psi+\sigma}}\cdot \nonumber \\
%&&
%p_{\nu+\chi+\tau+\psi+\sigma+1}
%^{{\rm l}_{\nu+\chi+\tau+\psi+\sigma+1}}\cdot
%p_{\nu+\chi+\tau+\psi+\sigma+2}
%^{{\rm l}_{\nu+\chi+\tau+\psi+\sigma+2}}\cdot \cdot \cdot
%p_{\nu+\chi+\tau+\psi+\sigma+\theta}
%^{{\rm l}_{\nu+\chi+\tau+\psi+\sigma+\theta}}\cdot \nonumber \\
%&&
%p_{\nu+\chi+\tau+\psi+\sigma+\theta+1}
%^{{\rm l}_{\nu+\chi+\tau+\psi+\sigma+\theta+1}}\cdot
%p_{\nu+\chi+\tau+\psi+\sigma+\theta+2}
%^{{\rm l}_{\nu+\chi+\tau+\psi+\sigma+\theta+2}}\cdot \cdot \cdot
%p_{\nu+\chi+\tau+\psi+\sigma+\theta+\phi}
%^{{\rm l}_{\nu+\chi+\tau+\psi+\sigma+\theta+\phi}}\cdot \nonumber \\
%&&
\label{8oy7tgygvhggvyut67i5r7i65tf}
%p_{\nu+\chi+\tau+\psi+\sigma+\theta+\phi+1}
%^{{\rm l}_{\nu+\chi+\tau+\psi+\sigma+\theta+\phi+1}}\cdot
%p_{\nu+\chi+\tau+\psi+\sigma+\theta+\phi+2}
%^{{\rm l}_{\nu+\chi+\tau+\psi+\sigma+\theta+\phi+2}}\cdot \cdot \cdot
%p_{\nu+\chi+\tau+\psi+\sigma+\theta+\phi+\xi}
%^{{\rm l}_{\nu+\chi+\tau+\psi+\sigma+\theta+\phi+\xi}}.
\end{eqnarray}
The non$-$negative integers
$\; \nu,\chi,\tau,\psi,\sigma,\theta,\phi,\xi \;$ are such that
$\; \nu+\chi+\tau+\psi+\sigma+\theta+\phi+\xi \leq {\rm s}
$.
Moreover,
$\;
(p_{\rm i}^{{\rm a}_{\rm i}},p_{\rm i}^{{\rm b}_{{\rm i}}})=
\mathcal P({\rm p}_{{\rm i}}^{a_{{\rm i}}},{\rm p}_{{\rm i}}
^{\beta_{{\rm i}}}),\;{\rm i}=1,2,...,{\rm s},\;$
for some permutation  $\mathcal P$ of the $s$ pairs of numbers
$\;({\rm p}_{1}^{a_{1}},{\rm p}_{1}^{\beta_{1}}),
({\rm p}_{2}^{a_{2}},{\rm p}_{2}^{\beta_{2}}),...,
({\rm p}_{{\rm s}}^{a_{{\rm s}}},{\rm p}_{{\rm s}}^{\beta_{{\rm s}}})$.
Furthermore, when $\;{\rm i}\in \left \{1,2,...,\nu \right \},\;$ then 
$\;{\rm r}_{\rm i} \in \left \{
0,1,2,...,p_{\rm i}^{{\rm k}_{\rm i}}-1 \right \} \; and \;
%\;,\;
{\rm k}_{\rm i}\leq {\rm min}({\rm a}_{\rm i},{\rm b}_{\rm i}),
\;$
and
when, $\;w\in \left \{\nu+\chi+1,\nu+\chi+2,...,\nu+\chi+\tau \right \},
\;$ then 
$\;{\rm j}_{w} \in \left \{
0,1,2,...,p_{w}^{{\rm a}_{w}}-1 \right \},
\; and \;
{\rm a}_{w}<{\rm k}_{w}\leq{\rm b}_{w}.
\;$
When $\;{\rm q}\in \left \{\nu+1,\nu+2,...,\nu+\chi \right \},\;$ then
$\;{\rho}_{\rm q} ~\in \left \{
0,1,2,...,p_{\rm q}^{{\rm k}_{\rm q}-1}-1 \right \}\;
\;$  \; and \; $\;{\rm k}_{\rm q}\leq {\rm min}({\rm a}
_{{\rm q}},{\rm b}_{{\rm q}}),$ and
when, $\;{y} \in \left \{\nu+\chi+\tau\!+\!1,
\nu+\chi+\tau+2,...,\nu+\chi+\tau+ \psi \right \} \;$ then
$\;{\rm j}_{y} \in \left \{
0,1,2,...,p_{y}^{{\rm b}_{y}}-1 \right \}
\; and \;
{\rm a}_{y} \geq {\rm k}_{y} >
{\rm b}_{y}.\;$
When $\;{\rm i}_{1}\in \left
\{\nu+\chi+\tau+\psi+1,\nu+\chi+\tau+\psi+2,\right. $  
%\newline ...,
$ \left. ...,\nu+\chi+\tau+\psi+\sigma \right \},\;$ then
$\;{\rm r}_{{\rm i}_{1}} \in \left \{
0,1,2,...,p_{{\rm i}_{1}}^{{\rm k}_{{\rm i}_{1}}-{\rm l}_{{\rm i}_{1}}}
-1 \right \}\; and \;
1\leq{\rm l}_{{\rm i}_{1}}\leq{\rm k}_{{\rm i}_{1}}
\leq {\rm min}({\rm a}_{{\rm i}_{1}},{\rm b}_{{\rm i}_{1}}),
\;$
and when,
$\;{\rm q}_{1}\in \left
\{\nu+\chi+\tau+
\right . $
$\left . \psi +\sigma+1,...,
\nu+\chi+\tau+\psi+\sigma+\theta \right \},\;$ then
$\;\rho_{{\rm q}_{1}} \in \left
\{0,1,2,...\right .$ 
$\left . ,p_{{\rm q}_{1}}^{{\rm k}_{{\rm q}_{1}}-{\rm l}_{{\rm q}_{1}}-1}-1
\right \} \; and \; 1 \leq {\rm l}_{{\rm q}_{1}} < {\rm k}_{{\rm q}_{1}} \leq
{\rm min}({\rm a}_{{\rm q}_{1}},{\rm b}_{{\rm q}_{1}})$.
When
$\;w_{1}
\in \left \{\nu+\chi+\tau+
\psi+\sigma+\theta+1,
\nu+\chi+\tau+\psi+\sigma+\theta+2,...,
\nu+ \right . $
$\left. + \chi +
\tau+\psi+\sigma+\theta+\phi \right \},
\;$ then
$\;{\rm t}_{{w}_{1}} \in \left \{
0,1,2,..
.,p_{{w}_{1}}^{{\rm a}_{{w}_{1}}-{\rm l}_{{w}_{1}}}-1 \right \}
\;$ \; and \;
$1\leq{\rm l}_{w_{1}}\leq{\rm a}_{{w}_{1}}<
{\rm k}_{w_{1}}\leq{\rm b}_{{w}_{1}}
.\;$
Finally, when $\;y_{1} \in \left \{\nu+\chi+\tau+\psi+\sigma+\theta
+\phi+1, \nu+\chi+\tau+\psi+\sigma+\theta+\phi+2,...,
\nu+ \right . $
$ \left .\chi+\tau+\psi+\sigma+\theta+\phi+\xi \right \}, \, \, then \,\,
{\rm j}_{{y}_{1}} \in \left \{
0,1,2,...,p_{{y}_{1}}^{{\rm b}_{{y}_{1}}-{\rm l}_{{y}_{1}}}-1
\right \}$ and
$1 \leq {\rm l}_{{y}_{1}} \leq {\rm b}_{{y}_{1}} <
{\rm k}_{{y}_{1}} \leq {\rm a}_{{y}_{1}} $.
\normalfont
\end{prop}

\noindent
We conclude therefore that every subgroup $\; \mathcal C \;$ of 
${\rm C}_{\rm n} \times {\rm C}_{\rm m}$ with two generators
%Therefore we
%conclude that the non-cyclic group $\; \mathcal C \;$
%(\ref{hjgvnmbgkugkbsd}) 
is the direct product of two cyclic groups
$\; \mathcal C_{1} \;$ and $\; \mathcal C_{2}$: \be
\label{yukitrftcgjhfd65tfdcccc} \mathcal C=\mathcal C_{1} \times
\mathcal C_{2}. \ee 
By construction, the order of the non$-$cyclic
group $\;\mathcal C\;$ is given by (\ref{lygvukygvuklygoyg}).
%This completes the proof.
As it was pointed out before
the choice of the cyclic subgroups
$\;\mathcal C_{1}\;$ and $\;\mathcal C_{2}\;$
in expression (\ref{yukitrftcgjhfd65tfdcccc}) is highly non$-$unique.
The choices displayed  in expressions
(\ref{ku6tr967rfujr6xikt5fklo})
and (\ref{8oly7tg67ri65rf6ui5}) are only two specific choices among the
many possible. What is common in all these choices is that the orders
of the cyclic groups $\;\mathcal C_{1}\;$ and
$\;\mathcal C_{2}\;$ are {\it not} relatively prime. In fact, in the
particular choice we made we have
\begin{eqnarray}
\label{lyuitgko6ygbk,u67tgo7g} |\mathcal C_{1}|&=&
\prod_{i=1}^{\nu+\chi+\tau+\psi+\sigma+\theta+\phi+\xi}
p_{i}^{{\rm k}_{i}}
\end{eqnarray}
and,
\begin{eqnarray}
\label{8oy7tgygvhggvyut67i5r7i65tf} |\mathcal C_{2}|&=&
\prod_{i=\nu+\chi+\tau+\psi+1}^{\nu+\chi+\tau+\psi+\sigma+\theta+\phi+\xi}
p_{i}^{{\rm l}_{i}}.
%p_{\nu+\chi+\tau+\psi+1}^{{\rm l}_{\nu+\chi+\tau+\psi+1}}\cdot
%p_{\nu+\chi+\tau+\psi+2}^{{\rm l}_{\nu+\chi+\tau+\psi+2}}\cdot
%\cdot \cdot p_{\nu+\chi+\tau+\psi+\sigma}
%^{{\rm l}_{\nu+\chi+\tau+\psi+\sigma}}\cdot \nonumber \\
%&&
%p_{\nu+\chi+\tau+\psi+\sigma+1}
%^{{\rm l}_{\nu+\chi+\tau+\psi+\sigma+1}}\cdot
%p_{\nu+\chi+\tau+\psi+\sigma+2}
%^{{\rm l}_{\nu+\chi+\tau+\psi+\sigma+2}}\cdot \cdot \cdot
%p_{\nu+\chi+\tau+\psi+\sigma+\theta}
%^{{\rm l}_{\nu+\chi+\tau+\psi+\sigma+\theta}}\cdot \nonumber \\
%&&
%p_{\nu+\chi+\tau+\psi+\sigma+\theta+1}
%^{{\rm l}_{\nu+\chi+\tau+\psi+\sigma+\theta+1}}\cdot
%p_{\nu+\chi+\tau+\psi+\sigma+\theta+2}
%^{{\rm l}_{\nu+\chi+\tau+\psi+\sigma+\theta+2}}\cdot \cdot \cdot
%p_{\nu+\chi+\tau+\psi+\sigma+\theta+\phi}
%^{{\rm l}_{\nu+\chi+\tau+\psi+\sigma+\theta+\phi}}\cdot \nonumber \\
%&&
%\label{8oy7tgygvhggvyut67i5r7i65tf}
%p_{\nu+\chi+\tau+\psi+\sigma+\theta+\phi+1}
%^{{\rm l}_{\nu+\chi+\tau+\psi+\sigma+\theta+\phi+1}}\cdot
%p_{\nu+\chi+\tau+\psi+\sigma+\theta+\phi+2}
%^{{\rm l}_{\nu+\chi+\tau+\psi+\sigma+\theta+\phi+2}}\cdot \cdot \cdot
%p_{\nu+\chi+\tau+\psi+\sigma+\theta+\phi+\xi}
%^{{\rm l}_{\nu+\chi+\tau+\psi+\sigma+\theta+\phi+\xi}}.
\end{eqnarray}
For the purposes of our study it is convenient to rewrite 
%It will prove convenient later on to rewrite
$\;\mathcal C_{1}\;$ and $\;\mathcal C_{2}\;$
as  subgroups of ${\rm S}{\rm O}(2) \times {\rm S}{\rm O}(2)\;$.
This is the content of the following Theorem.
%By collecting the previous results we have the following theorem
%about the finite non-cyclic subgroups of
%$\;{\rm S}{\rm O}(2) \times {\rm S}{\rm O}(2)\;$
\begin{thrm}
\label{uyghkgvyjtrfitfvkyjtggfi7k6g}
Let $\;{\rm n}\;$ and $\;{\rm m}\;$ be {\it any} non$-$negative integers.
Then all the finite non$-$cyclic subgroups $\; \mathcal C \;$  of
$\;{\rm S}{\rm O}(2) \times {\rm S}{\rm O}(2)\;$
can be written as the direct product of two cyclic groups
$\; \mathcal C_{1} \;$ and $\; \mathcal C_{2} \;$
whose orders are not relatively prime. Thus we have
\begin{equation}
\mathcal C=\mathcal{C}_{1} \times \mathcal{C}_{2}.
\end{equation}
The choice of $\; \mathcal C_{1} \;$ and $\; \mathcal C_{2} \;$
is highly non unique.  A possible choice for
$\; \mathcal C_{1} \;$ 
%(expression (\ref{ytgoylglyugol8y7tgolyugbolytg}))
is given by
\begin{eqnarray}
\mathcal{C}_{1} & = &
\left ( R \left ( \left ( \frac {2 \pi}{{\rm n}} {\mathcal A}_{1} \right )
 i
\right ) ,
 R  \left ( \left ( \left .
\frac {2 \pi}{{\rm m}}  {\mathcal B}_{1} \right ) i
\right ) \right . \right ), 
%\nonumber
\end{eqnarray}
and a possible choice for $\; \mathcal C_{2} \;$
%(expressions
%(\ref{;uilgbyukfvjytmgcvjkytfv}) and
%(\ref{lygfkygvhjgyukgflyuig}))
is given by
\begin{eqnarray}
\label{e2}
\mathcal C_{2} & = &
\left ( R \left ( \left ( \frac {2 \pi}{{\rm n}} {\mathcal A}_{2}
\right )
 i
\right ) ,
 R  \left ( \left ( \left .
\frac {2 \pi}{{\rm m}}  {\mathcal B}_{2} \right ) i
\right ) \right . \right ).  
%\nonumber \\
%&& \nonumber \\
%&=&
%\left ( R \left ( \left ( \frac {2 \pi}{{\rm n}'_{2}}
%\right )
% i'
%\right ) ,
% R  \left ( \left ( \left .
%\frac {2 \pi}{{\rm m}'_{2}}   \right ) i'
%\right ) \right . \right ). \nonumber
\end{eqnarray}
The meaning and the ranges of the parameters
$\; \mathcal A_{1}, \ \mathcal B_{1}, \ \mathcal A_{2}, \ \mathcal B_{2} \;$
appearing in these expressions  are displayed in Proposition
(\ref{kittttttygy5ygyk5j55}).
% and the values  of
%$\;{\rm n}'_{2}\;$ and $\;{\rm m}'_{2}\;$ are given by
%(\ref{bkuyrfu6jtrfcj}) and (\ref{luhgjkytdfvjk6ytir})
%correspondingly.
For each specific subgroup these parameters take specific values.
Different values of the parameters
correspond to different subgroups and vice versa.
The order of the group $\; \mathcal C \;$ is given by
(\ref{lygvukygvuklygoyg}).
\end{thrm}
%This completes our consideration of the non-cyclic subgroups of
%$\;{\rm C}_{{\rm n}} \times {\rm C}_{{\rm m}}\;$.
Using now Proposition \ref{kittttttygy5ygyk5j55}, Theorem \ref{uyghkgvyjtrfitfvkyjtggfi7k6g}
and the observation made at the beginning of this subsection
we give now in detail the little groups with two generators.

%The choices displayed  in expressions
%(\ref{ku6tr967rfujr6xikt5fklo})
%and (\ref{8oly7tg67ri65rf6ui5})

\begin{prop}
\label{kloikloikloikloiklo}
%\label{kuyfkuyfvvfyufkyuffku}
 Let ${\rm n}$ and ${\rm m}$ be any
positive even numbers. Then
\begin{equation}
{\rm C}_{\rm n} \times {\rm C}_{\rm m}= ({\rm C}_{2^{a_{1}}}
\times {\rm C}_{2^{\beta_{1}}}) \times ({\rm C}_{{\rm
p}_{2}^{a_{2}}} \times {\rm C}_{{\rm p}_{2}^{\beta_{2}}}) \times
({\rm C}_{{\rm p}_{3}^{a_{3}}} \times {\rm C}_{{\rm
p}_{3}^{\beta_{3}}}) \times \ ... \ \times ({\rm C}_{{\rm p}_{\rm
s}^{a_{\rm s}}} \times {\rm C}_{{\rm p}_{\rm s}^{\beta_{\rm s}}})
\end{equation}
where $\;a_{1} \geq 1, \;$ $\; \beta_{1}  \geq 1, \;$ $\; {\rm
p}_{2},{\rm p}_{3},...,{\rm p}_{\rm s} \; $ are odd primes,
$\;a_{\rm i} \geq 0, \;$ and, $\; \beta_{\rm i} \geq 0,$
${\rm i} \in \left \{ 2,3,...,{\rm s} \right \}. \;$ Every
subgroup $\; \mathcal C \;$ of $\;K=SO(2)\times SO(2)\;$ with two
generators which falls into the class $\;{\rm E}_{1} \;$
%which contains the element $\;(-I,-I)\;$
is written uniquely in the form
$$
\mathcal C={\rm A}_{1} \times {\rm A}_{2} \times {\rm A}_{3}\times
... \times {\rm A}_{\rm s},
$$
where $\;{\rm A}_{{\rm i}}\;$ is a subgroup of $\; {\rm C}_{{\rm
p}_{{\rm i}}^{a_{\rm i}}} \times {\rm C}_{{\rm p}_{{\rm
i}}^{\beta_{\rm i}}}, \;$ $\;{\rm i} \in \left \{ 1,2,...,{\rm s}
\right \}, \; ({\rm p}_{1}=2)$.
$\;{\rm A}_{1}\;$ is not the
identity element and must be cyclic.
%can be either cyclic or have
%two generators. We distinguish two cases
%\begin{enumerate}
%\item{
%When $\;{\rm A}_{1}\;$ is cyclic then at least one of the
%$\;{\rm A}_{{\rm i}}\;$ ,$\;{\rm i} \in \left \{ 2,...,{\rm s} \right \} \;$,
%which are not all necessarily different from the identity element,
%has two generators.
In particular, $\;{\rm A}_{1}\;$ is restricted to be one of the
following cyclic subgroups of $\;{\rm C}_{2^{a_{1}}} \times {\rm
C}_{2^{\beta_{1}}}\;$
$$
\left ( R \left ( \left ( \frac {2 \pi}{2^ {{\rm k}_{1}}} {\rm r}
\right )  i_{1} \right ) ,
 R  \left ( \frac {2 \pi}{2^{{\rm k}_{1}}}    i_{1} \right )
\right ),
$$
where $\;1 \leq {\rm k}_{1} \leq {\rm min}(a_{1},\beta_{1})$,
$\; {\rm r} \;$ parametrises the groups and takes values in the
set $\; \left \{ 1,2,...,2^{{\rm k}_{1}}-1 \right \}- \left \{2,2
\cdot 2,...,(2^{{\rm k}_{1}-1}-1)2 \right \} \;$ and $\;i_{1}\;$
enumerates the elements of each group and takes values in the set
$\; \left \{0,1,2,...,2^{{\rm k}_{1}}-1 \right \}.$ One of the
$\; {\rm A}_{{\rm i}}, \;$ $\; {\rm i} \in \left \{ 2,...,{\rm s}
\right \}, $ which are not all necessarily different from the
identity element, has two generators. Two generators of $\;
\mathcal C \;$ are given by 
(\ref{ku6tr967rfujr6xikt5fklo}) and (\ref{8oly7tg67ri65rf6ui5}),
%(\ref{kuyitgfui65r6u5dfu}) and
%(\ref{lyuigiktdcheyxuj5te4ujtfuj}), 
where, $\; (p_{\rm i}^{{\rm
a}_{\rm i}},p_{\rm i}^{{\rm b}_{{\rm i}}})= \mathcal P({\rm
p}_{{\rm i}}^{a_{{\rm i}}},{\rm p}_{{\rm i}} ^{\beta_{{\rm
i}}}), \;{\rm i} \left \{ 1,2,...,{\rm s} \right \},\;$ for some
permutation  $\mathcal P$ of the $s$ pairs of numbers $\;({\rm
p}_{1}^{a_{1}},{\rm p}_{1}^{\beta_{1}}), ({\rm p}_{2}^{a_{2}}
%\;$
%\newline
%$\;
,{\rm p}_{2}^{\beta_{2}}),..., ({\rm p}_{{\rm s}}^{a_{{\rm
s}}},{\rm p}_{{\rm s}}^{\beta_{{\rm s}}}),\;$ and where, $\; {\rm
n}\;$ and $\;{\rm m} \;$ are positive even numbers and one of the
primes $\;p_{1},p_{2},...,p_{\nu}\;$ is the prime number $\rm 2$. If
say, $\;p_{t}=2,\;$, $\; t \in \left \{ 1,2,...,\nu \right \}, \;$
then $\;{\rm r}_{t} \in \left \{1,2,...,2^{{\rm k}_{t}}-1 \right
\}- \left \{2,2 \cdot 2,...,(2^{{\rm k}_{t}-1}-1)2 \right \}. $
The rest of the indices $\;{\rm r}_{d},\;$ $\;d \in \left \{
1,2,...,\nu \right \}- \left \{ t \right \}, \;$ take values in
the sets $\left \{0,1,2,...,p_{d}^{{\rm k}_{d}}-1 \right \}. $ The
other indices which appear in
(\ref{ku6tr967rfujr6xikt5fklo}) and (\ref{8oly7tg67ri65rf6ui5}) 
%(\ref{kuyitgfui65r6u5dfu}) 
take
values in the sets which are displayed in Proposition
(\ref{kittttttygy5ygyk5j55}). In 
(\ref{ku6tr967rfujr6xikt5fklo}) and (\ref{8oly7tg67ri65rf6ui5}) 
%(\ref{kuyitgfui65r6u5dfu}) 
some
of the exponents $\;{\rm a}_{\rm i} \;$ and $\;{\rm b}_{\rm i}
,\;$ $\;{\rm i} \in \left \{ 1,2,3,...,\nu+\chi+\tau+\psi \right
\} - \left \{ t \right \}, \;$ or in fact all of them, can be equal
to zero. On the other hand, in
(\ref{ku6tr967rfujr6xikt5fklo}) and (\ref{8oly7tg67ri65rf6ui5})  
%(\ref{kuyitgfui65r6u5dfu}) 
at least
one of the products of primes $\;{\rm a}_{\rm i} \cdot {\rm
b}_{\rm i} \neq 0, \; \; {\rm i} \in \left \{
\nu+\chi+\tau+\psi+1,...,
\nu+\chi+\tau+\psi+\sigma+\theta+\phi+\xi \right \} $.
\end{prop}
For the purposes of this study it is convenient to rewrite
%It will prove convenient later on to rewrite
the little groups  with two generators $\mathcal C$
%$\;\mathcal C_{1}\;$ and $\;\mathcal C_{2}\;$
as  subgroups of ${\rm S}{\rm O}(2) \times {\rm S}{\rm O}(2)$.

\begin{thrm}
%The subgroup 
Every little group 
$\; \mathcal C \;$ with two generators is written as a direct product of
two cyclic groups in a highly non$-$unique way. A possible choice
%(expressions (\ref{yukitrftcgjhfd65tfdcccc}),
%(\ref{ytgoylglyugol8y7tgolyugbolytg}), and, (\ref{;uilgbyukfvjytmgcvjkytfv}))
is the following
$$
\mathcal C = \mathcal C_{1} \times \mathcal C_{2},
$$
where,
\be
\label{ex1}
\mathcal{C}_{1}  = \left ( R \left ( \left ( \frac {2 \pi}{{\rm
n}} {\mathcal A}_{1} \right )
 i_{1}
\right ) ,
 R  \left ( \left ( \left .
\frac {2 \pi}{{\rm m}}  {\mathcal B}_{1} \right ) i_{1} \right )
\right . \right ),
\ee
and where,
\be
\label{ex2}
\mathcal{C}_{2}  = \left ( R \left ( \left ( \frac {2 \pi}{{\rm
n}} {\mathcal A}_{2} \right )
 i_{2}
\right ) ,
 R  \left ( \left ( \left .
\frac {2 \pi}{{\rm m}}  {\mathcal B}_{2} \right ) i_{2} \right )
\right . \right ).
\ee
The coefficients $\; \mathcal A_{1} \;$ and $\; \mathcal B_{1} \;$
are given by (\ref{vyitd6uir5er7i65rfiky}) and
(\ref{luhytoro678olygbol786t}) correspondingly, and the
coefficients $\; \mathcal A_{2} \;$ and $\; \mathcal B_{2} \;$ are
given respectively by (\ref{liyugh6yrfhjtxfktug}) and
(\ref{kgd6u5rdf6ujtfk}). The index $\;i_{1}$ enumerates the
elements of the group $\;\mathcal C_{1} \;$ and takes values in
the set $\; \{0,1,...,| \mathcal C_{1} |-1 \}, \;$ where $\; | \mathcal
C_{1} | \;$ is given by (\ref{lyuitgko6ygbk,u67tgo7g}), and
$\;i_{2}\;$ enumerates the elements of the group $\;\mathcal C_{2}
\;$ and takes values in the set $\; \{0,1,...,| \mathcal C_{2} |-1
\} \;$, where $\; | \mathcal C_{2} | \;$ is given by
(\ref{8oy7tgygvhggvyut67i5r7i65tf}). 
The rest of the indices which appear in  (\ref{ex1}) and (\ref{ex2})
are given in Proposition \ref{kloikloikloikloiklo}.

%In the expression (\ref{vyitd6uir5er7i65rfiky}) $\; {\rm
%n}\;$ and $\;{\rm m} \;$ are positive even numbers and one of the
%primes $\;p_{1},p_{2},...,p_{\nu}\;$ is the prime number $\rm 2$.
%If
%say, $\;p_{t}=2,\;$ $\; t \in \left \{ 1,2,...,\nu \right \}, \;$
%then $\;{\rm r}_{t} \in \left \{1,2,...,2^{{\rm k}_{t}}-1 \right
%\}- \left \{2,2 \cdot 2,...,(2^{{\rm k}_{t}-1}-1)2 \right \}. $
%The rest of the indices $\;{\rm r}_{d},\;$ $\;d \in \left \{
%1,2,...,\nu \right \}- \left \{ t \right \}, \;$ take values in
%the sets $\left \{0,1,2,...,p_{d}^{{\rm k}_{d}}-1 \right \}. $
%The
%other indices which appear in (\ref{vyitd6uir5er7i65rfiky}) take
%values in the sets which are displayed in Proposition
%(\ref{kittttttygy5ygyk5j55}). In (\ref{vyitd6uir5er7i65rfiky})
%some of the exponents $\;{\rm a}_{\rm i} \;$ and $\;{\rm b}_{\rm
%i} ,\;$ $\;{\rm i} \in \left \{ 1,2,3,...,\nu+\chi+\tau+\psi
%\right \} - \left \{ t \right \}, $ or in fact all of them, can
%be equal to zero.
%%$\;p_{1},p_{2},...,p_{t-1},p_{t+1},...,p_{\nu},
%%p_{\nu+1},
%%...,p_{\nu+\chi},
%%p_{\nu+\chi+1},...,p_{\nu+\chi+\tau},
%%p_{\nu+\chi+\tau+1},...,
%%p_{\nu+\chi+\tau+\psi}\;$
%%\newline
%%,or in fact all of them,
%%can be equal to one.
%\nopagebreak On the other hand, in (\ref{vyitd6uir5er7i65rfiky})
%at least one of the products of primes \nopagebreak $\;{\rm
%a}_{\rm i} \cdot {\rm b}_{\rm i} \neq 0, \; \; {\rm i} \in \left \{
%\nu+\chi+\tau+\psi+1,...,
%\nu+\chi+\tau+\psi+\sigma+\theta+\phi+\xi \right \} $.
\end{thrm}

%\end{enumerate}
%njhuhyghgtfgfr

\section{Form of the induced representations}
\label{s6}
\label{xasertyudololkij} \indent
%the induced representations of $\; \mathcal H \mathcal B_{c}= A(
%mathcal N )   \bigcirc\!\!\;\!\!\!\;\!\!\!\!s \  _{T} \mathcal G
%$ , then , it is enough to provide the information cited in 1
%and 2 for each of the orbit types.
%It follows from the discussion in chapter \ref{chap2}, section
%\ref {l,yyf,,yjhjjhv} that 
In order to give explicitly the
operators of the representations of $\; \mathcal H \mathcal B
%_{c}
\;$ induced from finite little groups it is necessary to give the
following information \cite{Wigner,Mackey,Mackey1,Simms,Isham,Mackey2}:
\begin{enumerate}
\item{ An irreducible unitary representation $\;U\;$ of
$\;L(\zeta)\;$ on a Hilbert space $\;D\;$ for each
$\;L(\zeta).\;$} \item{ A $\; \mathcal G\;$$-$quasi$-$invariant
measure $\; \mu$, $\; \mathcal G=G \times G$,
 on each orbit $\; \mathcal G \zeta \approx
\mathcal G / L(\zeta)$; where $\; L(\zeta) \;$ denotes the
little group of the base point $\; \zeta \in \mathcal H({\rm
T}^{2}) \;$ of the orbit $\; \mathcal G \zeta$. }
\end{enumerate}
\noindent The actual little groups have been given in Theorem
\ref{nhmjknmnhbgvfcd}. The finite ones are those which fall into
the class $\; {\rm E}_{1}$. These are described in detail in
Propositions \ref{yhyhyuhvbhygygku} and \ref{kloikloikloikloiklo}.

 The information cited in 1 and 2 for each of the
aforementioned groups and the corresponding orbit types is now
provided.

\vspace{0.5cm} 1.
%\hspace{0.5cm} 01.
\hspace{0.2cm} The
 finite little groups
 are either cyclic or can be expressed as the direct product of
 two cyclic groups. The little groups which are cyclic are
described in detail in Proposition \ref{yhyhyuhvbhygygku},
whereas, the little groups which are the direct product of two
cyclic groups are given in detail in Proposition
\ref{kloikloikloikloiklo}. The irreducible unitary representations
$\;U_{N} \;$ of the cyclic groups $\; C_{N} \;$ are
%given in the class of little groups described in case 02. They are
indexed by an integer $\; \nu \;$ which, for distinct
representations, takes values in the set $\; \nu \in
\{0,1,2,...,N-1 \}$. The number of the representations equals to
the order of the group $\;N$. Denoting them by $\; D^{(\nu)}$, 
they are given by multiplication in one complex dimension $\;
D \approx C \;$ by \be
 \label{;ujujuunlwjq;jnqj} D^{(\nu)}\left (
\left (  R \left ( \left ( \frac {2 \pi}{{\rm n}} \mathcal A
\right ) {\rm  j} \right ) ,
 R  \left ( \left ( \frac {2 \pi}{{\rm m}} \mathcal B \right )   {\rm j} \right )
\right )\right )=e^{i\frac{2 \pi}{N}\nu{\rm j}}, \ee where, taking
into account the restrictions given in Proposition
\ref{yhyhyuhvbhygygku}, $\; \mathcal A \;$ and $\; \mathcal B \;$
are given respectively by (\ref{bycoxenasha}) and by
(\ref{bnkyuffdcikp;jp;0i9}).
%are restricted as it is explained in Proposition
%\ref{yhyhyuhvbhygygku} .
The order $\; N \;$ $\; (\equiv \mathcal C )\;$ of the group $\;
C_{N} \;$ is given by (\ref{lyuidfytde5rikyto9}).
%\label{xasertyudololkij}
%is given by (\ref{bycoxenasha}), $\; \mathcal B \;$
%is given by (\ref{bnkyuffdcikp;jp;0i9}) and the order of the group
%$\; N \;$ is given by (\ref{lyuidfytde5rikyto9}).
%%%%%%%%%%%%%%%%%%%%%%%%%%%%%%%%%%%%%%%%%%%%%%%%%%%%%%%%%%%%%%%%%%%%%%%%%%%%%%%%%%%%%%%%%%%%%%%%%%%%%%%%%%%%%%%%%%%%%%%%%%

Let $\; C_{N_{1}} \times C_{N_{2}} \;$ be one of the little groups
which can be expressed as the direct product of two cyclic groups
$\; C_{N_{1}} \;$ and $\; C_{N_{2}}$. The unitary irreducible
representations of  $\; C_{N_{1}} \times C_{N_{2}} \;$ are indexed
by two integers $\; \nu_{1} \;$ and $\; \nu_{2} \;$ which for
distinct representations take independently values in the sets $\;
\{0,1,2,...,N_{1}-1 \}\;$ and $\; \{0,1,2,...,N_{2}-1 \}\;$ 
(this is justified in the remark which follows).
% for
%a justification of this statement look at section
%\ref{jknmbhgsedr} , Remark 2 ). 
Denoting these representations by
$\; D^{(\nu_{1},\nu_{2})}$, they are given by multiplication in
one complex dimension $\; D \approx C \;$ by \begin{eqnarray}
D^{(\nu_{1},\nu_{2})} \left ( \left ( R \left ( \left ( \frac {2
\pi}{{\rm n}} \mathcal A_{1} \right ) {\rm j}_{1} \right ) ,
 R  \left ( \left ( \frac {2 \pi}{{\rm m}} \mathcal B_{1} \right )   {\rm j}_{1} \right )
\right ) \right .  & \times & \nonumber \\ \left . \left ( R \left
( \left ( \frac {2 \pi}{{\rm n}} \mathcal A_{2} \right ) {\rm
j}_{2} \right ) ,
 R  \left ( \left ( \frac {2 \pi}{{\rm m}} \mathcal B_{2} \right )   {\rm j}_{2} \right )
\right ) \right )  & = & \nonumber \\ \label{ir} e^{i\frac{2
\pi}{N_{1}}\nu_{1}{\rm j}_{1}} e^{i\frac{2 \pi}{N_{2}}\nu_{2}{\rm
j}_{2}}, && \end{eqnarray} where, taking into account the
restrictions given in Proposition \ref{kloikloikloikloiklo}, $\;
\mathcal A_{1}$,
%is given by (\ref{vyitd6uir5er7i65rfiky}),
$\; \mathcal B_{1}$,
%is given by (\ref{luhytoro678olygbol786t}),
$\; \mathcal A_{2}$,
%is given by (\ref{liyugh6yrfhjtxfktug}),
and $\; \mathcal B_{2}$, are
%given by (\ref{kgd6u5rdf6ujtfk}) are
given respectively by (\ref{vyitd6uir5er7i65rfiky}),
(\ref{luhytoro678olygbol786t}), (\ref{liyugh6yrfhjtxfktug}), and
(\ref{kgd6u5rdf6ujtfk}).
 The order $\; N_{1} \;$$\; (\equiv \mathcal C_{1} ) \;$ of the
group $\;   \left ( R \left ( \left ( \frac {2 \pi}{{\rm n}}
\mathcal A_{1} \right ) {\rm j}_{1} \right ) ,
 R  \left ( \left ( \frac {2 \pi}{{\rm m}} \mathcal B_{1} \right )   {\rm j}_{1} \right )
\right ) \;$ is given by (\ref{lyuitgko6ygbk,u67tgo7g}) and the
order $\; N_{2} \;$$\;(\equiv \mathcal C_{2}) \;$  of the group
$\;
 \left ( R \left
( \left ( \frac {2 \pi}{{\rm n}} \mathcal A_{2} \right ) {\rm
j}_{2} \right ) ,
 R  \left ( \left ( \frac {2 \pi}{{\rm m}} \mathcal B_{2} \right )   {\rm j}_{2} \right )
\right ) \;$ is given by (\ref{8oy7tgygvhggvyut67i5r7i65tf}).

%\vspace{0.5cm}
\noindent
A remark now is in order regarding the unitary irreducible
representations of  $\; C_{N_{1}} \times C_{N_{2}}$. 
The 
%typical 
problem we encounter in 
this case
%all cases  01 to 05 
is
the determination of the unitary irreducible representations of
the direct product $\; A \times B$, where $\; A \;$ and $\; B
\;$ are abelian groups, and where the unitary irreducible
representations of $\; A \;$ and $\; B \; $are known. The group
$\; A \times B \;$ is abelian and therefore its irreducible
representations are one-dimensional. Moreover, since by assumption
they are complex, they operate in one complex dimension $\; D
\approx C$. Let $\; U(\xi,\omega) \;$ be a unitary irreducible
complex representation of the group $\; A \times B$; the
parameters $\; \xi \;$ and $\; \omega \;$ enumerate the elements
of the groups $\;A\;$ and $\;B\;$ correspondingly, and they  are
either continuous or discrete depending on the groups $\;A\;$ and
$\;B$, them being either continuous or finite. Then $\;
U(\xi,\omega) \;$ will have the form \be  U(\xi,\omega)
=u(\xi,\omega) \emph{I}, \ee where $\; \emph{I} \;$ denotes the
identity operator in one complex dimension. Let $\; \xi \;$ take a
specific value $\; \xi=\xi_{{\rm o}}$. Then the complex number
of modulus one $\; u(\xi, \omega)\;$ takes the form \be
\label{lylyhglhglhglih} \; u(\xi_{{\rm o}},
\omega)=\alpha(\xi_{{\rm o}})u_{B}(\omega), \;\ee where $\;
\alpha(\xi_{{\rm o}}) \;$ is a complex number of modulus one which
is a function of the specific choice $\; \xi=\xi_{{\rm o}} \;$ we
made, and $\; u_{B}(\omega) \;$ is a representation
(unitary, irreducible) of the group $\; B$. We repeat the same
argument with $\; \omega \;$ now. So, let $\; \omega \;$ take a
specific value $\; \omega=\omega_{{\rm o}}$. Then $\;
u(\xi,\omega) \;$ equals to \be \label{uughuglygoygy}
\;u(\xi,\omega_{{\rm o}})=u_{A}(\xi)\beta(\omega_{{\rm o}}),\; \ee
where, $\; \beta(\omega_{{\rm o}}) \;$ is a complex number of
modulus one which is a function of the specific choice $\;
\omega=\omega_{{\rm o}} \;$ we made, and $\; u_{A}(\xi)\;$ is a
representation (unitary, irreducible) of the group $\; A$. When
$\; \xi=\xi_{{\rm o}} \; $ and $\; \omega=\omega_{{\rm o}} \;$ the
right hand side of  Eqs. (\ref{lylyhglhglhglih}) and
(\ref{uughuglygoygy}) are identical and therefore we obtain \be
u(\xi_{{\rm o}},\omega_{{\rm o}})=u_{A}(\xi_{{\rm
o}})u_{B}(\omega_{{\rm o}}). \; \ee Consequently, {\it all} the
unitary irreducibles $\; u(\xi,\omega) \;$ of the direct product
$\; A \times B \;$ have the form \be \label{ygvuyhytfgytytfitf}
u(\xi,\omega)=u_{A}(\xi)u_{B}(\omega). \ee It is equation
(\ref{ygvuyhytfgytytfitf}) which was used in 
%cases 01-05 in order
(\ref{ir})
to determine  the unitary irreducibles of the non$-$cyclic little groups.
%direct product $\; A \times
%B. \;$

\vspace{0.5cm}
\noindent
%\vspace{0.5cm} 
2.
%\hspace{0.5cm} 01.
\hspace{0.2cm} We now proceed to give the information cited in 2.
Although a $\; \mathcal G$$-$quasi$-$invariant
measure  is all what is needed, a $\; \mathcal G$$-$invariant
measure  will be provided in all cases.

\vspace{0.3cm} 01. \hspace{0.2cm}
 The orbit 01 is homeomorphic to
$\; 01 \approx  ({\rm S}{\rm L}(2,R) \times {\rm S}{\rm L}(2,R))/
\mathcal C$, where the group $\; \mathcal  C \;$ is either
cyclic or is the direct product of two cyclic groups. The coset
space$\; ({\rm S}{\rm L}(2,R) \times {\rm S}{\rm L}(2,R))/
\mathcal C \;$ is the space of orbits of the
%natural
right action $\; {\rm R}_{\mathcal C} \;$
\begin{eqnarray} {\rm R}_{\mathcal C} & : & {\rm
S}{\rm L}(2,R) \times {\rm S}{\rm L}(2,R)\longrightarrow {\rm
S}{\rm L}(2,R) \times {\rm S}{\rm
L}(2,R) \nonumber \\
\label{ugliygt6tvtyfrfttfvtfu}
 {\rm R}_{\mathcal C}((g,h)) & := &
(g,h)c
 \end{eqnarray}
 of $\; \mathcal C \;$ on $\;{\rm S}{\rm L}(2,R)
 \times {\rm S}{\rm L}(2,R)$,
  where $\;
(g,h) \in {\rm S}{\rm L}(2,R) \times {\rm S}{\rm L}(2,R) \;$ and
$\; c \in \mathcal C$. Since the group $\; \mathcal C \;$ is
finite and since the action (\ref{ugliygt6tvtyfrfttfvtfu}) is
fixed point free the coset space $\; ({\rm S}{\rm L}(2,R) \times
{\rm S}{\rm L}(2,R))/ \mathcal C \;$ inherits the measure of $\;
{\rm S}{\rm L}(2,R) \times {\rm S}{\rm L}(2,R)$.

\vspace{0.6cm} \noindent This completes the necessary information
in order to construct  representations of $\; \mathcal H \mathcal
B
%_{c} 
\;$ induced  from finite little groups. The two remarks made
at the Discussion of \cite{Mel1}
%end of section \ref{l,yyf,,yjhjjhv} 
are also  relevant
here.

\section{Discussion}

\label{s8}

%\noindent Two remarks are in order regarding the
%representations of $\; \mathcal H \mathcal B
%%_{c}
%\;$
%obtained by
%the above construction
%\begin{enumerate}

%\item{
\noindent
By using Propositions \ref{yhyhyuhvbhygygku} and
\ref{kloikloikloikloiklo}
%the results of this section
one could try to give more concrete results than those of
Proposition \ref{nncnncncnmjdk}. The hope is that   ``nice
looking'' connected elementary regions, as opposed to those
described in Proposition \ref{nncnncncnmjdk}, can be found for
the finite actual little groups. It turns out that this is not
such an easy task. To illustrate the difficulty involved, we
firstly consider the cyclic finite actual little groups. If $\;
\mathcal C \;$  is one of them, then, $ \; \mathcal C =\mathcal
C_{N_{1}} \times \mathcal C_{N_{2}}, \;$  where, according to
Proposition \ref{yhyhyuhvbhygygku}, the relatively prime numbers
$ \; N_{1} \; $ and  $\; N_{2} \;$ can be chosen as follows
$$
 N_{1} = \prod_{i=1}^{\nu}p_{i}^{{\rm k}_{i}} \times \prod_{i=\nu + \chi + 1}^{\nu + \chi + \tau}p_{i}^{{\rm k}_{i}}
 \qquad {\rm and} \qquad N_{2} = \prod_{i=\nu+1}^{\nu+\chi}p_{i}^{{\rm k}_{i}} \times
\prod_{i=\nu + \chi + \tau + 1}^{\nu + \chi + \tau +
\psi}p_{i}^{{\rm k}_{i}}, \quad
$$
where one of the primes $\; p_{1},p_{2},...,p_{\nu} \;$ is the
number 2 and the ranges of the exponents $\; {\rm k}_{i} \;$ are
given in Proposition \ref{yhyhyuhvbhygygku}. One can prove that
an elementary region for the group $ \; \mathcal C_{N_{1}} \;$ is
given by $ \; {\rm E}_{N_{1}} = \{ (\rho , \sigma) \in P_{1}(R)
\times P_{1}(R) \ | \ 0 \leq \rho < 2 \pi, \  \ 0 \leq \sigma < 2
\pi / ( N_{1}/2) \} \;$ and an elementary region for the group $
\; \mathcal C_{N_{2}} \;$ is given by $ \; {\rm E}_{N_{2}} = \{
(\rho , \sigma) \in P_{1}(R) \times P_{1}(R) \ | \ 0 \leq \rho < 2
\pi /  N_{2}, \  \ 0 \leq \sigma < 2 \pi  \}. \;  $ One is tempted
to conjecture that an elementary region for the cyclic group $ \;
\mathcal C \;$ is given by the intersection $ \; {\rm E}_{N_{1}}
\cap {\rm E}_{N_{2}}. \;$ It turns out that this is wrong. In
specific cases one can easily  construct counterexamples where $
\; {\rm E}_{N_{1}} \cap {\rm E}_{N_{2}} \;$ is not an elementary
domain for the $\; \mathcal C$$-$action. To illustrate the subtlety
of the problem a few more remarks are here in order. The
definition of elementary domain which is given in section
\ref{uyfvfvjytddxhtesxe} is equivalent to the following one:  (Def 2)
Let $M$ be any topological space, and let $G$ be any
finite group which acts on $M$ from the right.  An \it{elementary
domain} \normalfont for the given action is an open subset ${\rm
E}\subset M$ such that {\it every } $G$$-$orbit intersects $\; {\rm
E} \;$ at only one point. In turn this last definition is
equivalent to the following one: (Def 3) Let $M$ be any
topological space, and let $G$ be any finite group which acts on
$M$ from the right. \ An \it{elementary domain} \normalfont for
the given action is an open subset ${\rm E}\subset M$ which
satisfies: $ \; \forall x \in {\rm E}, \qquad xg \in {\rm E}
\Rightarrow g={\rm I}, \;$ where $ \; {\rm I} \;$ denotes the
identity element of the group $\; G. \;$ In our problem, it can
be proved that the area $ \; \mathcal E = {\rm E}_{N_{1}} \cap
{\rm E}_{N_{2}} \;$ satisfies the following \be
\label{byjukilokjuioljh} \mathcal E g = \mathcal E \Rightarrow g =
{\rm I}. \quad \ee Now, one can show  that Def 3 implies Eq.
(\ref{byjukilokjuioljh}). But Eq. (\ref{byjukilokjuioljh}) does
not imply Def 3, and in fact, as we have already said the region
$ \; \mathcal E = {\rm E}_{N_{1}} \cap {\rm E}_{N_{2}} \;$ is {\it
not} an elementary region for the $\; \mathcal C $$-$action on $
\mathcal P.$
%\; P_{1}(R) \times P_{1}(R) \;$.
Similar problems are encountered
when one tries to find elementary regions for the finite actual
little groups with two generators.

%}

%\end{enumerate}

\bibliographystyle{plain}
{}

\end{document}